\documentclass[11pt,letterpaper]{amsart}
\usepackage{amssymb,amsmath,amsthm}
\usepackage{mathrsfs}
\usepackage{srcltx}
\usepackage[english]{babel}
\usepackage{graphicx}
\usepackage{color}
\usepackage{enumitem}
\usepackage{mathtools}
\usepackage{xcolor}
\usepackage{amsmath,amsthm,amssymb,amsfonts,mathtools,csquotes}
\usepackage{tikz-cd}
\usepackage[all]{xy}
\usepackage[backend=biber,style=alphabetic,maxalphanames=5,sorting=nyt, maxbibnames=99]{biblatex}
\usepackage{comment}
\date{\today}
\usepackage[colorlinks,citecolor=blue,urlcolor=,bookmarks=false,hypertexnames=true]{hyperref}

\newtheorem{theorem}{Theorem}[section]
\newtheorem{corollary}[theorem]{Corollary}
\newtheorem{proposition}[theorem]{Proposition}
\newtheorem{lemma}[theorem]{Lemma}

\theoremstyle{definition}
\newtheorem{assumption}[theorem]{Assumption}
\newtheorem{remark}[theorem]{Remark}
\newtheorem{example}[theorem]{Example}
\newtheorem{definition}[theorem]{Definition}
\newtheorem{notation}[theorem]{Notation}

\newtheorem{Maximal negativity of transverse Lyapunov exponent when c(G)=2 cases}{Case}
\newtheorem{Strict Maximal Inequality of transverse Lyapunov exponent cases}{Case}
\newtheorem{case negativity of transversal lyapunov exponent on positive weyl chambers}{Case}
\newtheorem{case no transverse manifolds in bad directions of generic points}{Case}
\newtheorem{case empirical measures sequence}{Case}
\newtheorem{step main}{Step}
\newtheorem{step bundle from very nice transversal}{Step}
\newtheorem{step no resonance of lo when strict fails}{Step}
\newtheorem{intersection all orbits wT when rG=2 case}{Case}
\newtheorem{refined error recurrence case}{Case}
\newtheorem{Q2 is Q1 and Q0 case}{Case}
\newtheorem{transverse pesin intersects all G-orbits when rG=2 case}{Case}
\newtheorem{coherence proposition case}{Case}

\newenvironment{geometric measure rigidity}{\emph{Proof of Proposition \ref{geometric measure rigidity}:}}{\hfill$\square$}
\newenvironment{Maximal negativity of transverse Lyapunov exponent when c(G)=2}{\emph{Proof of Proposition \ref{Maximal negativity of transverse Lyapunov exponent}:}}{\hfill$\square$}
\newenvironment{Strict Maximal Inequality of transverse Lyapunov exponent for main theorem}{\emph{Proof of Proposition \ref{Strict Maximal Inequality of transverse Lyapunov exponent}:}}{\hfill$\square$}
\newenvironment{proof of P-equivariance case c(G)>2}{\emph{Proof of Proposition \ref{P-equivariance of transverse Pesin manifolds}:}}{\hfill$\square$}
\newenvironment{proof of entropic measure rigidity}{\emph{Proof of Proposition \ref{entropic measure rigidity}:}}{\hfill$\square$}
\newenvironment{proof of coherence proposition}{\emph{Proof of Proposition \ref{coherence proposition}:}}{\hfill$\square$}
\newenvironment{proof of coherence proposition when rG=2} {\emph{Proof of Proposition \ref{coherence proposition} when $c(G)=2$:}}{\hfill$\square$}
\newenvironment{proof of Pertubations by generic elements of P have LU factorization}{\emph{Proof of Proposition \ref{Pertubations by generic elements of $P$ have $LU$ factorization}:}}{\hfill$\square$}
\newenvironment{proof of N intersects all orbits}{\emph{Proof of Proposition \ref{N intersects all orbits}:}}{\hfill$\square$}
\newenvironment{proof of main}{\emph{Proof of Theorem \ref{main}:}}{\hfill$\square$}
\newenvironment{proof of entropic measure rigidity negative transverse exponents}{\emph{Proof of Proposition \ref{entropic measure rigidity} when all elements in $\Lambda_\mu^T(s)$ are negative:}}{\hfill$\square$}
\newenvironment{sketch of proof}{\emph{Sketch of Proof}:}{\hfill$\square$}
\newenvironment{proof unique ergodicity}{\emph{Proof Proposition of \ref{unique P ergodicity}:}}{\hfill$\square$}
\newenvironment{proof of main non resonance}{\emph{Proof of Theorem \ref{main} when $\lop$ is not resonant to a root and $c(G)\geq 3$:}}{\hfill$\square$}
\newenvironment{proof of geometric high entropy method}{\emph{Proof of Theorem \ref{geometric high entropy method}:}}{\hfill$\square$}
\newenvironment{proof of main when c(G)>2 and resonance}{\emph{Proof of Theorem \ref{main} when $\lop$ is resonant to a root:}}{\hfill$\square$}
\newenvironment{Proof of Coherence of Transverse Pesin manifolds with respect to Q0 when delta >2}{\emph{Proof of Proposition \ref{Coherence of Transverse Pesin manifolds with respect to Q0 when delta >2}:}}{\hfill$\square$}
\newenvironment{proof of main mixing}{\emph{Proof of Theorem \ref{main mixing}:}}{\hfill$\square$}
\newenvironment{Proof of Theorem main when c(G)geq 3 and lop is resonant to a root:}{\emph{Proof of Theorem \ref{main} when $\lop$ is resonant to a root}:}{\hfill$\square$}
\newenvironment{Proof of Theorem main mixing when strict maximal inequality:}{\emph{Proof of Theorem \ref{main mixing} when inequality of Proposition \ref{Strict Maximal Inequality of transverse Lyapunov exponent} holds}:}{\hfill$\square$}
\newenvironment{maximal entropy along transverse Pesin}{\emph{Proof of Proposition \ref{maximal entropy along transverse Pesin}} :}{\hfill$\square$}
\newenvironment{coherence theta3}{\emph{Proof of Lemma \ref{coherence theta3}}: }{\hfill$\square$}
\newenvironment{proof almost continuity}{\emph{Proof of Proposition \ref{almost continuity birkhoff averages}}: }{\hfill$\square$}
\newenvironment{proof of equivariant factor map for main mixing when rG>2}{\emph{Proof of Proposition \ref{equivariant factor map for main mixing when rG>2}}: }{\hfill$\square$}
\newenvironment{proof of coherence main mixing rG>2}{\emph{Proof of Proposition \ref{coherence main mixing rG>2}}: }{\hfill$\square$}
\newenvironment{proof of coherence main mixing rG=2}{\emph{Proof of Proposition \ref{coherence mixing when rG=2}}: }{\hfill$\square$}
\newenvironment{proof  main}{\emph{Proof of Theorem \ref{main}}: }{\hfill$\square$}
\newenvironment{proof of main mixing when rG>2}{\emph{Proof of Theorem \ref{main mixing} when $c(G)\geq 3$}:}{\hfill$\square$}
\newenvironment{proof of main mixing when rG=2}{\emph{Proof of Theorem \ref{main mixing} when $c(G)=2$}:}{\hfill$\square$}
\newenvironment{proof of resonance and rG=2 case}{\emph{Proof of Proposition \ref{resonance and rG=2 case}}:}{\hfill$\square$}
\newenvironment{proof of main when c(G)>2}{\emph{Proof of Theorem \ref{main} when $G$ is not locally isomorphic to $\SL(3,\R)$:}}{\hfill$\square$}
\newenvironment{Maximal negativity of transverse Lyapunov exponent when r(G)=2}{\emph{Proof of Proposition \ref{Maximal negativity of transverse Lyapunov exponent}:}}{\hfill$\square$}

\newcommand{\N}{\mathbb{N}}
\newcommand{\Z}{\mathbb{Z}}

\newcommand{\R}{\mathbb R}
\newcommand{\C}{\mathbb C}

\newcommand{\s}{\vspace{0.2cm}}

\newcommand{\supp}{\operatorname{supp}}

\newcommand{\lo}{\lambda_\mu^T}
\newcommand{\haar}{\operatorname{Haar}}
\newcommand{\SL}{\operatorname{SL}}

\newcommand{\locstable}{\w^s_{\operatorname{loc}}}
\newcommand{\lie}{\operatorname{Lie}}
\newcommand{\lop}{\lambda_{\mu^P}^T}
\newcommand{\w}{\mathcal W}
\newcommand{\Pop}{P^{\operatorname{op}}}

\def\acts{\curvearrowright}

\makeatletter
\@namedef{subjclassname@2020}{%
  \textup{2020} Mathematics Subject Classification}
\makeatother

  \DeclareFontFamily{U}{wncy}{}
    \DeclareFontShape{U}{wncy}{m}{n}{<->wncyr10}{}
    \DeclareSymbolFont{mcy}{U}{wncy}{m}{n}
    \DeclareMathSymbol{\Sha}{\mathord}{mcy}{"58}

\title{Classification of Higher Rank Codimension One Actions}
\author{Camilo Arosemena Serrato}
\address{Department of Mathematics\\University of Michigan\\Ann Arbor,  MI 48109, USA}
\email{juancars@umich.edu}
\begin{document}
\begin{abstract}
    Let $G$ be a simple, connected, non-compact, higher rank Lie group, with finite center(e.g., $\SL(n,\R), n\geq 3$). We prove that any smooth and locally free action of $G$ on a closed and connected manifold $M$ of dimension $\dim(G)+1$ is $G$-equivariantly and finitely covered by $G/\Gamma \times S^1$ for some cocompact lattice $\Gamma \le G$, where $G$ acts trivially on $S^1$. This is proved by showing all such actions have a $G$-invariant probability measure combined with work of Stuck and Zimmer.  
    
    This gives a complete classification of all such actions. This result is in the spirit of the Zimmer program.
\end{abstract}
\maketitle
\section{Introduction}

In his 1986 ICM address Robert Zimmer, see \cite{zimmericm}, laid out a series of questions, now called the Zimmer program,  asking whether actions of semisimple, higher rank Lie groups(e.g., $\SL(n,\R)$ for $n\geq 3$) and their lattices(e.g. $\SL(n,\Z)$ in $\SL(n,\R)$) on closed manifolds can be classified, under dynamical and/or geometric assumptions, more specifically, he asked whether all such actions are of algebraic nature, or essentially trivial when the dimension of the manifold is  low with respect to the Lie group. This program has had many important successes, notably, recent work of Brown, Fisher and Hurtado, \cite{bfh}, \cite{bfhinventiones} and \cite{BFHgafa}, on the resolution of various cases of Zimmer's conjecture proving finiteness of smooth actions of higher rank lattices on closed manifolds, when the dimension of the manifold is small compared to the simple Lie group;  see \cite{MR2807830} for results on the Zimmer program up to 2008. Another important recent work in the Zimmer program regarding continuous actions is the work of Deroin and Hurtado, see \cite{dh}, showing that continuous actions of higher rank lattices on the real line by homeomorphisms are finite, i.e. that such lattices are not left-orderable.

Within the spirit of the Zimmer program,  we give a complete classification of smooth actions of higher rank simple Lie groups on closed manifolds that are locally free, i.e., whose stabilizers are discrete,  having codimension one orbits.

\begin{theorem}\label{main smooth}
    Let $G$ be a non-compact, connected, simple higher rank Lie group, with finite center. Let $G\acts M$ be a smooth, locally free action on a closed connected manifold $M$ with $\dim(M)=\dim(G)+1$. Then $M$ is $G$-equivariantly and finitely covered by $G/\Gamma \times S^1$, for some cocompact lattice $\Gamma \le G$, where $G$ acts trivially on $S^1$. 
\end{theorem}

A discrete subgroup $\Gamma\leq G$ is a \emph{cocompact lattice} if $G/\Gamma$ is compact. Given a locally free action of a group $G$ on a manifold, we say it is \emph{of codimension k}, if its orbits are immersed submanifolds of codimension $k$.

Theorem \ref{main smooth} will follow from the main result of our paper, Theorem \ref{main}, by work of Stuck and Zimmer, \cite[Theorem~5.1]{sz}.

\begin{theorem}\label{main}
    Let $G$ be a connected, non-compact, simple higher rank Lie group, with finite center. Let $G\acts M$ be a locally free and smooth action, where $M$ is a closed and connected manifold of dimension $\dim(G)+1$. Then there exists a $G$-invariant probability measure on $M$.
\end{theorem}
The proof of Theorem \ref{main} is the main work of this document. Theorem \ref{main} and Theorem \ref{main smooth} both hold for $C^2$ actions. 

Theorem \ref{main} improves a result from an earlier preprint of the author, see \cite{serrato2025globalrigiditycodimensionactions}. In this document we weaken the hypothesis, by not assuming the existence of a $P$-invariant and $A$-mixing probability measure; and the conclusion is stronger, namely, we prove that there are always $G$-invariant probability measures.

Theorem \ref{main smooth} answers in the negative a question posed by Bertrand Deroin during the 2022
American Institute of Mathematics workshop \emph{Global rigidity of actions by
higher-rank groups}~\cite{AIM2022}, asking whether there could be smooth locally free actions of $\SL(3,\R)$ on a closed manifold of dimension $9$ with a dense orbit. 

Theorem \ref{main smooth} proves that actions satisfying its hypothesis are induced from a finite action of a lattice on $S^1$, see \cite[Section~2.3]{MR2807830} for the notion of induced action. Without the local freeness assumption, actions of such $G$'s might not be induced from the action of a closed subgroup, by an example of Stuck \cite[Section~4]{stuckminimal}.

Theorem \ref{main smooth} does not hold for simple Lie groups of real rank one without Property (T), see the counterexamples in Section \ref{rank one counterexamples}.  

Locally free actions play an important role in the recent breakthroughs on Zimmer's conjecture by Brown, Fisher and Hurtado, see \cite{bfh}, \cite{bfhinventiones} and \cite{BFHgafa}; and in the recent proof by Deroin and Hurtado on non-left-orderability of higher rank lattices. Given a lattice $\Gamma$ of a Lie group $G$, along with a smooth action $\alpha:\Gamma\acts M$, $M$ a manifold, there is an induced smooth action $G\acts M^\alpha$, see \cite[Section~14.1]{brown2019entropylyapunovexponentsrigidity} and \cite[Section~4.1]{bfh}, encoding the original action $\alpha$, such that all the stabilizers of  $G\acts M^\alpha$ are discrete subgroups of $G$, i.e., such action is locally free. 

One of the main steps in the proofs of these results is constructing $G$-invariant probability measures on $M^\alpha$, whenever the dimension of $M$ is small enough with respect to $G$, by getting suitable entropy estimates by using a smooth and $G$-equivariant map $M^\alpha\rightarrow G/\Gamma$(coming from the induced construction), this is one of the main steps in \cite{brhw}, which is used in \cite{bfh},  \cite{bfhinventiones} and \cite{BFHgafa}. Note that $G$ is not amenable, see \cite[Example~9.1]{ZimmerMorris2008}, and thus, such $G$-invariant probability measures might, apriori, not exist.

It is natural to ask whether smooth locally free actions of higher rank simple Lie groups on closed manifolds of small dimension preserve a probability measure, This is not trivial, as G is not amenable, \cite[Example~9.1]{morris2009latticesslnractcircle}. For semisimple G with Property (T), which is the case for every simple Lie group of real rank at least two, an invariant probability measure on a locally free action with codimension one orbits forces the action to be as simple as possible: by work of Stuck and Zimmer, \cite[Example~5.1]{sz}, it is finitely and equivariantly covered by $G/\Gamma \times S^1$ for a cocompact lattice $\Gamma\leq G$, so in particular they are induced by finite actions of such lattices on the circle. The locally free hypothesis cannot be dropped by work of Katok-Lewis, \cite{MR1380646},  and  Benveniste-Fisher, \cite{MR2154667}, which construct actions of higher rank simple Lie groups with invariant ergodic measures that are not induced from an action of a closed subgroup of such $G$'s.

Theorem \ref{main} proves that any smooth, locally free, codimension one action of a simple higher rank Lie group admits an invariant probability measure, recall that $G$ is not amenable.




Theorem \ref{main smooth} continues the line of work on classification of smooth locally free actions of Lie groups. Rosenberg, Roussarie and D. Weil, see \cite{rrw}, prove that the only  connected, oriented, closed $3$-manifolds admitting a smooth locally free action of $\R^2$ are fiber bundles over a $k$-torus with fiber a $(3-k)$-torus; Ch\^atelet and Roussarie \cite{chateletrosenberg} prove that, for all $n\geq 2$, any closed $n$-manifold admitting a smooth locally free action of $\R^{n-1}$ is a fiber bundle over a $k$-torus with fiber a $(n-k)$-torus. For smooth locally free actions of the group of orientation preserving affine transformations of the line, $\operatorname{Aff}^+(\R)$, on closed $3$-manifolds, Ghys proved in his PhD thesis, see \cite{ghys85}, that if such an action admits a continuous volume form, invariant under the action, then the action is smoothly conjugate to a homogeneous action; Asaoka proved, see \cite{asaokaannals}, the hypothesis on the continuous invariant volume form cannot be dropped.

For locally free actions of higher rank Lie groups, recent important breakthroughs, in the framework of the Katok-Spatzier conjecture, by Spatzier and Vinhage, see \cite{spatziervinhage}; and by the same authors  with Damjanovi\'c and Xu, see \cite{damjanovic2025zimmerprogrampartiallyhyperbolic}, under ergodic/dynamical assumptions, more specifically,  for actions having enough elements in the group acting in a uniformly partially hyperbolic manner and under the existence of an invariant volume form.

For analytic actions of $G=\SL(n,\R)$ with $n\geq 3$ on closed manifolds of low dimension, recent work of Fisher and Melnick, see \cite{fishermelnick},  classifies all such actions on closed manifolds of dimension $n$; and recent work of Son, see \cite{son2025realanalyticmathrmslnmathbbractionsclosed}, classifies all such actions on manifolds of dimension $n\leq m\leq 2n-3$.


\subsection{Outline of the proof}\label{outline of the proof}

The proof of Theorem \ref{main} is done by contradiction, by fixing, for some simple higher rank Lie group $G$ and closed connected manifold $M$ of dimension $\dim(G)+1$,  a locally free and $C^2$ action $G\acts M$ that has no $G$-invariant probability measure. During this outline, we refer to this assumption as ``Main Assumption".

Fix $A$ a maximal $\R$-split torus of $G$, along with a minimal parabolic subgroup $P\leq G$ such that $A\leq P$. 

Since $P$ is amenable, fix a $P$-invariant and $P$-ergodic probability $\mu$. Such measures are $A$-ergodic by the Mautner phenomenon, see \cite[Chapter~3]{bm}, and since $\dim(M)=\dim(G)+1$, and the action is locally free, by the higher rank Osedelec's Theorem, see Theorem \ref{Higher rank osedelec} and \cite[Theorem~2.4]{brhwzd}, there exists a Lyapunov exponent $\lo\in A^*$ of the ergodic action $A\acts (M,\mu)$ whose Osedelec distribution is transversal to the orbits of $G\acts M$ for $\mu$-a.e. $x$.

We show that on the Weyl chamber of $P$ in $A$, denoted by $\w^P$, $\lo$ is strictly less than any root of $(G,A)$, by using entropy tools in \cite{brhwzd}. This yields, for $\mu$-a.e. $x$, a one dimensional Pesin manifold at $x$, denoted by $\w^T(x)$, see Definition \ref{definition transverse Pesin manifolds}, corresponding to the least Lyapunov exponent of the diffeomorphism $a:(M,\mu)\rightarrow(M,\mu)$, this manifold is independent of $a\in\w^P$.

\subsubsection{Case when $G$ is \emph{not} locally isomorphic to $\SL(3,\R)$}

In the case when $G$ is not locally isomorphic to $\SL(3,\R)$, we prove that $\lo$ is less than all roots of $(G,A)$ in the interior of the codimension one faces of $\overline{\w^P}$, using again \cite{brhwzd}. Applying this we can show that the partition $\mathscr W^T$ into transverse Pesin manifolds, seen as the partition into strong stable manifolds with respect to any $a\in\w^P$, jointly integrates with the partition by $P$-orbits, which can be seen as a central unstable partition for any $a\in\w^P$. In fact, something stronger happens, $\w^T(x)$ is equivariant under the action of $p$ for $\haar_P$-a.e. $p$, see Proposition \ref{P-equivariance of transverse Pesin manifolds}. We denote the unique $P$-invariant probability measure by $\mu^P$. 

We then show that for $\mu$-a.e. $x$, $\w^T(x)$ intersects all orbits of $G\acts M$, by showing that the restricted action $P\acts M$ is uniquely ergodic, see Proposition \ref{unique P ergodicity}.  

The aforementioned $P$-equivariance property, along with the fact that generic Pesin manifolds $\w^T(x)$ intersect all $G$-orbits, are combined to show that for $\mu^P$-a.e. $x$, $\{g\cdot\w^T(x)\}_{g\in G}$ is a $\lie(G)$-flow of $M$, see Lemma \ref{g-flow from generic pesin} and \cite{carriere1} for this notion. This in turn yields a smooth map $M\rightarrow G/\Gamma$ for some discrete subgroup $\Gamma\leq G$, whose fiber over the trivial coset is $\w^T(x)$. This yields a contradiction to ``Main Assumption", as under that assumption, one can prove that for $\mu^P$-a.e. $x$, $\w^T(x)$ is diffeomorphic to $\R$. 

\subsubsection{Case when $G$ \emph{is} locally isomorphic to $\SL(3,\R)$}

For this case we need to prove unique ergodicity of the restricted action $P\acts M$.  

We reduce this case to one where  there exists a $P$-invariant and $P$-ergodic measure $\mu$
such that for $\mu$-a.e. $x$, $\mu$ disintegrates equivalently to the Riemannian volume of $\w^T(x)$ on this manifold, see Lemma \ref{maximal parabolic stabilizer when SMI fails}. We use this property to show that the restricted action $P\acts M$ is uniquely ergodic, see Lemma \ref{unique P-ergodicity when rG=2}. If $\mu^P$ is this unique measure, from this unique ergodicity property, it is shown that for $\mu^P$-a.e. $x$, $\w^T(x)$ intersects every orbit of $G\acts M$.  

It is also proved that for $\mu$ as above, then $Q_0:=\operatorname{Stab}(\mu^P)$ has the following coherence property: for $\mu^P$-a.e. $x$, if $g\in G$ and $g\w^T(x)\cap\w^T(x)\neq\emptyset$, then $g\in Q_0$, see Proposition \ref{coherence proposition}; we prove this coherence property using a geometric version, see Theorem \ref{geometric high entropy method}, of the High Entropy Method of Eisindler and Katok, see \cite{ek}.

Combining the aforementioned coherence proposition along with $\w^T(x)$ intersecting all $G$-orbits of $\mu^P$-a.e. $x$, yields a $G$-equivariant and smooth map $M\rightarrow G/Q_0$, see Section \ref{proof of main theorem final section}. Using the argument given in the proof of \cite[Theorem~5.1]{sz} we prove that having such map implies that our original action $G\acts M$ has a closed orbit, which yields a contradiction to ``Main Assumption".

\subsection{Rank one Lie groups without Property (T) }\label{rank one counterexamples}

    \begin{example}
        Let $\Gamma_g$ be the fundamental group of the closed oriented surface of genus $g$, realized as a cocompact lattice in $\operatorname{SO}_0(2,1)$. Then there is an onto homomorphism $\Gamma_g\rightarrow \Z$, which then yields an action $\Gamma_g\acts S^1$, preserving Lebesgue measure,  all whose orbits are dense, by using an irrational rotation. We can then induce this action to an action $\operatorname{SO}_0(2,1)\acts M$ on a closed $4$-manifold, in particular with codimension one orbits, that has an invariant probability measure, and all of whose orbits are dense. Therefore, such an action cannot be covered by an action of the form $G\acts G/\Gamma\times S^1$, as in Theorem \ref{main smooth} with $G=\operatorname{SO}_0(2,1)$.
    \end{example}

    \begin{example}
        If $\Gamma\leq\operatorname{SO}_0(3,1)$ is the fundamental group of a closed hyperbolic manifold that fibers over the circle, then there is homomorphism from $\Gamma$ onto $\Z$. Then we can repeat the same argument as in the previous example, and obtain a locally free codimension one action of $\operatorname{SO}_0(3,1)$ on a closed manifold, with dense orbits.

        Furthermore, any cocompact lattice $\Gamma\leq \operatorname{SO}_0(3,1)$ has a finite index subgroup without torsion, by a classic result by Selberg, see \cite[Lemma~4.8.2]{wittemorrisbook}, so the latter is the fundamental group of a closed hyperbolic $3$-manifold, and such a group has a finite index subgroup which is the fundamental group of a closed hyperbolic manifold fibering over the circle, because of Agol's virtual fibering theorem, see \cite{MR3104553}. Hence, such examples can be built easily with any closed compact hyperbolic manifold.
    \end{example}

    \begin{example}
        By results of Millson, see \cite{MR422501}, for every $n\geq 4$ there is a lattice $\Gamma\rightarrow \operatorname{SO}_0(n,1)$ with a homomorphism onto $\Z$. Thus, we can construct examples where Theorem \ref{main smooth} fails for $\operatorname{SO}_0(n,1)$ for all $n\geq 4$.
    \end{example}

    \begin{example}
        There are cocompact lattices in $\operatorname{SU}(n,1)$, with homomorphisms onto $\Z$, see for example \cite[Lemma~3.5]{MR4021068}.
    \end{example}

\subsection{AI disclosure}

A previous version of this paper gave only a dichotomy for Theorem \ref{main}, namely, either the action $G\acts M$ had an $G$-invariant probability measure or there existed a $G$-equivariant and smooth map $M\rightarrow G/Q$, which gave a smooth version of the results of Nevo and Zimmer, \cite{nzann}, for codimension one actions of higher rank simple Lie groups. 

Homin Lee indicated to the author that the map $\sigma:G\times \w^T(x)\rightarrow M$ is a covering map(in the case where $G$ is not locally isomorphic to $\SL(3,\R)$), see Lemma \ref{sliding generic transverse pesin}, and then the author used an LLM model to ask whether such covering maps appeared in the literature, then the LLM model indicated that such covering map gave $M$ the structure of a $G$-flow, and referenced the work of Carri\'ere, see \cite{carriere1}. Carrier\'e defines this notion using deck transformations of the universal cover of a closed manifold, and this notion was all the author needed to build a map $M\rightarrow \tilde G/\Gamma$, see Corollary \ref{map psi_x}, in order to get the conclusion of Theorem \ref{main}, see Section \ref{proof of main when not SL3R Section}, in the case when $G$ is not locally isomorphic to $\SL(3,\R)$. 

LLM models were also used to double check basic facts on representation theory of Lie groups. All main ideas and contributions of the paper are human made.

\subsection*{Acknowledgments}
This work is an extension of the author's PhD thesis.

The author thanks his advisor, David Fisher, for his unwavering guidance and support. The author is also grateful to Sebastián Hurtado, whose enthusiasm and feedback were very important to this project.

Appreciation is extended to Aaron Brown, Ralf Spatzier, and Kurt Vinhage for their helpful suggestions, as well as to Rose Elliot Smith and Homin Lee for their comments. Finally, the author thanks Yale University, Northwestern University, The University of Michigan, and The University of Utah for their hospitality

\tableofcontents

\section{Preliminaries}

\subsection{Some Lie group preliminaries}\label{Some Lie group preliminaries}

\subsubsection{Root decomposition of semisimple Lie groups}

Let $G$ be a semisimple real Lie group with finite center and denote by $\mathfrak g$ its Lie algebra; let $A$ be a maximal $\R$-split torus of $G$. 

We denote by $\Sigma_{(G,A)}$ the of all roots of $G$ with respect to $A$. For a root $\alpha\in A^\ast$ of $(G,A)$, consider the corresponding root space
\begin{equation*}
    \mathfrak g_\alpha=\{X\in\mathfrak g: \operatorname{Ad}(s)X=\alpha(s)X, \forall s\in S\},
\end{equation*}
see \cite[Section~IV.2]{knapprep}. Set $m_\alpha=\dim\mathfrak g_\alpha$. Denote by $G^\alpha$ the corresponding analytic subgroup of $G$ associated to $\mathfrak g_\alpha$. Furthermore, we have the \emph{root decomposition} of $\operatorname{Lie}(G)$, 
\begin{equation}\label{root decomposition}
    \operatorname{Lie}(G)=\bigoplus_{\alpha\in\Sigma_{(G,A)}}\mathfrak g_\alpha.
\end{equation}

\subsubsection{Stable and Unstable subgroups}\label{Stable and Unstable subgroups section}
Fix a maximal compact subgroup $K$ of $G$. We fix a right $G$-invariant, bi-$K$-invariant metric $d_G$ on $G$. For a semisimple element $a$ of $G$ we consider the following subgroups of $G$:
\begin{equation}\label{Na and La}
\begin{split}
    N_a:=\{g\in G:\lim_{n\rightarrow-\infty} d_G(a^nga^{-n},e)=0\}\\
    L_a:=\{g\in G:\lim_{n\rightarrow\infty} d_G(a^nga^{-n},e)=0\}.
\end{split}
\end{equation}
called the \emph{unstable} and \emph{stable} subgroups of $G$ with respect to $a$, respectively.

We denote by $C_a$ the (group-theoretic) center of $a$ in $G$. If $a\neq e$, denote by $Q_a$ the subgroup of $G$ generated by $N_a$ and $C_a$, and call it the \emph{center unstable} subgroup of $G$ with respect to $a$.

For any closed subgroup $H$ of $G$ and any $\delta>0$, we set 
\begin{equation*}
    H^{<\delta}:=\{h\in H:d_G(h,e)<\delta\},\text{ }H^{\leq\delta}:=\{h\in H:d_G(h,e)\leq\delta\}.
\end{equation*}

\begin{remark}
     From the definitions above it follows that $C_a$ normalizes both $N_a$ and $L_a$.
\end{remark}


\begin{lemma}\label{multiplication maps}
    The maps $C_a\times N_a\rightarrow Q_a$ and $L_a\times C_a\times N_a\rightarrow G$ given by $(c,n)\mapsto cn$ and $(l,c,n)\mapsto lcn$ are diffeomorphisms of a neighborhood of the identity of $C_a\times N_a$ and $L_a\times C_a\times N_a$, respectively, onto a neighborhood of the identity of $Q_a$ and $G$, respectively.
\end{lemma}
\begin{proof}
    This follows by showing the Jacobian of each of the maps above has full rank, as $\lie(G)=\lie(L_a)\oplus\lie(C_a)\oplus\lie(N_a)$ and $\lie(Q_a)=\lie(C_a)\oplus\lie(N_a)$, by using the root space decomposition of $\lie(G)$ with respect to a maximal torus of $G$ containing $a$. 
\end{proof}

\subsubsection{Parabolic subgroups}\label{section parabolic subgroups}

Let $P$ be a minimal parabolic of $G$ containing $A$. We denote by $\Delta_P$ the set of simple roots of $(G,A)$ with respect to $P$, and by $\Sigma_P$ the set of all positive roots of $(G,A)$ with respect to $P$.

Given a subset $\theta$ of $\Delta_P$, denote by $P_\theta$ the parabolic subgroup of $G$ generated by $P$ along with all root groups $G^{-\alpha}$ of all $\alpha\in\Delta_P$, see \cite[Sections~I.1, I.2]{margulisbook}.

Given a parabolic subgroup $Q\leq G$ with $P\leq Q$, denote by $\Sigma_Q$ the set of roots of $(G,A)$ whose associated unipotent subgroups are contained in $Q$; i.e., if $Q=P_\theta$, then $\Sigma_Q=\Sigma_P\cup\{-\alpha:\alpha\in[\theta]\}$, where $[\theta]$ is the set of all roots of $(G,A)$ that are non-negative integer combinations of elements in $\theta$.

\subsection{Some smooth ergodic theory}

\subsubsection{Osedelec's multiplicative ergodic theorem and Lyapunov exponents}

For any $C^1$ diffemorphism $f:M\rightarrow M$, and any $f$-invariant  and ergodic probability measure $\mu$ on $M$ we have Osedelets' theorem, see \cite{osedelec}, and also  \cite[Theorem~7.1]{brown2019entropylyapunovexponentsrigidity}.

\begin{theorem}[Osedelets' Multiplicative Ergodic Theorem \cite{osedelec}]\label{Osedelets theorem}
There are:
\begin{enumerate}
    \item an $f$-invariant and measurable set $\Lambda$ with $\mu(\Lambda)=1$;
    \item real numbers $\lambda^1>\ldots>\lambda^p$;
    \item a $\mu$-measurable, $Df$-invariant splitting $T_xM=\bigoplus_{i=1}^r E^i(x)$, for all $x\in\Lambda$,
\end{enumerate}
    such that for all $x\in\Lambda$ and all $v\in E^i(x)\setminus\{0\}$ we have
    \begin{equation*}
        \lim_{n\rightarrow\pm\infty}\frac1n\log||D_xf^n(v)||=\lambda^i.
    \end{equation*}
\end{theorem}

The real numbers $\lambda^i$ are called \emph{Lyapunov exponents}, and the spaces $E^i(x)$ are called Osedelets' distributions. We call the $\mu$-almost surely value $m^i:=\dim E^i(x)$ the multiplicity of $\lambda^i$, for $\mu$-a.e. $x$.

\subsubsection{Unstable Pesin manifolds}\label{unstable pesin manifolds}

Given a compact manifold $M$, a $C^2$ diffeomorphism $f:M\rightarrow M$, and a $f$-ergodic probability measure $\mu$, with corresponding Lyapunov exponents $\lambda^1>\ldots>\lambda^d$, and Osedelets' distributions $E_i(-)$; there are associated global $j$th unstable manifolds, which are  $C^2$ injectively  immersed 
submanifolds, given by 
\begin{equation}\label{Pesin manifolds}
    \w^j(x)=\left \{y\in M:\limsup_{n\rightarrow\infty}\frac1n\log d(f^{-n}(x),f^{-n}(y))\leq -\lambda^j\right\}
\end{equation}
whenever $\lambda^j>0$, for $\mu$-a.e. $x$. These manifolds satisfy $T_x\w^j(x)=\bigoplus_{\lambda^j\leq\lambda^i}E^i(x)$ for $\mu$-a.e. $x$. See \cite[Section~4]{brhwzd}.

The set 
\begin{equation*}
    \w^u(x)=\left \{y\in M:\limsup_{n\rightarrow\infty}\frac1n\log d(f^{-n}(x),f^{-n}(y))< 0\right\}
\end{equation*}
is equal to the $r_0$th unstable manifold, where $1\leq r_0\leq p$ is such that $\lambda_{r_0}>0$, but $\lambda_{r_0+1}\leq 0$, if it exists; so that $\w^u(x)$ is a $C^2$ injectively immersed submandifold, for $\mu$-a.e. $x$, called \emph{the global unstable Pesin manifold} at $x$. It satisfies $T_x\w^u(x)=\bigoplus_{\lambda^i>0}E^i(x)=:E^u(x)$, for $\mu$-a.e. $x$.

We denote by $\mathscr W^j$ and $\mathscr W^u$ the corresponding partitions by Pesin manifolds.  

There are also corresponding stable and central distributions $E^s(x):=\bigoplus_{\lambda^i<0}E^i(x)$, $E^0(x):=\bigoplus_{\lambda^i=0}E^i(x)$, only the former  has Pesin manifolds $\w^s(x)$ tangent to $E^s(x)$, called \emph{stable Pesin manifolds}.

\subsubsection{The Ledrappier-Young theorem}

Let $f$ be a diffeomorphism of a closed manifold $M$. Let $\mu$ be an $f$-invariant and ergodic, probability measure on $M$. The metric entropy of $f$ with respect to $\mu$, denoted by $h_\mu(f)$, is equal to the entropy conditional to any measurable increasing  partition $\xi$ of $M$, subordinate to $\mathscr W^u$ with respect to $\mu$, see \cite[Section~8.3.1]{brown2019entropylyapunovexponentsrigidity}, this value is defined in \cite[Section~8.1.2]{brown2019entropylyapunovexponentsrigidity}; such measurable partitions always exists, and this value is independent of the measurable partition satisfying this properties, see \cite[Section~8.3]{brown2019entropylyapunovexponentsrigidity}.

For the $f$-ergodic probability measure $\mu$, consider  corresponding Lyapunov exponents $\lambda^1>\cdots>\lambda^p$, along with distributions $E^i(-)$ corresponding to $\lambda^i$, for $1\leq i\leq p$, as in Osedelec's Theorem, see Theorem \ref{Osedelets theorem}. Denote by $m_i$ the $\mu$-almost surely value $\dim E^i(x)$.

We have the following theorems relating the metric entropy of $f$ with respect to $\mu$, and the $\mu$-almost surely constant values given by Osedelec's Theorem.  

\begin{theorem}[Ruelle-Margulis inequality]\label{Ruelle-Margulis inequality}
    \begin{equation*}
        h_\mu(f)\leq \sum_{\lambda^i>0}m^i\lambda^i.
    \end{equation*}
\end{theorem}

$\mu$ is called an SRB measure if for any measurable partition $\xi$ of $M$ subordinate to $\mathscr W^u$ with respect to $\mu$, we have that for $\mu$-a.e. $x$, the conditional measure $\mu^\xi_x$ is absolutely continuous with respect to the Riemannian volume on $\w^u(x)$.

The Ledrappier-Young Theorem, see \cite[Theorem~A]{ledrappieryoung1} says that the $f$-ergodic probability measures $\mu$ having equality in Theorem \ref{Ruelle-Margulis inequality}, are exactly the SRB measures. That SRB measures have such an equality was first proved by Ledrappier and Strelcyn, see \cite{ledrappierstrelcyn}. Furthermore, when such equality holds, the corresponding Radon-Nikodym densities are $C^1$ and strictly positive, see the remark after \cite[Theorem~A]{ledrappieryoung1}.

\begin{theorem}[Ledrappier-Young theorem, \cite{ledrappieryoung1}]\label{ledrappier young}
    We have that
    \begin{equation*}
        h_\mu(f)= \sum_{\lambda^i>0}m^i\lambda^i,
    \end{equation*}
    if and only if for any measurable partition $\xi$ subordinate to $\mathscr W^u$ with respect to $\mu$ we have that $\mu_x^\xi$ is equivalent to the Riemannian volume on  $\w^u(x)$, for $\mu$-a.e. $x$. Furthermore,  the corresponding densities are positive  $C^1$ functions.
\end{theorem}

\subsubsection{Entropy conditional on measurable foliations} We follow \cite[Chapter~7]{brhwzd}.
Let $f:M\rightarrow M$ be a $C^2$ diffeomorphism, and $\mu$ an ergodic $f$-invariant probability measure on $M$. 

Let us consider measurable $f$-invariant foliations $\mathcal F$ of $M$ by $C^2$ leaves, which are \emph{tame measurable foliations}, that is, $\mathcal F$ is a partition of $M$ by $C^2$ manifolds, such that restricting the foliation to sets of large measure, $\mathcal F$ has the structure of a continuous family of $C^2$ disks. Examples of such partitions are the partition into Pesin manifolds, $\mathscr W^j$ and $\mathscr W^u$ constructed above, see \cite[Chapter~7]{bpintro}. We say $\mathcal F$ is \emph{expanding} if $\mathcal F(x)\subseteq\w^u(x)$ for $\mu$-a.e. $x$. 

For any $\xi$ measurable partition subordinate to an expanding tame measurable foliation $\mathcal F$ of $M$, we define the entropy of $f$ with respect to $\mu$ conditional on $\mathcal F$ to be
\begin{equation*}
    h_\mu(f\mid\mathcal F):=h_\mu(f\mid\xi);
\end{equation*}
such measurable partitions always exist, and the value above does not depend on the particular measurable partition subordinate to $\mathcal F$.

If $\mathcal F$ is a tame measurable foliation, we set 
\begin{equation*}
    h_\mu(f,\mathcal F)=h_\mu(f,\mathcal F^u),
\end{equation*}
where $\mathcal F^u$ is the expanding, tame measurable partition given by 
\begin{equation*}
    \mathcal F^u:=\mathcal F\vee\mathscr W^u:=\{A\cap B: A\in\mathcal F, B\in \mathscr W^u\}.
\end{equation*}

\subsubsection{Conditional entropy along measurable partitions and the Ledrappier-Young Theorem}

We follow \cite[Section~4, Section~7]{brhwzd}.

We have the following version of the Ledrappier-Young inequality for conditional entropies, see \cite[Theorem~7.2]{brhwzd} and \cite[Corollary~6.1.4]{ledrappieryoung1}. Define the multiplicity of $\lambda^i$ relative to $\mathcal F$ to be the $\mu$-almost surely constant value of
\begin{equation*}
    m_i(\mathcal F):=\dim(E^i(x)\cap T_x[\mathcal F(x)])
\end{equation*}

\begin{theorem}\cite[Corollary~6.1.4]{ledrappieryoung1}\cite[Theorem~7.2]{brhwzd}\label{ledrappier young conditional entropy}
    We have
\begin{equation*}
    h_\mu(f\mid\mathcal F)\leq\sum_{1\leq i\leq r}\lambda_i m_i(\mathcal F).
\end{equation*}
Moreover, equality holds if and only if for any measurable parition $\eta$ the following equivalent conditions hold:

\begin{enumerate}
    \item for any measurable partition $\xi$ subordinate to $\mathcal F^u$, with respect to $\mu$, we have that $\mu_x^\xi$ is absolutely continuous with respect to the volume measure on $\xi(x)$
    \item for any measurable partition $\xi$ subordinate to $\mathcal F^u$ we have that $\mu_x^\xi$ is equivalent to the Riemannian volume on  $\xi(x)$, for $\mu$-a.e. $x$. Furthermore,  the corresponding densities are positive  $C^1$ functions.
\end{enumerate}
\end{theorem}

\subsubsection{Geometric characterization of the defect in the entropy formula}\label{Geometric characterization of the defect in the entropy formula section}

We follow \cite[Chapter~7]{brhwzd}. Let $f:M\rightarrow M$ be a $C^2$ diffeomorphism, and $\mu$ an ergodic, $f$-invariant probability measure on $M$. Let $\lambda^1>\ldots>\lambda^p$ be the 

Let $\eta$  be any measurable partition of $(M,\mu)$. For $1\leq i\leq r$, let $\xi^i$ be a measurable partition subordinate to the partition associated to the $i$th unstable Pesin manifolds, i.e., $\mathscr W^i$. We define the \emph{$i$th upper and lower pointwise dimensions of $\mu$ relative to $\eta$} at $x$ to be the limits
\begin{equation}\label{pointwise dimensions}
\begin{split}
    \overline{\dim}^i (\mu, x\mid\eta) &:= \limsup_{\delta \to 0} \frac{\log \left( \mu_x^{\xi^i \vee \eta}(B(x, \delta)) \right)}{\log \delta},\\ \underline{\dim}^i (\mu, x\mid\eta) &:= \liminf_{\delta \to 0} \frac{\log \left( \mu_x^{\xi^i \vee \eta}(B(x, \delta)) \right)}{\log \delta}
\end{split}
\end{equation}
where $B(x,\delta)$ is the ball of radius $\delta$ centered at $x$ of a Riemannian
metric on $M$, these values are independent of such metrics, and of the chosen measurable partitions. These values are constant $\mu$-almost surely, and we  denote, respectively, these values by $\overline{\dim}^i(\mu\mid\eta)$ and $\underline{\dim}^i(\mu\mid\eta)$

For the partition $\eta^+:=\bigvee_{j=0}^\infty f^j\eta$, the values $\overline{\dim}^i(\mu\mid\eta^+)$ and $\underline{\dim}^i(\mu\mid\eta^+)$ coincide, see \cite[Proposition~7.4]{brhwzd}, and we denote this common value by $\dim^i(\mu\mid\eta^+)$. 

Set $\dim^0(\mu\mid\eta^+) = 0$. For $1 \le i \le r$ the \textit{$i$th transverse dimension of $\mu$ relative to $\eta^+$} is
$$
\gamma^i(\mu\mid\eta^+):= \dim^i(\mu\mid\eta^+) - \dim^{i-1}(\mu\mid\eta^+).
$$
The $i$th transverse dimension is not bigger than the multiplicity $m_i$ of the Osedeltec's distributions $E^i$, see \cite[Claim~7.5]{brhwzd}.

Using the pointwise transverse dimensions above we can indicate the exact failure for equality to hold in the Ruelle-Margulis inequality (\ref{Ruelle-Margulis inequality}), this was proved by Ledrappier and Young.

\begin{theorem}\cite{ledrappieryoung2}\label{ledrappieryoung2}
    Let $\eta$ be any measurable partition of $(M,\mu)$. Then
    \begin{equation*}h_\mu(f\mid\eta)=\sum_{\lambda_i>0}\lambda_i\gamma^i(\mu\mid\eta^+).
    \end{equation*}
\end{theorem}

\subsection{Smooth ergodic theory of nonuniform partially hyperbolic diffeomorphisms}\label{Smooth ergodic theory of nonuniformly partially hyperbolic diffeomorphisms}

We follow \cite{bpintro} and \cite{bp}. Let $f$ be a $C^1$ diffeomorphism of a closed manifold $M$. 

\subsubsection{Nonuniformly partially hyperbolic sets}

\begin{definition}\label{definition of nonuniformly partially hyperbolic subset}

We call an $f$-invariant nonempty measurable subset $Z \subset M$ \textit{nonuniformly partially hyperbolic} if there are  numbers $\varepsilon\geq 0,\lambda_1, \lambda_2, \mu_1, \mu_2 > 0$; measurable functions $ C, K : Z \to (0, \infty)$; subspaces $E^s(x)$, $E^u(x)$, and $E^0(x)$, which depend measurably on $x \in Z$, such that the following conditions hold for all $x\in Z$:

\begin{itemize}
    \item[(1)]
    \[
    \lambda_1e^{\varepsilon} < \lambda_2e^{-\varepsilon} < 1 < \mu_2e^{\varepsilon} < \mu_1e^{-\varepsilon};
    \]
    
    \item[(2)] for $m\in \Z$
    \begin{equation*}
        C(f^m(x))\leq C(x)e^{\varepsilon|m|}, K(f^m(x))\geq K(x)e^{-\varepsilon|m|}
    \end{equation*}
    \item[(3)] $T_x M = E^s(x) \oplus E^0(x) \oplus E^u(x)$ and for $a = s, u, 0$,
    \[
    d_x f E^a(x) = E^a(f(x));
    \]
    
    \item[(4)] for $v \in E^s(x)$, $w \in E^u(x)$, and  $n \ge 0$ we have
    \begin{equation*}
        ||d_xf^nv||\leq C(x)\lambda_1^ne^{\varepsilon n}||v||, ||d_xf^{-n}w||\leq C(x)\mu_1^{-n}e^{\varepsilon n}||v|| 
    \end{equation*}
    
    \item[(5)] for $v \in E^0(x)$
    \[
    C(x)^{-1}\lambda_2^n e^{-\varepsilon n}\|v\| \le \|d_x f^n v\| \le C(x)\mu_2^n e^{\varepsilon n}\|v\|, \quad n \ge 0,
    \]
    \[
    C(x)^{-1}\mu_2^ne^{\varepsilon n} \|v\| \le \|d_x f^n v\| \le C(x)\lambda_2^n e^{-\varepsilon n}\|v\|, \quad n \le 0;
    \]
    
    \item[(6)] the angles satisfy $\angle(E^a(x), E^b(x)) \ge K(x)$ where $a, b \in \{s, u, 0\}$ and $a \ne b$.
\end{itemize}
\end{definition}

\subsubsection{Regular sets of nonuniformly partially hyperbolic sets}

\begin{definition}\label{definition regular sets}
For $Z$  a nonuniformly partially hyperbolic set of a $C^1$ diffeomorphism $f:M\rightarrow M$, and any $\ell>0$, we can consider \emph{regular sets} given by 
\begin{align*}
 Z^\ell = \{ & x \in Z : \ C(x) \leq \ell, \ K(x) \geq \ell^{-1} \}\label{regular sets page}.
\end{align*}
\end{definition}

Such sets are compact; satisfy for $0<\ell\leq\ell'$, $ Z^\ell\subseteq Z^{\ell'}$; and are such that
\begin{equation*}
    Z:= \bigcup_{\ell>0} Z^\ell.
\end{equation*}

\begin{remark}\label{nonuniform partially hyperbolic from multiplicative}
   If $f$ is a diffeomorphism $M\rightarrow M$ and $\mu$ is an $f$-invariant and ergodic probability measure, then the Oseledets' Multiplicative Ergodic Theorem, see Theorem \ref{Osedelets theorem}, yields a nonuniformly partially hyperbolic subset of $M$ for $f$, of full $\mu$ measure, see \cite[Section~2.2, 4.5]{bp}.
\end{remark}

\begin{definition}\label{local stable manifolds}
Each $x$ in such regular sets
$ Z^\ell$ has a local stable manifold, denoted by $\w^s_{\operatorname{loc}}(x;f)$(depending on $\ell$; we drop the $f$ for $\w^s_{\operatorname{loc}}(x;f)$ when context is clear), which is a connected, injectively immersed submanifold of $M$ having the dimension of $E^s(x)$, tangent to this distribution at $x$; such that for some $T(x)>0$ and $0<\lambda_x<1$(depending on $\ell$), for all $y\in\w^s_{\operatorname{loc}}(x;f)$ we have $$d_M(f^n(x),f^n(y))<T(x)\lambda_x^nd_M(x,y),$$
for all $n>0$. Furthermore, $\mathcal W_{\operatorname{loc}}^s(x;f), T(x)$ and $\lambda_x$ depend continuously in the $C^1$-norm on $x\in  Z^\ell$, see \cite[Exercise~7.9.4]{bpintro}. 
    
\end{definition}

\subsubsection{Absolute continuity of holonomy along stable manifolds}\label{Absolute continuity along stable manifolds section}
 
Given $x\in Z^\ell$, consider the union 
\begin{equation}\label{union of family of local stables}
    \mathcal{Q}^\ell(x)=\bigcup_{w\in  Z^\ell\cap B^M_{r^l}(x)}\locstable(w, f)
\end{equation}
of the family of local stable manifolds(with respect to $ Z^\ell$)
\begin{equation}\label{family of local stable manifolds}
    \mathcal L^\ell(x)=\left \{\locstable(w):w\in Z^\ell\cap B^M(x,r_l)\right \},
\end{equation}
for any $r_l>0$ lower bound of the diameter of these stable manifolds.

A \textit{local smooth submanifold} of $M$ is the image of a smooth injective function $\mathbb{R}^k \to M$, with precompact image. A \textit{local transversal} to the family $\mathcal{L}^l(x)$ is a local smooth submanifold which is uniformly transverse to every manifold in this family, i.e., the angle between $T$ and $\mathcal{W}_{\text{loc}}^s(w)$ is uniformly away from zero for every $w\in Z^\ell\cap B^M(x,r^l)$; and also, $T\cap\locstable(w)$ consists of a single point. For any $T^1$ and $T^2$ such local transversals to $\mathcal{L}^l(x)$, of dimension $\dim(E^0(x))+\dim(E^u(x))$, we have a corresponding holonomy map
\begin{equation}\label{holonomy between transversals along stable}
    \mathcal{Q}^\ell(x) \cap T^1 \to \mathcal{Q}^\ell(x) \cap T^2, 
\end{equation}
by setting $\locstable(w)\cap T^1= z_1 \mapsto z_2=\locstable(w)\cap T^2$, for every $w\in Z^\ell\cap B^M(x,r^l)$.

The local transversal submanifolds we will consider are ones that are the $(u0)$-\emph{canonical local transversal} through $x$, given by 
\begin{equation*}
    T^{u0}(x):=\exp_x(B^u(0,r^l/2)\oplus B^0(0,r^l/2)),
\end{equation*}
where $B^\sigma(0,r^l/2)$ is the ball of radius $r^l/2$ in $E^\sigma(x)\subseteq T_xM$ centered at the origin, for $\sigma\in\{s,0,u\}$. Then, any local smooth submanifold of $M$ close enough in the $C^1$-norm to $T^{u0}(x)$ is a local transverse to the family $\mathcal L^l(x)$.

For a local smooth submanifold $T$ of $M$ we denote by $\nu_T$ the volume measure on $T$ inherited from the restriction of the Riemannian metric of $M$ to $T$.

\begin{theorem}\cite[Theorem~8.2, Remark~8.19]{bpintro} \label{absolute continuity along stable}
    If in (\ref{holonomy between transversals along stable}) $\nu_{T^1}(T^1\cap \mathcal{Q}^l(x))>0$, the holonomy map (\ref{holonomy between transversals along stable}) is absolutely continuous, with respect to the volume measures $\nu_{T_1}$ and $\nu_{T_2}$.
\end{theorem}

\subsection{Smooth ergodic theory of higher rank abelian actions}\label{Smooth ergodic theory of higher rank abelian actions}

\subsubsection{Higher-rank Osedelec's Theorem }

For the rest of this section, we follow \cite{brhwzd}.

Let $M$ be a closed manifold and $\alpha:\Z^d\rightarrow\operatorname{Diff}^1(M)$ an action. Let $\mu$ be an invariant and $\alpha$-ergodic probability measure. 

\begin{theorem}[Higher Rank Osedelec's Theorem]\cite[Theorem~2.4]{brhwzd}\label{Higher rank osedelec}
    There exist
    \begin{enumerate}
        \item an $\alpha$-invariant meausurable set $\Lambda\subseteq M$ such that $\mu(\Lambda)=1$;
        \item linear functionals $\lambda^1,\ldots,\lambda^p:\R^d\rightarrow\R;$
        \item a $\mu$-measurable, $D\alpha$ invariant splitting $T_xM=\bigoplus_{i=1}^p E^i(x)$, defined for $\mu$-a.e. $x$;
        \item for $\mu$-a.e. $x$ we have that for every $v\in E^i(x)\setminus\{0\}$, 
        \begin{equation*}
            \lim_{|n| \to \infty} \frac{\log \| D_x \alpha(n)(v) \| - \lambda^i(n)}{|n|} = 0
        \end{equation*}
    \end{enumerate}
\end{theorem}
We call the functionals $\lambda^1,\ldots,\lambda^p$ the Lyapunov exponents of $\alpha:\Z^d\acts (M,\mu)$; $E^i$ the Osedelec's distribution or vector space associated to $\lambda^i$; along with the almost surely constant multiplicity of $\lambda^i$, given by $m_i:=\dim E^i(x)$, for $\mu$-a.e. $x$.

\subsubsection{Global unstable Pesin manifolds in higher rank actions}\label{Global unstable Pesin manifolds in higher rank actions}
Let $\alpha:\Z^d\times (M,\mu)\rightarrow(M,\mu)$ be an ergodic action on a finite measure space. Let $\mathcal L=\{\lambda_i:\R^d\rightarrow\R:1\leq i\leq p\}$ be the set of corresponding Lyapunov exponents of this action. Given any $\vec n\in\Z^d$ there is a permutation $\Delta(\vec n)$ of $\{1,\ldots,p\}$ and $u(\vec n)$ in this set so that
\begin{equation*}
    \lambda_{\Delta(\vec n)(1)}(\vec n)\geq \lambda_{\Delta(\vec n)(2)}(\vec n)\geq\cdots\geq \lambda_{\Delta(\vec n)(u(\vec n))}(\vec n)>0\geq\cdots\geq\lambda_{\Delta(\vec n)(p)}(\vec n).
\end{equation*}

\begin{definition}\label{definition ith unstable}
    For any $x\in M$ and any $\vec n\in\Z^d$, and $1\leq i\leq u(\vec n)$, define the \emph{$i$th unstable manifold through $x$ for $\alpha(\vec n)$} as the set
\[
W^{u,i}_{\vec{n}}(x) = \left\{ y \in M \,\middle|\, \limsup_{k \to -\infty} \frac{1}{k} \log d\big(\alpha(k\vec{n}, x), \alpha(k\vec{n}, y)\big) \leq -\lambda_{\Delta(\vec{n})(i)}(\vec{n}) \right\}.
\]

The \textit{unstable manifold through \( x \) for \( \alpha(\vec{n}) \)} is the set
\[
W^u_{\vec{n}}(x) := \left\{ y \in M \,\middle|\, \limsup_{k \to -\infty} \frac{1}{k} \log d\big(\alpha(k\vec{n}, x), \alpha(k\vec{n}, y)\big) < 0 \right\}.
\]

\end{definition}

It can be shown that the unstable manifold $W^u_{\vec n}(x)$ with respect to the diffeomorphism $\alpha(\vec n)$ coincides with $\w^{u,u(\vec n)}_{\vec n}(x)$ for $\mu$-a.e. $x$. The manifold $W^{u,i}_{\vec{n}}(x)$ is tangent to $\bigoplus_{1\leq j\leq i}E_{\lambda_{\Delta(\vec n)(j)}}(x)$ for $\mu$-a.e. $x$.

We have the following uniqueness property of the $i$th unstable manifolds defined above.

\begin{lemma}\cite[Lemma~4.7]{brhwzd}\label{uniqueness of ith unstable higher rank}
    Let \( \vec{n_1}, \vec{n_2} \in \mathbb{Z}^d \). If for some \( 1 \leq i \leq \min\{u(\vec{n_1}), u(\vec{n_2})\} \) we have the following properties:
\begin{enumerate}
  \item \( \{\Delta(\vec{n_1})(j) : 1 \leq j \leq i\} = \{\Delta(\vec{n_2})(j) : 1 \leq j \leq i\} \);
  \item \( \lambda_{\Delta(\vec{n_1})(i)}(\vec{n_1}) > \lambda_{\Delta(\vec{n_1})(i+1)}(\vec{n_1}) \);
  \item \( \lambda_{\Delta(\vec{n_2})(i)}(\vec{n_2}) > \lambda_{\Delta(\vec{n_2})(i+1)}(\vec{n_2}) \),
\end{enumerate}
then \( \mathscr{W}^{u,i}_{\vec{n_1}}\overset{\circ}{=}\mathscr{W}^{u,i}_{\vec{n_2}} \).
\end{lemma}

For two tame measurable foliations $\mathcal F$ and $\mathcal G$, the relation $\mathcal F\overset{\circ}{=}\mathcal G$ means that $\mathcal F(x)=\mathcal G(x)$ for $\mu$-a.e. $x$.

\subsubsection{Coarse Lyapunov exponents and coarse Pesin manifolds}\label{Coarse Lyapunov exponents and Pesin manifolds for higher rank actions}

For  $\alpha$, $\mu$ and $\mathcal L$ as above, we say $\lambda^i, \lambda^j\in\mathcal L$  are \emph{resonant} if there exists some $c>0$ such that $\lambda^i=c\lambda^j$. The corresponding equivalence classes are called \emph{coarse Lyapunov exponets}, and we denote the corresponding set of equivalence classes by $\hat {\mathcal L}$. 

Given $\chi\in\hat{\mathcal L}$, with $\chi\neq\{0\}$(here $0$ is the zero element of $A^\ast$), the \emph{coarse Lyapunov foliation} corresponding to $\chi$ is given by 
\begin{equation*}
    \mathscr W^\chi:=\bigvee_{\{\vec n\in\Z^d:\chi(\vec n)>0\}}\mathscr W_{\vec n}^u.
\end{equation*}
With corresponding leaves denoted by $\w^\chi(x)$. $\mathscr W^\chi$ is a tame measurable foliation.

\subsubsection{Product structure of entropy for higher rank actions}

We have the following product structure theorem for higher rank actions $\alpha:\Z^d\acts (M,\mu)$.

\begin{theorem}\cite[Corollary~13.2]{brhwzd}\label{entropy product structure}
    For any $\vec n\in\Z^d$ we have 
    \begin{equation*}
        h_\mu(\alpha(\vec n))=\sum_{\{\chi\in\tilde{\mathcal L}:\chi(\vec n)>0\}}h_\mu(\alpha(\vec n)\mid\mathscr W^\chi)
    \end{equation*}
\end{theorem}

The previous result is an immediate consequence of the following version of it for conditional entropy.

\begin{theorem}\cite[Theorem~3.1]{brhwzd}\label{conditional entropy product structure}
    Let $\mathcal F$ an $\alpha$-invariant, $C^{1+\beta}-$ tame measurable foliation. For $\vec n\in\Z^d$ we have 
    \begin{equation*}
        h_\mu(\alpha(\vec n)\mid\mathcal F)=\sum_{\{\chi\in\tilde{\mathcal L}:\chi(\vec n)>0\}}h_\mu(\alpha(\vec n)\mid\mathscr W^\chi\vee\mathcal F)
    \end{equation*}
\end{theorem}

\subsection{Invariance from entropy considerations}\label{Invariance from entropy considerations}

In this subsection we follow \cite[\S~6]{el}, \cite{brown2019entropylyapunovexponentsrigidity}, and \cite {brhw}.

Let $G$ be a Lie group, and suppose $G$ acts, to the left, smoothly and locally freely on a closed manifold $M$. Suppose $H\leq G$ is a closed subgroup, and suppose $s\in G$ normalizes $H$. Let $\mu$ be an ergodic, $s$-invariant probability measure on $M$(we mean ergodicity and invariance with respect to the diffeomorphism of $M$ given by $x\mapsto sx$). Suppose the orbit $H\cdot x$ is contained in the unstable manifold $\w^u(x;s)$ for $\mu$-a.e. $x$. A measurable partition $\eta$ of $M$ is \emph{subordinate to the $H$-orbits of $M$}, with respect to $\mu$, if for $\mu$-almost every $x$ we have that $\eta(x)$ is contained in $H\cdot x$; and $\eta(x)$ contains a neighborhood of $x$(in the immersed topology) of the orbit $H\cdot x$. Then, see \cite[Definition~8.1]{brown2019entropylyapunovexponentsrigidity}, there is family of Borel probability measures $\{\mu_x^\eta\}_{x \in X}$, called a \emph{family of conditional measures} of $\mu$ relative to $\eta$, with the following properties: For $\mu$-a.e. $x$,
\begin{enumerate}
    \item[(1)] $\mu_x^\eta$ is a Borel probability measure on $X$ with $\mu_x^\eta(\eta(x)) = 1$;
    \item[(2)] if $y \in \eta(x)$ then $\mu_y^\eta = \mu_x^\eta$.
\end{enumerate}
Moreover, if $D \subset X$ is a Borel subset then
\begin{enumerate}
    \item[(3)] $x \mapsto \mu_x^\eta(D)$ is measurable and
    \item[(4)] $\mu(D) = \int \mu_x^\eta(D) \, d\mu(x)$.
\end{enumerate}
Such a family is unique modulo $\mu$-null sets.

\begin{lemma}\cite[Lemma~8.4] {brown2019entropylyapunovexponentsrigidity}\label{zero entropy atomic}
The following are equivalent:
\begin{enumerate}
    \item[\textit{(1)}] $h_\mu(s) = 0$;
    \item[\textit{(2)}] for $\mu$-a.e. $x$, the conditional measure $\mu_x^\eta$ has at least one atom;
    \item[\textit{(3)}] for $\mu$-a.e. $x$, the conditional measure $\mu_x^\eta$ is a single atom supported at $x$;
    \item[\textit{(4)}] the partition of $(M, \mu)$ into full $W^u$-manifolds is a measurable partition.
\end{enumerate}
\end{lemma}

\begin{lemma}\cite[Theorem~9.4]{brown2019entropylyapunovexponentsrigidity}\label{invariance for subgroups with orbits in unstable}
$\mu$ is $H$-invariant if and only if for any measurable partition $\eta$ subordinate to the partition into $H$-orbits, and for $\mu$-a.e. $x$, the conditional measure $\mu^\eta_x$ coincides(up to normalization) with the restriction of the left-Haar measure on $H\cdot x$ to $\eta(x)$.
\end{lemma}

We define the entropy of $s$, seen as an element of $\operatorname{Diff}(M)$, conditioned on $H$-orbtis, as the entropy conditioned on any measurable partition $\eta$ subordinate to the $H$-orbits, denoted by $h_\mu(s\mid H)$, namely,
\begin{equation*}
    h_\mu(s\mid H):=h_\mu(s\mid \eta).
\end{equation*}

Let $\lambda^i$, $E^i(x)$, and $m^i$ be as in Theorem \ref{Osedelets theorem} for the dynamics of $s$ and the measure $\mu$.  We define the \emph{multiplicity of $\lambda^i$ relative to $H$} to be (the almost surely constant value of)  $$m^{i,H} = \dim (E^i(x) \cap T_x (H\cdot x)).$$

We have the following version of the Ledrappier-Young Theorem, see \cite{ledrappieryoung1}, characterizing $H$-invariant measures. 

\begin{lemma}\cite[Theorem~9.5]{brown2019entropylyapunovexponentsrigidity}\label{invariance from maximal entropy}
    The following are equivalent:
		\begin{enumerate} 
			\item $h_\mu(f\mid H) = \sum_{\lambda^i>0} \lambda^i m^{i,H}$;
			\item $\mu$ is $H$-invariant.  
		\end{enumerate}
\end{lemma}

To any Borel probability measure $\mu$ on $M$, there is a family $\{\mu_x^H\}_{x\in M}$ of measures of $H$, called \emph{leafwise measures along the orbits of the action} $H\acts M$, see \cite[Section~6]{el}, having the following property: for any measurable partition $\eta$ of $(M,\mu)$ subordinate to the $H$-orbits, with respect to $\mu$, there exists a measurable function $c^\eta:X\rightarrow (0,\infty)$ such that if $\{\mu_x^\eta\}_{x\in M}$ is a family of conditional measures on $(M,\mu)$ associated with $\eta$, then 
\begin{equation*}
    \mu_x^\eta=c^\eta(x)(h\mapsto h\cdot x)_\ast(\mu_x^H)\upharpoonright_{\eta(x)}.
\end{equation*}
We have the following immediate consequence of Lemma \ref{invariance from maximal entropy}:
\begin{proposition}\cite[Problem~6.28]{el}\label{haarconditionals}
    A probability measure $\mu$ on $M$ is $H$-invariant if and only if $\mu_x^H$ coincides with the (positive proportionally class of the) left Haar measure on $H$ for $\mu$-almost every $x$.
\end{proposition}

On the opposite side of behavior, i.e., when $\mu^H_x$ is atomic for $\mu$-a.e. $x$, we have the following.

\begin{lemma}\label{zero entropy atomic leafwise}
If $\mu_x^H$ is atomic for $\mu$-a.e. $x$, then $h_\mu(s\mid H)=0$.
\end{lemma}

\begin{proof}
     This follows directly from the definition of metric entropy, see \cite[Section~8.1]{brown2019entropylyapunovexponentsrigidity}.
\end{proof}






\section{Transverse Pesin Manifolds and Equivariance}\label{Transverse Pesin Manifolds and Equivariance}

 We begin the proof of Theorems \ref{main} in this section. The proof is by contradiction  by assuming that $M$ does not have an invariant measure under $G\acts M$.

To achieve this, from this point of the document on, we work under the following assumption: 

\begin{assumption}\label{assumption no invariant probability}
   Let $G$ be a connected, simple Lie group, with finite center, all of whose simple factors are higher rank; let $A\leq G$ is an $\R$-split maximal torus; let $P$ be a minimal parabolic of $G$ containing $A$. Fix $G\acts M$  a locally free action, where $M$ is connected and closed and $\dim(M)=\dim(G)+1$. \emph{From this point of the paper on}, assume $G\acts M$ has no invariant probability measure.
\end{assumption}


\subsection{Some observations on $P$-ergodic measures}
\subsubsection{Birkhoff ergodic theorem and measures invariant under unstable subgroup}

Let $a\in A$. We need the following observation about $a$-ergodic measures that are also invariant under the unstable subgroup of $a$ in $G$, see also \cite[Proposition~8.4]{dh}. Recall the definition of the group $N_a$ in Section \ref{Stable and Unstable subgroups section}.

\begin{lemma}\label{birkhoff for almost all px}
    Let $a\in A$, and suppose $\mu$ is a $N_a$-invariant and $a$-ergodic probability measure on $M$. For $\mu$-a.e. $x$ and a.e. $u\in N_a$, the Birkhoff ergodic theorem holds for $ux$ with respect to $\mu$ and $a$.
\end{lemma}

\begin{proof}  
Since $N_a$ is, by construction, the unstable subgroup of $a$ in $G$, we obtain from Proposition \ref{haarconditionals} that the Birkhoff ergodic theorem holds for $ux$ with respect to $\mu$ and $a$, for $\mu$-a.e. $x$ and $\haar_{N_a}$-a.e. $u$, since this theorem holds for $\mu$-a.e. $x$. 
\end{proof}

\subsubsection{$a$-ergodic measures and transverse Lyapunov exponents}

\begin{definition}[Transverse Lyapunov exponent]
    Take $a\in A$. Any $a$-ergodic and $a$-invariant probability measure on $\mu$ on $M$ has a corresponding Osedelets' decomposition of $TM$, on a subset of full measure of $M$, along with corresponding Lyapunov exponents, see Theorem \ref{Osedelets theorem}. Note that we have canonical Lyapunov exponents, i.e., roots of the maximal torus containing $a$, coming from the Cartan flow of the locally free action $G\acts M$. Since $\dim(M)=\dim(G)+1$, there exists a unique Lyapunov exponent of $a\acts (M,\mu)$ whose corresponding Osedelets' distribution is not tangent to the $G$-orbits, which we denote by $\lambda_\mu^T(a)$, and call \emph{the transverse Lyapunov exponent} of the $a$-ergodic measure $\mu$. 
\end{definition}

\subsubsection{$A$-ergodic measures and transverse Lyapunov exponents}
\begin{definition}\label{definition transverse Lyapunov A-ergodic}
    Since $G\acts M$ is locally free, and its orbits are of codimension $1$, we obtain that for any $A$-invariant and ergodic probability measure $\mu$ on $M$ there is a unique Lyapunov exponent, which we shall denote by $\lambda_\mu^T\in A^\ast$ and call \emph{transversal Lyapunov exponent(functional)}, with respect to the $A$-ergodic measure $\mu$, for which its corresponding Osedelets' distribution is not tangent to the $G$-orbits.
\end{definition}



\subsection{$P$-ergodic measures}

\begin{notation}
    For the rest of this section, Section \ref{Transverse Pesin Manifolds and Equivariance}, $\mu$ denotes a fixed arbitrary $P$-invariant and $P$-ergodic Radon probability measure on $M$, where we are considering the corresponding restrictions of the action $G\acts M$ to $P$.
\end{notation}

\begin{remark}\label{mu is A-ergodic}
    By the Mautner Phenomenon, see \cite[Theorem~III.1.4]{bm}, $\mu$ is $a$-ergodic for all $a$ in the Weyl chamber $\w^P\subseteq A$ associated to $P$. In particular, $\mu$ is $A$-ergodic as well. 
\end{remark}


\begin{remark}\label{transverse lyapunov P-ergodic measures}
We get from Definition \ref{definition transverse Lyapunov A-ergodic} that there is a corresponding transverse Lyapunov exponent $\lo\in A^\ast$.    
\end{remark}

\subsubsection{Invariance of ergodic components of $\mu$}

Take any non-trivial $s\in\overline{\w^P}$, and let $\mu=\int\nu_xd\mu(x)$ be the $s$-ergodic decomposition of $\mu$. Recall we denote by $N_s$ the unstable subgroup of $s$ in $G$, see Section \ref{Stable and Unstable subgroups section}. Note that $N_s\leq P$, and thus $\mu$ is $N_s$-invariant.

The following fact uses the Hopf argument, and is inspired by \cite[Proposition~8.4]{dh}.

\begin{lemma}\label{invariance of ergodic components}
    For $\mu$-a.e. $x$ we have that $\nu_x$ is $N_s$-invariant.
\end{lemma}
\begin{proof}
    
    For any measurable partition $\xi$ of $M$ subordinate to the $N_s$-orbits, with respect to $\mu$, the conditional measure of $\mu$ along $\xi(x)$ is equivalent to the Riemannian volume on the $N_s$-orbit $N_s\cdot x$, by Lemma \ref{invariance for subgroups with orbits in unstable}, for $\mu$-a.e. $x$, as $\mu$ is $N_s$-invariant. Because of the Hopf argument, see \cite[Proposition~8.12]{brown2019entropylyapunovexponentsrigidity}, the partition into unstable manifolds with respect to $s$ refines the ergodic decomposition of $\mu$ with respect to $s$, so in particular the partition of $M$ into $N_s$-orbits refines this ergodic decomposition. 
    
    Thus, for $\mu$-a.e. $x$, the conditional measure of $\nu_x$ along $\xi(y)$ is equal to the conditional measure of $\mu$ for this element, for $\nu_x$-a.e. $y$. Hence, applying Lemma \ref{invariance for subgroups with orbits in unstable}, it follows that $\nu_x$ is $N_s$-invariant for $\mu$-a.e. $x$.
\end{proof}

\subsection{Negativity of transversal Lyapunov exponent on $\w^P$} 

\subsubsection{Minimal codimension}\label{Minimal  codimension}
We make use of the notion of minimal resonant codimension of $G$ defined in \cite{brhw}.

\begin{definition}\label{minimal resonant codimension definition}
    The \emph{minimal codimension} of $G$, which we denote by $c(G)$, is the minimal dimension of the manifold $G/Q$, for all proper parabolic subgroups $Q\leq G$.
\end{definition}

Now we follow \cite[Section~1]{brhw}.

Two roots $\alpha,\beta\in\Sigma_{(G,A)}$ are resonant, if there exists $c>0$ such that $\alpha=c\beta$. Denote by $\hat \Sigma_{(G,A)}$ the corresponding set of equivalence classes. For $\xi\in\hat \Sigma_{(G,A)}$, consider the corresponding Lie algebra $\mathfrak g^\xi=\bigoplus_{\beta\in\xi}\mathfrak g^\beta$, and its corresponding analytic group in $G$ by $G^\xi$. Given any $\beta\in\Sigma_{(G,A)}$, denote its equivalence class in $\hat\Sigma_{(G,A)}$ by $\hat \beta$.

Given a parabolic subalgebra $\mathfrak q$ of $\lie(G)$ containing $\operatorname{Lie}(P)$, the \emph{minimal resonant codimension} of $\mathfrak q$ is the cardinality of the set
\begin{equation*}
    \{\xi\in \hat \Sigma_{(G,A)}:\mathfrak g^\xi\not\subset\mathfrak q\}.
\end{equation*}

\begin{definition}
    We define the \emph{minimal resonant codimension} of $\lie(G)$, denoted by $r(G)$, to be the minimal resonant codimension of $\mathfrak q$ as $\mathfrak q$ varies among the set of maximal parabolic subalgebras of $\lie(G)$ containing $\lie(P)$.
\end{definition}

\begin{lemma}
    We have $c(G)\geq 2$. Furthermore, $c(G)=2$ if and only if $G$ is locally isomorphic to $\SL(3,\R)$. 
\end{lemma}

\begin{proof}

    If $r(G)$ denotes the resonant codimension of $G$, defined in \cite[Defintion~1.2]{brhw}, then it is not hard to see that $r(G)\leq c(G)$. Since $G$ is simple and higher rank, we have $r(G)\geq 2$, see \cite[Example~1.3]{brhw}, and thus $c(G)\geq 2$. Hence $c(G)=2$ only if $r(G)=2$, and the latter happens only when $G$ is locally isomorphic to $\SL(3,F)$ for $F=\R, \C, \mathbb H$, see \cite[Example~1.3]{brhw}, but $c(G)=2$ only in the $F=\R$ case.
\end{proof}

\subsubsection{Kernel of transverse Lyapunov exponent of $\mu$ does not intersect $\w^P$}

The following lemma  is inspired by ideas appearing in \cite{bfh} and \cite{brhw} to construct $G$-invariant measures on induced spaces of actions of lattices on manifolds. Recall  the definition of $\lo$ for $\mu$, see Definition \ref{definition transverse Lyapunov A-ergodic}. This argument is inspired by one of the arguments given in the recent progress on Zimmer's conjecture by Sebastian Hurtado for cocompact lattices, see the proof of \cite[Theorem~15.1]{brown2019entropylyapunovexponentsrigidity}.

\begin{lemma}\label{kernel of transverse exponent touches not weyl chamber}
     $\ker(\lambda_\mu^T)\cap\w^P=\{0\}$.
\end{lemma}
\begin{proof}
    Suppose, to get a contradiction, that there exists some $a\in \ker(\lambda_\mu^T)\cap\w^P$. Let us show $\mu$ is $G^{-\beta}$ invariant for all $\beta\in\Sigma_P$, which implies $\mu$ is also invariant under the parabolic $\Pop$ opposite to $P$, which implies $\mu$ is $G$-invariant as it is $P$-invariant as well, yielding a contradiction to Assumption \ref{assumption no invariant probability}.

    As $G\acts M$ has codimension one orbits, and because of our choice of $a$, we get that the Lyapunov exponents of $a$ with respect to $\mu$ which are positive or negative at $a$ are only roots of $(G,A)$, and hence 
    \begin{equation*}
        \sum_{\beta\in\Sigma_P}m_\beta\beta(a)\geq h_{\mu}(a^{-1})=h_{\mu}(a)=\sum_{\beta\in\Sigma_P}m_\beta\beta(a),
    \end{equation*}
    where $m_\beta$ is the \emph{multiplicity} of $\beta$, i.e., the dimension of the root group $G^\beta$; the first inequality follows from Theorem \ref{Ruelle-Margulis inequality}; and the last equality follows from Lemma \ref{invariance from maximal entropy} and $P$-invariance of $\mu$. Hence, $h_{\mu}(a^{-1})=\sum_{\beta\in\Sigma_P}m_\beta\beta(a)$, and therefore $\mu$ is invariant under $G^{-\beta}$ for all $\beta\in\Sigma_P$ applying Lemma \ref{invariance from maximal entropy}, which contradicts Assumption \ref{assumption no invariant probability}. 
    \end{proof}

    \begin{definition}\label{defintion Q0}
Let $Q_0$ be the stabilizer of $\mu$ in $G$, i.e., the subgroup $\{g\in G:g_\ast\mu=\mu\}$.    
\end{definition}

$Q_0$ is a proper parabolic subgroup of $G$, because of Assumption \ref{assumption no invariant probability}, and as it  contains $P$. Denote by $\Sigma_{Q_0}$ the set of all nonzero roots $\beta\in A^\ast$ such that the associated unipotent group $G^\beta$ is contained in $Q_0$ in $G$. The following argument is based on some invariance arguments appearing in the work of Brown, Rodr\'iguez-Hertz and Wang, see \cite{brhw}, using the resonant codimension invariant $r(G)$.

For the following argument, we use the higher rank hypothesis on $G$. The argument is based on some of the ideas behind the proof of the ``non-resonance implies invariance" philosophy used in \cite{brhw}, see for instance Proposition 5.1 in that paper.

    \begin{lemma}\label{kernel of transversal no faces}
        For each $\alpha\in\Delta_P$ we have $\ker(\lambda_\mu^T)\cap\operatorname{Int}(\ker(\alpha)\cap\overline{\w^P})=\emptyset$.
    \end{lemma}

    \begin{proof}
    Suppose, to get a contradiction, that there is $\alpha\in\Delta_P$ along with some $s_0\in\operatorname{Int}(\ker(\alpha)\cap\overline{\w^P})$ such that $\lo(s)=0$.
    We have
    \begin{equation*} h_\mu(s_0)=\sum_{\beta\in\Sigma_P\setminus\hat\alpha}m_\beta\beta(s_0)\geq h_\mu(s_0^{-1})
    \end{equation*}
    where the equality holds as $\mu$ is $P$-invariant and $P$-ergodic, and $\lo$ is the only Lyapunov exponent of the ergodic action $A\acts (M,\mu)$ which is not a root; and the inequality follows from Theorem \ref{Ruelle-Margulis inequality}. Therefore, 
    \begin{equation*}       h_\mu(s_0^{-1})=\sum_{\beta\in\Sigma_P\setminus\hat\alpha}m_\beta\beta(s_0),
    \end{equation*}
    this proves that $\mu$ is $G^{\beta}$-invariant for all $\beta\in\Sigma_{(G,A)
    }\setminus\hat\alpha$ as this measure is $P$-invariant as well, hence, as $r(G)\geq 2$, since $G$ is simple and higher rank, see \cite[Example~1.3]{brhw}, we conclude that $\mu$ is $G$-invariant, which contradicts Assumption \ref{assumption no invariant probability}.
    \end{proof}

\begin{corollary}\label{resonant root can't be in root complement of Q0}
$\lambda_\mu^T$ cannot be resonant to a root in $\Sigma_{Q_0}$
\end{corollary}

\begin{proof}
    Suppose $\lambda_\mu^T$ is resonant to some $\alpha\in\Sigma_{Q_0}$. There exists some parabolic $P'$ contained in $Q_0$ such that $\alpha\in\Sigma_{P'}$, this is not hard to check by considering the Levi decomposition of $Q_0$. 
    
    Since $\mu$ is $A$-ergodic and $P'\leq Q_0$, by Remark \ref{mu is A-ergodic}, we get that $\mu$ is $P'$-invariant and $P'$-ergodic, so we can apply Lemma \ref{kernel of transversal no faces} and obtain a contradiction, as $\lo$ is being assumed to be resonant to a root in $\Delta_{P'}$. 
\end{proof}


\subsection{The subgroup $Q_1$}

\subsubsection{The High Entropy Method}
Let $G$ be a semisimple Lie group, $A\leq G$ be a maximal subtorus.

For an $A$-invariant probability measure $\mu$ on $M$, and a root $\beta$ of $(G,A)$, we denote by $\mu^{G^\beta}_x$ the corresponding leafwise measures with respect to the partition by $G^\beta$ orbits, see \cite[Section~6]{el}. We have the following immediate consequence of \cite[Theorem~8.5]{ek}.

\begin{theorem}\label{high entropy method}
    Suppose we have a locally free action $G\acts X$, $X$ a locally compact, second countable, metric space. Let $\mu$ be an $A$-invariant and ergodic measure on $X$. Let $\alpha,\beta\in\Sigma_{(G,A)}$ be such that they are not proportional to each other. Then for $\mu$-a.e. $x$, $\mu$ is invariant under all elements of the bracket subgroup $[\supp(\mu_x^{G^\alpha}),\supp(\mu_x^{G^\beta})]$.
\end{theorem}

\subsubsection{The subgroup $Q_1$}

\begin{definition}\label{definition Q1}
Let $Q_1$ be the subgroup of $G$ generated by $A$ along with all root groups $G^\alpha$, for which $[\mu]_x^{G^\alpha}$ is not atomic for $\mu$-a.e. $x$. 
\end{definition}

\begin{remark}\label{Q0 subgroup of Q1}
    Notice that $P\leq Q_0\leq Q_1$, by Lemma \ref{haarconditionals}. See also \cite[Section~5.3]{brhw}.     
\end{remark}

\begin{lemma}\label{nonatomics Q1}
$\Sigma_{Q_1}$ equals the set of roots $\alpha\in\Sigma_{(G,A)}$ such that $[\mu]_x^{G^\alpha}$ is non-atomic for $\mu$-a.e. $x$. 
\end{lemma}

\begin{proof}
Any $\alpha\in\Sigma_{Q_1}$ equals a sum $\alpha_1+\cdots+\alpha_r$, $r\geq 1$, for $\alpha_j\in\Sigma_{(G,A)}$ such that for every $1\leq j\leq r$, $\alpha_1+\cdots+\alpha_j\in\Sigma_{(G,A)}\setminus\{0\}$ and $[\mu]_x^{G^{\alpha_j}}$ is not atomic for $\mu$-a.e. $x$ for every $1\leq j\leq r$. Hence, by applying Proposition \ref{high entropy method} iteratively, we obtain that $[\mu]_x^{G^\alpha}$ is not atomic for $\mu$-a.e. $x$.
\end{proof}

For the proof of Proposition \ref{Strict Maximal Inequality of transverse Lyapunov exponent} we need the following observation.  

\begin{lemma}\label{Q1 is proper}
     $Q_1$ is a proper parabolic subgroup of $G$. 
\end{lemma}
\begin{proof}
    Suppose otherwise. Then for every simple roots $\alpha,\beta\in\Delta_P$,  such that $\alpha+\beta$ is a root of $(G,A)$ as well, we have that $\mu$ is $G^{-\alpha-\beta}$ invariant as $-\alpha,-\beta\in\Sigma_{Q_1}$ by assumption and because of Theorem \ref{high entropy method}. Hence, as $\mu$ is both $G^\alpha$ and $G^\beta$ invariant, we get that $\mu$ is both $G^{-\alpha}$ and $G^{-\beta}$ invariant. Therefore, since for any $\alpha\in\Delta_P$ there is such a $\beta$(as Dynkin diagrams of simple Lie algebras are connected, and as $G$ is simple and higher rank), it follows that $\mu$ is $G^{-\alpha}$ invariant for every $\alpha\in\Delta_P$.

    Hence, since the root groups $\{G^{\pm\alpha}:\alpha\in\Delta_P\}$ generate $G$, and $\mu$ is $P$-invariant, this implies $\mu$ is $G$-invariant, contradicting Assumption \ref{assumption no invariant probability}.
\end{proof}

\subsection{Negativity of Transverse Lyapunov exponent on Weyl chamber of $P$}

Given any $\xi\in\tilde\Sigma(G,A)$, we denote by  $\mathscr W^\xi$ the partition of $M$ by the $G^\xi$ orbits, with respect to our action $G\acts M$.

If $\chi$ is the coarse Lyapunov exponent of the action $A\acts (M,\mu)$, recall we denote by $\mathscr W^\chi$ the corresponding partition by coarse Pesin manifolds with respect to $\chi$, see Section \ref{Coarse Lyapunov exponents and Pesin manifolds for higher rank actions}.

In the next lemma we prove that the transversal Lyapunov exponent is negative in the interior of the Weyl chamber corresponding to $P$.

\begin{lemma}\label{negativity of transversal lyapunov exponent on positive weyl chambers}
    For any $a\in\w^P$ we have $\lambda_{\mu}^T(a)<0$.
\end{lemma}
\begin{proof}  
    Suppose, to get a contradiction, that $\lambda_\mu^T(a)\geq 0$ for some $a\in\w^P$. 
    
    Let us first prove that
    \begin{equation}\label{mbjuhe}
        h_\mu(a)\leq \sum_{\beta\in\Sigma_P}m_\beta\beta(a).
    \end{equation}

    If $\lo$ is not resonant to a root in $\Sigma_{(G,A)}$, then (\ref{mbjuhe}) follows immediately from Theorem \ref{entropy product structure}, and from the $P$-invariance and $a$-ergodicity of $\mu$, see Remark \ref{mu is A-ergodic}.

    If $\lo$ is resonant to a root in $\Sigma$, because of our codimension one assumption on $G\acts M$, it follows that there exists a unique $\xi\in\hat\Sigma_P$ such that $\chi=\xi\cup\{\lo\}$. We then have from the $P$-invariance and $a$-ergodicity of $\mu$ that
    \begin{equation*}                h_\mu(a)=\sum_{\beta\in\Sigma_P\setminus\xi}m_\beta\beta(a)+h_\mu(a\mid\mathscr W^\chi);
    \end{equation*}
    however, from Theorem \ref{ledrappieryoung2} and the $P$-invariance of $\mu$ we have
    \begin{equation*}
        h_\mu(a\mid\mathscr W^\chi)\geq \sum_{\beta\in\xi}m_\beta\beta(a),
    \end{equation*}
    and hence we obtain (\ref{mbjuhe}) in this case as well. 

    However, since all Lyapunov exponents of $A\acts (M,\mu)$ which are positive at $a^{-1}$ are exactly $\{-\beta:\beta\in\Sigma_P\}$, from Theorem \ref{Ruelle-Margulis inequality} we get
    \begin{equation*}
        h_\mu(a^{-1})\leq \sum_{\beta\in\Sigma_P}m_\beta\beta(a);
    \end{equation*}
    hence, combining this with (\ref{mbjuhe}), along with $h_\mu(a)=h_\mu(a^{-1})$, it follows that $\mu$ has maximal entropy with respect to both $a$ and $a^{-1}$, which implies that $\mu$ is $G$-invariant, contradicting Assumption \ref{assumption no invariant probability}.
    \end{proof}

\begin{corollary}\label{negativity on faces of weyl chamber}
    For any $\alpha\in\Delta_P$ we have that $\lo(s)<0$ for all $s\in\operatorname{Int}(\ker(\alpha)\cap\overline{\w^P})$
\end{corollary}
\begin{proof}
    This is an immediate consequence of Lemmas \ref{kernel of transversal no faces} and \ref{negativity of transversal lyapunov exponent on positive weyl chambers}.
\end{proof}


\subsection{Maximal and Strict Maximal Inequalities of transversal Lyapunov exponent}

Let $\mu$ be a $P$-invariant and ergodic probability measure on $M$.

\begin{proposition}[Maximal Inequality of transverse Lyapunov exponent]\label{Maximal negativity of transverse Lyapunov exponent}
    For any  $a\in\w^P$ we have
    \begin{equation}
        \beta(a)<\lambda_\mu^T(a^{-1})
    \end{equation}
    for all $\beta\in\Sigma_P$
\end{proposition}

\begin{example}
    We need $G$ to have higher rank simple factors for Proposition \ref{Maximal negativity of transverse Lyapunov exponent}  to hold. Indeed, if $G=\SL(2,\R)$, and $\Gamma$ is a torsion free uniform lattice of $G$, and $P$ is the subgroup of upper triangular matrices of $G$, consider the canonical diagonal action $G\acts M:=G/\Gamma\times G/P$, which is smooth, locally free, and with codimension one orbits. Then, $P\acts M$ is uniquely ergodic, with corresponding invariant probability measure, $\mu:=\haar_{G/\Gamma}\times\delta_{eP}$, where $eP$ is the identity coset of $G/P$. The transverse Lyapunov exponent of $\mu$ is the zero functional in this case.
\end{example}

When the minimal codimension $c(G)$ is at least $3$, we have, in fact,  the following sharper version of Proposition \ref{Maximal negativity of transverse Lyapunov exponent}.

\begin{proposition}[Strict Maximal Inequality of transverse Lyapunov exponent]\label{Strict Maximal Inequality of transverse Lyapunov exponent}
    Suppose $c(G)\geq 3$. Let  $\alpha\in\Delta_P$. Then for any $s\in\operatorname{Int}(\ker(\alpha)\cap\overline{\w^P})$ we have
    \begin{equation}
        \beta(s)<\lambda_\mu^T(s^{-1})
    \end{equation}
    for all $\beta\in\Sigma_P$.
\end{proposition}

\begin{remark}\label{locally isometric to SL(3,R) when Strict Inequality fails}
When $c(G)\geq 3$, Proposition \ref{Maximal negativity of transverse Lyapunov exponent} follows immediately from Proposition \ref{Strict Maximal Inequality of transverse Lyapunov exponent} as $\w^P$ is contained in the convex hull of the codimension one faces of $\overline{\w^P}$. For the case $c(G)=2$ Proposition \ref{Strict Maximal Inequality of transverse Lyapunov exponent}, a priori, might not necessarily hold, this proposition is the reason why we need a separate proof of Theorem \ref{main} when $c(G)=2$.
\end{remark}

\begin{remark}
    We prove Proposition \ref{Maximal negativity of transverse Lyapunov exponent} in the case $c(G)=2$ in Section \ref{unique ergodicity sl3 section}, see Corollary \ref{maximal negativity case sl3r}.
\end{remark}

\subsection{Proof of Strict Maximal Inequality}

Denote by $\beta_P^H$ the \emph{highest root} with respect to $P$, i.e., $\beta_P^H$ is greater than or equal to all elements in $\Sigma_P^+$, that is, for any $\beta\in\Sigma_P$, $\beta_P^H-\beta$ is a sum of elements in $\Sigma_P$, with non-negative integer coefficients.

\begin{remark}\label{highest root can't be in Q1}
    We have that $-\beta_P^H\notin\Sigma_{Q_1}$. Otherwise, since for any $\alpha\in\Sigma_P$ we have that $\beta_P^H-\alpha$ is a sum of roots in $\Delta_P$, then $-\alpha$ is obtained from $-\beta_P^H$ by adding to it this sum, and as $\mu$ is $P$-invariant, it follows from using the bracket relation $[G^{-\beta_P^H},G^{\beta_P^H-\alpha}]=G^{-\alpha}$ that $\mu$ is $G^{-\alpha}$-invariant. Since $\alpha\in\Sigma_P$ was arbitrary, this contradicts Assumption \ref{assumption no invariant probability}. 
\end{remark}

Denote by $\Pop$ the minimal parabolic of $G$ opposite to $P$, so its set of simple positive roots is denoted by $\Delta_{\Pop}$. Recall that $Q_0=Q_0(\mu)$ is the stabilizer of $\mu$ in $G$, and that the construction of the subgroup $Q_1=Q_1(\mu)$ of $G$ is given in Defintion \ref{definition Q1}.

\begin{lemma}\label{Q1 and Q0 difference}
    We have that $\Sigma_{Q_1}\setminus\Sigma_{Q_0}\subset\Delta_{\Pop}$. 
\end{lemma}
\begin{proof}
    Let $\alpha\in\Sigma_{\Pop}\setminus\Delta_{\Pop}$. Suppose that $\alpha\in\Sigma_{Q_1}$. Then all elements in $\Delta_{\Pop}$ making $\alpha$ up belong to $\Sigma_{Q_0}$ as $\mu$ is $P$-invariant and by applying Theorem \ref{high entropy method} and using a similar argument as the one given in the previous remark; this in turn implies $\mu$ is $G^\alpha$-invariant, i.e., $\alpha\in\Sigma_{Q_0}$.
\end{proof}

\begin{corollary}\label{lambda not resonant to root in Q1}
    Suppose $\lo$ is resonant to some $\beta\in\Delta_{(G,A)}$. Then $\beta\notin \Sigma_{Q_1}$.
\end{corollary}

\begin{proof}
    Assume, to get a contradiction, that $\beta\in\Sigma_{Q_1}$. Since $\beta\notin\Sigma_{Q_0}$ by Corollary \ref{resonant root can't be in root complement of Q0}, we get that $-\beta\in\Delta_P$ because of Lemma \ref{Q1 and Q0 difference}; but then $\ker(\lo)=\ker(-\beta)$, contradicting Lemma \ref{kernel of transversal no faces}. 
\end{proof}

\begin{lemma}\label{no full entropy along coarse Pesin manifold}
    Let $\chi$ be the coarse Lyapunov exponent of $\lo$ with respect to $\mu$. For all $s\in\overline{\w^P}\setminus\{e_G\}$,
    \begin{equation*}h_\mu(s^{-1}\mid\w^\chi)\leq\lambda_\mu^T(s^{-1}).
    \end{equation*}
\end{lemma} 
\begin{proof} 
    If $\lo$ is not resonant to a root in $\Sigma_{(G,A)}$, then $\chi=\{\lo\}$, and $\w^\chi(x)$ is one-dimensional for $\mu$-a.e. $x$.  Then, in this case, our lemma follows immediately from Theorem \ref{ledrappier young conditional entropy}. 
    
    Assume $\lo$ is resonant to an element in $\Sigma_{(G,A)}$. Then there is one coarse root $\xi$ of the root system $\Sigma_{(G,A)}$, all of whose elements are resonant to $\lo$, denote by $\chi$ the corresponding coarse Lyapunov for $\lo$, so that $\chi=\xi\cup\{\lo\}$. 
    
    Then, for $\mu$-a.e. $x$, $\w^\chi(x)$ is a $(m_{\hat{\xi}}+1)$-dimensional injectively immersed submanifold of $M$ containing $x$ and our given action $G\acts M$ restricts to an action $G^{\xi}\acts \w^\chi(x)$. Hence, there exists a diffeomorphism $\phi:(-1,1)^{n_{\xi}+1}\rightarrow W$, where $W$ is an open neighborhood of $x$ in $\w^\chi(x)$(immersed topology), such that:
    \begin{enumerate}[label=(\roman*)]
        \item $\phi(0,0)=x$;
        \item $\phi$ restricted to $[-1,1]^{n_{\xi}}\times\{t\}$ is a diffeomorphism onto the connected component of  $\left[G^{\xi}\cdot\phi(0,t)\right]\cap W$ containing $\phi(0,t)$.
    \end{enumerate}

    Consider the disintegration, $\{\mu_x^\chi\}_{x\in M}$ of $\mu$ along a measurable partition subordinate to the partition $\mathscr W^\chi$, see \cite[Definition~8.1]{brown2019entropylyapunovexponentsrigidity}, normalized so that they all are probability measures.

    Set $\mu^I:=(\pi_2\circ \phi^{-1})_\ast[\mu^\chi_x\upharpoonright W]$, where $\pi_2:[0,1]^{n_{\xi}+1}\rightarrow [0,1]$ is the canonical projection map onto the second coordinate, and $\mu\upharpoonright W$ is the restriction of the measure $\mu$ to $W$. 

    Let $\{\mu_t\}_{t\in (-1,1)}$ be the Rokhlin disintegration of $\mu_x^\chi\upharpoonright W$ with respect to the factor map $\pi_2\circ\phi^{-1}:(W,\mu_x^\chi\upharpoonright W)\rightarrow((-1,1),\mu^I)$. 
    
    By construction of leafwise measures, see Section \ref{Invariance from entropy considerations}, we get that $\mu_t$ is a  measure with a single atom for $\mu^I$-a.e. $t$, as these measures are, $\mu^I$ almost surely, the  pushforward of $\mu_{\phi(0,t)}^{G^{\xi}}$(normalized to be a probability measure) under the map $u\in G^{\xi}\mapsto u\cdot\phi(0,t)$, and $\mu_x^{G^{\xi}}$ is an atom for $\mu$-a.e. $x$, since $\xi\cap\Sigma_{Q_1}=\emptyset$, which follows by applying Lemma \ref{nonatomics Q1} along with \cite[Lemma~8.4, (1),(3)]{brown2019entropylyapunovexponentsrigidity}.

    Thus, for any small enough $\delta>0$, we have
    \begin{equation*}
    \begin{split}
        \mu^\chi_x(B_{\w^\chi(x)}(x,\delta))&=\int_{B_I(0,\delta)} \mu_t(B_{\w^\chi(x)}(x,\delta))d\mu^I(t)\\&=\mu^I(B_I(0,\delta)).
    \end{split}
    \end{equation*}

    Hence, for $\mu$-a.e. $x$ we have that
    \begin{equation*}
        \limsup_{\delta\rightarrow 0}\frac{\log(\mu^\chi_x(B_{\w^\chi(x)}(x,\delta)))}{\log\delta}=\limsup_{\delta\rightarrow 0}\frac{\log(\mu^I(B_I(0,\delta)))}{\log\delta}\leq 1,
    \end{equation*}
    where the last inequality follows from $I$ being a one-dimensional manifold. Since the partition $\mathscr W^\chi$ is $s$-invariant, for any $s\in\overline{\w^P}\setminus\{0\}$, we get that for any measurable partition $\eta$ subordinated to the partition $\mathscr W^\chi$, with respect to $\mu$, the pointwise dimension of $\mu$ relative to $x$, see (\ref{pointwise dimensions}), is at most $1$ for $\mu$-a.e. $x$.

    Therefore, since $\lo(s)<0$ for all $s\in\overline{\w^P}\setminus\{e_G\}$ because of Lemmas \ref{kernel of transversal no faces} and \ref{negativity of transversal lyapunov exponent on positive weyl chambers}, we can apply Theorem \ref{ledrappieryoung2} and obtain the inequality of our lemma. 
\end{proof}

\begin{Maximal negativity of transverse Lyapunov exponent when r(G)=2}
Take any $a\in\w^P$. 

Since $\lo$ is not resonant to any root in $\Sigma_{(G,A)}$, we have by Theorem \ref{entropy product structure}, Corollary \ref{lambda not resonant to root in Q1}, and Lemma \ref{non resonance when failure of strict} that
\begin{equation*}
h_\mu(a)=\sum_{\beta\in\Sigma_P}m_\beta\beta(a)\text{, } h_\mu(a^{-1})=h_\mu(a^{-1}\mid\w^{\lo})+\underset{-\beta\in\Sigma_{Q_1}}{\sum_{\beta\in\Sigma_P}}h_\mu(a^{-1}\mid\w^{-\beta}),
\end{equation*}
which implies 
\begin{equation}\label{rfghhj}
\begin{split}
    \lo(a^{-1})\geq h_\mu(a^{-1}\mid\w^{\lo})&=\sum_{\beta\in\Sigma_P}m_\beta\beta(a)-\underset{-\beta\in\Sigma_{Q_1}}{\sum_{\beta\in\Sigma_P}}h_\mu(a^{-1}\mid\w^{-\beta})\\   &\geq\underset{-\beta\notin\Sigma_{Q_1}}{\sum_{\beta\in\Sigma_P}}m_\beta\beta(a)>\beta(a),
\end{split}
\end{equation}
where the first inequality follows from Theorem \ref{Ruelle-Margulis inequality}; the second inequality follows from Lemma \ref{zero entropy atomic leafwise}; and the last inequality follows from the fact that $m_{\beta_P^H}\beta_P^H(a)$ is a summand of the sum appearing in the second line of (\ref{rfghhj}), and in this case we have a strict inequality since $r(G)\geq 2$, as $G$ is higher rank, there is some $\beta\in\Sigma_P\setminus\{\beta_P^H\}$ such that $-\beta\notin\Sigma_{Q_1}$.  
\end{Maximal negativity of transverse Lyapunov exponent when r(G)=2}




We now prove Proposition \ref{Strict Maximal Inequality of transverse Lyapunov exponent}.
\s

\begin{Strict Maximal Inequality of transverse Lyapunov exponent for main theorem}
    Let $\alpha\in\Delta_P$ and $s\in\operatorname{Int}(\ker(\alpha)\cap\overline{\w^P})$. We have $\lambda_\mu^T(s)<0$ by Lemmas \ref{kernel of transversal no faces} and \ref{negativity of transversal lyapunov exponent on positive weyl chambers}.

Denote by $\chi$ the coarse Lyapunov exponent of $\lambda_\mu^T$.

We have by Theorem \ref{entropy product structure} and Corollary \ref{lambda not resonant to root in Q1} that
\begin{equation*}
h_\mu(s)=\sum_{\beta\in\Sigma_P}m_\beta\beta(s)\text{, } h_\mu(s^{-1})=h_\mu(s^{-1}\mid\w^\chi)+\underset{-\beta\in\Sigma_{Q_1}}{\sum_{\beta\in\Sigma_P}}h_\mu(s^{-1}\mid\w^{-\beta});
\end{equation*}
hence, we obtain from $h_\mu(s)=h_\mu(s^{-1})$ and Lemma \ref{nonatomics Q1} that
\begin{equation}\label{useful inequality for the future}
h_\mu(s^{-1}\mid\w^\chi)=\sum_{\beta\in\Sigma_P}m_\beta\beta(s)-\underset{-\beta\in\Sigma_{Q_1}}{\sum_{\beta\in\Sigma_P}}h_\mu(s^{-1}\mid\w^{-\beta})\geq \underset{-\beta\notin\Sigma_{Q_1}}{\sum_{\beta\in\Sigma_P}}m_\beta\beta(s),
\end{equation}
since $h_\mu(s^{-1}\mid\w^{-\beta})\leq m_\beta\beta(s)$ for all $\beta\in\Sigma_P$, and by Lemmma \ref{zero entropy atomic leafwise}. As $\lo(s^{-1})\geq h_\mu(s\mid\w^\chi)$ because of Lemma \ref{no full entropy along coarse Pesin manifold}, it follows that 
\begin{equation}\label{mmmn}
    \lo(s^{-1})\geq\underset{-\beta\notin\Sigma_{Q_1}}{\sum_{\beta\in\Sigma_P}}m_\beta\beta(s); 
\end{equation}
however, by Remark \ref{highest root can't be in Q1}, we get that $\beta_P^H(s)$ is a summand of the right hand side of the inequality above; but the right hand sum in the inequality above contains at least one summand different to $\alpha(s)$ and $\beta_P^H(s)$, as $c(G)\geq 3$. We conclude from the previous inequality that
\begin{equation*}
\lo(s^{-1})>m_{\beta_P^H}\beta_P^H(s)\geq\beta_P^H(s),
\end{equation*}
 which proves the assertion.
\end{Strict Maximal Inequality of transverse Lyapunov exponent for main theorem}




\subsection{Transverse Pesin manifolds}\label{transverse pesin manifolds section}

\begin{definition}\label{definition transverse Pesin manifolds}
 By Propositions \ref{Maximal negativity of transverse Lyapunov exponent} we have that $\lambda_\mu^T(a^{-1})>\beta(a)$ for all $a\in\w^P$ and all $\beta\in\Sigma_P$. Since the set of Lyapunov exponents of the measure $\mu$ is $\Sigma_{(G,A)}\cup\{\lo\}$ and $\dim(M)=\dim(G)+1$, this implies that for $\mu$-a.e. $x$ there exists a one-dimensional global Pesin manifold associated to the Lyapunov exponent $\lambda_\mu^T$, which equals the first unstable Pesin manifold  with respect to $a^{-1}$, $\w^{u,1}(x;a^{-1})$ for all $a\in\w^P$, see Definition \ref{definition ith unstable}. By Proposition \ref{Maximal negativity of transverse Lyapunov exponent} and Lemma \ref{uniqueness of ith unstable higher rank}, $\w^{u,1}(x;a^{-1})$ does not depend on the particular choice of $a\in\w^P$, hence we denote this injectively immersed submanifold by $\w^T(x)$, and call it \emph{transverse Pesin manifold at $x$}.
\end{definition}

\begin{notation}
    A priori, it could be the case that for a point $x\in M$ we could have two different transverse Pesin manifolds at $x$ coming from $P$-invariant and $P$-ergodic probability measure different to $\mu$, however, because of Proposition \ref{Maximal negativity of transverse Lyapunov exponent}, we have for any $a\in\w^T$ and $\mu^P$-a.e. $x$ that 
    \begin{equation*}
    \begin{split}
        \w^T(x)&=\w^{u,1}(x;a^{-1})\\&=\left\{ y \in M \,\middle|\, \limsup_{n \to \infty} \frac{1}{k} \log d\big(a^{n}x, a^ny)\big) <-\beta_P^H(a) \right\},
    \end{split}
    \end{equation*}
    so $\w^T_\mu(x)$ does not depend on the particular probability measure we use. 
\end{notation}

\begin{remark}\label{why name transverse}
The name in the previous definition comes from the observation that $\w^T_\mu(x)$ is transversal to the $G$-orbits of the action $G\acts M$ at $x$, for $\mu$-a.e. $x$, as $\w^{u,1}(x;a^{-1})$ is tangent to the one dimensional Osedelets distribution associated to $\lo$, and this distribution is transversal to $G\cdot x$ at $x$, since $\lo\notin\Sigma_{(G,A)}$.    
\end{remark}

\begin{remark}\label{first unstable for regular and not so singular elements}
    When $r(G)\geq 3$, by Proposition \ref{Strict Maximal Inequality of transverse Lyapunov exponent} we have, for any $\alpha\in\Delta_P$, any $s\in\operatorname{Int}(\ker(\alpha)\cap \overline{\w^P})$ and any $a\in\w^P$, that  the first unstable manifold $\w^{u,1}(x;s^{-1})$ exists and equals $\w^{u,1}(x;a^{-1})$, i.e. $\w^T(x)$, for $\mu$-a.e. $x$,  because of Lemma \ref{uniqueness of ith unstable higher rank}.
\end{remark}

\subsection{Essential $P$-equivariance of Transverse Pesin Manifolds}

Under Assumption \ref{assumption no invariant probability}, we have the following equivariance property of transversal Pesin manifolds, with respect to the action of $P$.

\begin{proposition}\label{P-equivariance of transverse Pesin manifolds}
    Suppose $r(G)\geq 3$. For $\mu$-a.e. $x$ and $\haar_P$-a.e. $p$ we have 
    \begin{equation*}
        p\cdot \w^T(x)=\w^T(px).
    \end{equation*}
\end{proposition}

By Proposition \ref{Maximal negativity of transverse Lyapunov exponent}, we can consider the partition of $M$ by transverse Pesin manifolds,  $\mathscr W^T$, to be a partition by strong stable manifolds with respect to $a^{-1}$; since $G\acts M$ has codimension one orbits, we get that we can consider the partition of $M$ into $P$ orbits as the partition into center stable manifolds with respect to $a^{-1}$; hence, the previous statement tells us that these strong stable and center unstable partitions jointly integrate. 

\begin{proof}
Take any $\alpha\in\Delta_P$. Fix any $a\in\w^P$ and $s_\alpha\in \operatorname{Int}(\ker(\alpha)\cap\overline{\w^P})$.

As $\mu$ is $G^\alpha$-invariant, we obtain from Proposition \ref{haarconditionals} that for $\mu$-a.e. $x$ and $\haar_{G^\alpha}$-a.e. $u$ the transverse Pesin manifolds $\w^T(x)$ and $\w^T(ux)$ exist; because of Remark \ref{first unstable for regular and not so singular elements} they equal the manifolds $\w^{u,1}(x;s_{\alpha}^{-1})$ and $\w^{u,1}(ux;s_\alpha^{-1})$, respectively; $u$ and $s_\alpha$ commute, hence for such $x$ and $u$ we have that 
    \begin{equation}\label{equivariance simple roots}
    u\cdot\w^T(x)=u\cdot\w^{u,1}(x;s_{\alpha})=\w^{u,1}(ux;s_{\alpha})=\w^T(ux)        
    \end{equation}
by definition of the first unstable manifold and of transverse Pesin manifold. Hence, we obtain the equality of the lemma for $\mu$-a.e. $x$ and $\haar_{G^\alpha}$-a.e. $u$.

For all roots $\alpha$ and $\beta$ in $\Sigma_P$ such that $\alpha+\beta\in\Sigma_P$, the map 
\begin{equation*}
    \begin{split}
        G^\alpha\times G^\beta&\rightarrow G^{\alpha+\beta}\\
        (u,v) &\mapsto [u,v]
    \end{split}
\end{equation*}
is  smooth and onto $G^{\alpha+\beta}$, and hence it pushes forward the product of the Lebesgue measures on $G^\alpha$ and $G^\beta$ to the Lebesgue measure on $G^{\alpha+\beta}$; hence this map pushes forward $\haar_{G^\alpha}\times \haar_{G^\alpha}$ to a measure in the class of $\haar_{G^{\alpha+\beta}}$. Hence, by using (\ref{equivariance simple roots}), the fact that all (non-zero) elements of $\Sigma_P$ are non-negative linear combinations of elements of $\Delta_P$, along with iterative applications of Fubini's theorem, we obtain that $u\cdot\w^T(x)=\w^T(ux)$ for $\mu$-a.e. $x$ and $\haar_{U(P)}$-a.e. $u$.

Finally, since $\w^T(x)$ is $C_K(A)$-equivariant for $\mu$-a.e. $x$, because $C_K(A)$ and commutes with every $a\in A$, we obtain from the Iwasawa decomposition of $P$ along with Fubini's theorem the statement of our proposition. 
\end{proof}

Not only we have essential equivariance of the transverse Pesin partition $\mathscr W^T$ with respect to the action of $P$, we also have essential equivariance for the same partition with respect to the whole action $G\acts M$. To prove this, we need the following observation.

\begin{lemma}\label{big cell}
   Let $G$ be a connected semisimple Lie group with finite center, $P\leq G$ a minimal parabolic subgroup, and $\tilde U$ the unipotent radical of the parabolic opposite to $P$. Let $X\subseteq G$ be a Borel subset. Then $X$ has full $\haar_G$ measure, if and only if, for $\haar_P$-a.e. $p$ and $\haar_{\tilde U}$-a.e. $u$ we have $up\in X$. 
\end{lemma}

\begin{proof}
    The ``big cell" lemma, see \cite[Theorem~21.84]{milne} and \cite[Lemma~5.1.4]{zimmerbook}, proves that the map $m:\tilde U\times P\rightarrow G$ given by $(u,p)\mapsto up$ is an open immersion, with dense image(see \cite[Theorem~21.73.b]{milne}). Thus, since Haar measures on $G$ are volume forms, as $G$ is a Lie group, any open dense subset of $G$ has full Haar measure, so $X\subseteq G$ has full $\haar_G$ measure if and only if its intersection with the image of $m$ has full Haar $\haar_G$ measure. 

    Furthermore, since the Haar measures of $U$ and $P$ are also given by volume forms, it follows that the measure induced by the push-forwarded form, via $m$, of the wedge of these two volume forms is equivalent to $\haar_G$ with a nowhere zero analytic density on the image of $m$, this yields the second statement of the claim.
\end{proof}

Suppose $f:N\rightarrow N$ is a diffeomorphism of a closed manifold $N$, and $\lambda>0$. For any $x\in N$, we consider the \emph{$\lambda$-unstable subset} of $M$ at $x$ with respect to $f$, to be the set given by

\begin{equation}\label{definition unstable set}
    \w^u_\lambda(x;f):=\left \{y\in N:\limsup_{n\rightarrow\infty}\frac1n\log d(f^{-n}(x),f^{-n}(y))\leq-\lambda\right \}.
\end{equation}

\begin{remark}\label{transverse pesin manifold and unstable of transverse Lyapunov}
    For any $a\in\w^P$, $\w_{\lo(a^{-1})}^u(x;a^{-1})$ coincides with $\w^T(x)$ for $\mu$-a.e. $x$, as for such $x$, the latter set coincides with the first unstable Pesin manifold at $x$ with respect to $a^{-1}$, by construction, see Definition \ref{definition transverse Pesin manifolds}.
\end{remark}

\begin{proposition}\label{essential G equivariance of transverse Pesin}
    For $\mu^P$-a.e. $x$ and $\haar_G$-a.e. $g$ we have $$g\cdot\w^T(x)=\w^u_{\lo(a^{-1})}(gx;a^{-1}),$$
    for any $a\in\w^P$.
\end{proposition}

\begin{proof}
    If $x$ is $\mu$ generic enough for Proposition \ref{P-equivariance of transverse Pesin manifolds} to hold, we get by Remark \ref{transverse pesin manifold and unstable of transverse Lyapunov} that for any $a\in\w^P$, any $u\in \tilde U$($\tilde U$ the unipotent radical of the opposite of $P$), and for $\haar_P$-a.e. $p$, that
    \begin{equation*}up\cdot\w^T(x)=u\cdot\w^T(px)=u\cdot\w^u_{\lo(a^{-1})}(px,a^{-1})=\w^u_{\lo(a^{-1})}(upx,a^{-1}).
    \end{equation*}
    The last equality follows from the fact that $d(a^nuy_1,a^nuy_2)\leq C_u\cdot d(a^ny_1,a^ny_2)$ for all $y_1,y_2\in M$ and all $n\geq 0$, where $C_u>0$ is given by
    \begin{equation*}
        C_u:=\sup_{n\geq0}\sup_{x\in M}||D_xu^{(n)}||
    \end{equation*}
    where $u^{(n)}:=a^nua^{-n}$ for all $n\geq 0$;  $C_u<\infty$ since $u^{(n)}\rightarrow e$ as $n\rightarrow\infty$ and $M$ is compact. 

    Hence, applying Lemma \ref{big cell}, we get the result.
    \end{proof}

\section{Unique ergodicity of $P\acts M$ when $r(c)>2$}\label{Unique ergodicity of P section}

In this section we work under Assumption \ref{assumption no invariant probability}. In this section we prove that the  $P\acts M$ is uniquely ergodic, when $c(G)\geq 3$. All the statements in this section hold in the $c(G)=2$ case, unless it is indicated we are working in the  $c(G)\geq 3$ case.

\begin{remark}\label{special case of cG=2 case}
    The assumption $c(G)\geq 3$ is used in this section to ensure that any $P$-invariant and $P$-ergodic probability measure satisfies the inequalities of Proposition \ref{Strict Maximal Inequality of transverse Lyapunov exponent}. However, if $c(G)=2$ and $G\acts M$ is such that any $P$-invariant and $P$-ergodic probability measure satisfies the inequalities of Proposition \ref{Strict Maximal Inequality of transverse Lyapunov exponent}, then the proof we give of unique ergodicity in the $c(G)\geq 3$ goes through in that case.
\end{remark}

\begin{proposition}\label{unique P ergodicity}
  If $c(G)\geq 3$, then $P\acts M$ is uniquely ergodic. 
\end{proposition}

\subsection{$P\acts M$ is uniquely ergodic}

Let $\mu$ be a $P$-ergodic, $P$-invariant Radon probability measure on $M$, where $G\acts M$ satisfies the hypothesis of Proposition \ref{unique P ergodicity}. By Definition \ref{definition transverse Pesin manifolds}, there exist transverse Pesin manifolds $\w^T(x)$ for $\mu$-a.e. $x$.

\begin{lemma}\label{transverse Pesin manifold is uniformly locally nice}
For $\mu$-a.e. $x$,  $\w^T_\mu(x)$ intersects  $G$-orbits transversely.
\end{lemma}
\begin{proof}
    
    Pick $a\in\w^P$, then $\mu$ is $a$-ergodic because of Remark \ref{mu is A-ergodic}. Fix $x\in M$ and suppose Birkhoff's ergodic theorem holds for $x$ with respect to $a$ and $\mu$.

    There exists $\Lambda\subseteq M$ of positive $\mu$ measure such that for some $r_0>0$ we have that for every $y\in\Lambda$ the following holds: 
    \begin{enumerate}
        \item $\w^T(y)$ exists;
        \item the closed ball of radius $r_0$ in $\w^T(y)$ centered at $y$ is an embedded submanifold that is everywhere transversal to the $G$-orbits.
    \end{enumerate}
    The latter property holds for some $r_0>0$ in a set of positive $\mu$ measure as for $\mu$-a.e. $y$, $\w^T(y)$ is transversal to $G\cdot y$ at $y$, see Remark \ref{why name transverse}

    Suppose, to get a contradiction, that $\w^T(x)$ intersects a $G$-orbit non-transversely. Let $z$ be a point of $\w^T(x)$ for which this submanifold is not trasnversal to $G\cdot z$ at $z$. As the Birkhoff ergodic theorem holds for $x$ with respect to $a$ and $\mu$, we get that there exists a large enough $N>0$ such that $a^Nx\in\Lambda$ and $d(a^Nx,a^Nz)<r_0$, the latter property holds since $a\in\w^P$ and because of Lemma \ref{negativity of transversal lyapunov exponent on positive weyl chambers}. Then the ball of radius $r_0$ of $\w^T(a^Nx;a)$ centered at $a^Nx$ intersects $G\cdot a^Nz$ non-transversely, which contradicts (2) above.    
\end{proof}

\begin{lemma}\label{Pattraction}
    If $c(G)\geq 3$, then for $\mu$-a.e. $x$ we have that 
    \begin{equation*}
        \lim_{N\rightarrow\infty}\frac{1}{N}\sum_{j=0}^{N-1}\delta_{a^jgy}=\mu,
    \end{equation*}
    in the weak-$\star$ topology, for any $a\in\w^P$, any $y\in\w^T(x)$ and $\haar_G$-a.e. $g\in G$. 
\end{lemma}

\begin{proof}
    Pick $x\in M$ such that Proposition \ref{P-equivariance of transverse Pesin manifolds} and Birkhoff's ergodic theorem to hold for $px$, with respect to $\mu$ and $a$, for $\haar_P$-a.e. $p$. Then, for any $y\in\w^{u,1}(x;a)$ we have $py\in p\cdot \w^T(x)=\w^T(px)\subseteq\w^s(px;a)$ for $\haar_P$-a.e. $p$, and thus 
    \begin{equation*}
        \lim_{N\rightarrow\infty}\frac{1}{N}\sum_{j=0}^{N-1}f(a^jpy)=\lim_{N\rightarrow\infty}\frac{1}{N}\sum_{j=0}^{N-1}f(a^jpx)=\int fd\mu,
    \end{equation*}
    for all $f\in C^0(M)$, as such functions are uniformly continuous since $M$ is compact. Hence, for all $u\in U(P^{\operatorname{op}})$ we have
    \begin{equation*}
        \lim_{N\rightarrow\infty}\frac{1}{N}\sum_{j=0}^{N-1}f(a^jupy)=\lim_{N\rightarrow\infty}\frac{1}{N}\sum_{j=0}^{N-1}f(a^jpy)=\int fd\mu
    \end{equation*}
    for $\haar_P$-a.e. $p$, as $U(P^{\operatorname{op}})$ is contained in the stable subgroup of $G$ with respect to $a$, see Section \ref{Stable and Unstable subgroups section}, since $a\in\w^P$.

    Since $\haar_G$-a.e. $g$ can be written uniquely in the form $up$, for all $u$ and $p$ in subsets of $\tilde U$ and $P$ respectively, of full, respective, Haar measures, by Lemma \ref{big cell}, we obtain our lemma. 
\end{proof}

Let $\{F_n\}_{n\geq 0}$ be a F\o lner sequence of $P$. Denote by $\lambda_P$ a left Haar measure of $P$. 

\begin{lemma}\label{weak stars of folner}
    Take any $x_0\in M$. Let $\nu$ be any weak-$\ast$ limit point of the sequence of probability measures
    \begin{equation*}
        \left\{\frac{1}{\lambda_P(F_n)}\int_{F_n}\delta_{px_0}d\lambda_P(p)\right\}_{n\geq 0}.
    \end{equation*}
     Then for any $y\in\supp(\nu)$ and any $ \epsilon>0$ there exists $p\in P$ such that $d_M(px_0,y)< \epsilon$.
\end{lemma}

\begin{proof}
    If $U$ is any neighborhood of $y$ in $M$, then $\nu(U)>0$, so that for one of the measures of the sequence appearing in the lemma, $U$ has positive measure for such measure, which yields some $p\in P$ such that $px_0\in U$.
\end{proof}

Let $\nu=\int\mu_yd\nu(y)$ be the $P$-ergodic decomposition of the $P$-invariant probability measure $\nu$, where $\nu$ is as in the previous lemma. 

\begin{corollary}\label{ergodic components of weak star}
    For $\nu$-a.e. $y$ we have the following: for $\mu_y$-a.e. $z$ and any $ \epsilon>0$ there exists $p\in P$ such that $d_M(px_0,z)< \epsilon$
\end{corollary}
\begin{proof}
    We have
    \begin{equation*}1=\nu(\supp(\nu))=\int\mu_y(\supp(\nu))d\nu(y),
    \end{equation*}
    which implies $\mu_y(\supp(\nu))=1$ for $\nu$-a.e. $y$. Thus, our corollary follows from the previous Lemma.
\end{proof}


\begin{lemma}\label{generic transverse Pesin manifold intersects all Gorbits}
     If $c(G)\geq 3$, then $\w^T(x)$ intersects all $G$-orbits for $\mu$-a.e. $x$
\end{lemma}


\begin{proof}
    Let $x\in M$, whose genericity with respect to $\mu$ will be indicated throughout the proof. Suppose Proposition \ref{P-equivariance of transverse Pesin manifolds} and Lemma \ref{transverse Pesin manifold is uniformly locally nice} holds for $x$. Consider the set
    \begin{equation*}
        X=\{y\in M: gy\in \w^T(x), \text{ for some }g\in G\}
    \end{equation*}
    $X$ is open by our choice of $x$ and by Lemma \ref{transverse Pesin manifold is uniformly locally nice}. We want to show $X$ is closed as well, as $M$ is connected, this proves our lemma. 

    Take any $z\in\overline{X}$. Let $\nu$ be a weak-$\star$ limit of the sequence 
    \begin{equation*}
        \left\{\frac{1}{\lambda_P(F_n)}\int_{F_n}\delta_{pz}d\lambda_P(p)\right\}_{n\geq 0},
    \end{equation*}
    as in Lemma \ref{weak stars of folner}, built with respect to $z$. Choose $\mu_1$ a $P$-ergodic component of $\nu$ for which Corollary \ref{ergodic components of weak star} holds with respect to $z$. $\mu_1$ is $a$-ergodic for all $a\in\w^P$, see Remark \ref{mu is A-ergodic}.

    We claim that $\mu_1=\mu$.

    Fix any $a\in\w^P$. Choose $y_1\in M$ such that $\w^{u,1}(y_1;a)=\w^T(y_1)$ exists, and such that Lemma \ref{Pattraction} holds for $y_1$ with respect to $a$ and $\mu_1$. We may assume that the same is true for $x$ with respect to $a$ and $\mu$. Assume Corollary \ref{ergodic components of weak star} holds for $y_1$ with respect to $z$, as well.  

    By our choice of $y_1$ with respect to Corollary \ref{ergodic components of weak star} there exists some $p\in P$ such that $pz$ is close enough to $y_1$ so that $\w^T(y_1)$ intersects the open set $X$, which holds as $X$ is $G$-invariant and by considering Remark \ref{why name transverse}.

    Thus, because of our choices of $x$ and $y_1$ with respect to Lemma \ref{Pattraction}, for $\mu$ and $\mu_1$ respectively, it follows that $\mu_1=\mu$. 

    Finally, let us assume that $x\in \supp(\mu)$ and that Birkhoff ergodic theorem holds for $y_1$ with respect to $a$ and $\mu_1=\mu$. These two assumptions imply that for any $ \epsilon>0$ there exists some $N>1$ such that $d_M(a^Ny_1,x)< \epsilon$. By our choices above, there exists some $p\in P$ such that $pz$ is close enough to $y_1$ so that there exists some $g_1\in G$ such that $g_1z\in \w^T(y_1)$.

    Since $d_M(a^ng_1z,a^ny_1)\rightarrow 0$ as $n\rightarrow\infty$, there exists some large enough $N>1$ such that $a^Ny_1$ is close enough to both $x_1$ and $a^Ng_1z$; so that $G\cdot z$ intersects $\w^T(x)$, by considering Remark \ref{why name transverse}, which means $z\in X$.

    Therefore, $X$ is both open and closed, hence $X=M$, and this proves the lemma. 
\end{proof}
We can now prove that the restriction $P\acts M$ of our action $G\acts M$ is uniquely ergodic. 

\begin{proof unique ergodicity}
    Let $\mu_1$ and $\mu_2$ be $P$-ergodic and invariant Radon probability measures on $M$. Fix any $a\in\w^P$, so that both $\mu_1$ and $\mu_2$ are $a$-ergodic because of Remark \ref{mu is A-ergodic}.

    Pick $x_j$ such that Proposition \ref{Pattraction} holds for $x_j$ with respect to $a$ and $\mu_j$, for $j=1,2$; and such that Lemma \ref{generic transverse Pesin manifold intersects all Gorbits} holds for $\w^T(x_j)$, for $j=1,2$. Then the orbit $G\cdot x_1$ intersects $\w^T(x_2)$, this implies by our choices above and Proposition \ref{Pattraction} that $\mu_1=\mu_2$. This proves the unique ergodicity statement. 
\end{proof unique ergodicity}

\begin{notation}\label{notation muP}
    From this point of the paper on, denote by $\mu^P$ the unique $P$-invariant Radon probability measure on $M$ of the restricted action $P\acts M$ of our action $G\acts M$, for the $c(G)\geq 3$ case.
\end{notation}

\section{Proof of main result when $c(G)>2$}\label{Proof of main result when G is not locally isomorphic to SL(3,R) Section}

In this section we combine Proposition \ref{P-equivariance of transverse Pesin manifolds} and Lemma \ref{generic transverse Pesin manifold intersects all Gorbits} to construct, in the $c(G)\geq 3$ case,  a one dimensional and $G$-invariant foliation on $M$ such that all of its leaves are everywhere transversal to the $G$-orbits, along with a smooth map $M\rightarrow G/\Gamma$
for some $\Gamma\leq G$ discrete subgroup, whose fibers are the leaves of this one-dimensional foliation; however the leaves of this one dimensional foliation are $G$-translates of transverse Pesin manifolds, which cannot be compact, this will give a contradiction to Assumption \ref{assumption no invariant probability}. 

\begin{assumption}\label{assumption strict maximal}
    In this section, Section \ref{Proof of main result when G is not locally isomorphic to SL(3,R) Section}, we work under Assumption \ref{assumption no invariant probability}, and we also assume that the strict inequality of Proposition \ref{Strict Maximal Inequality of transverse Lyapunov exponent} holds for $\mu^P$. 
\end{assumption}

\subsection{For generic $x$, $\{g\w^T(x)\}_{g\in G}$ is a foliation of $M$}

\begin{lemma}\label{equality from intersections}
    For $\mu^P$-a.e. $x$ and any $g\in G$, if $g\cdot\w^T(x)\cap \w^T(x)\neq\emptyset$, then $g\cdot \w^T(x)=\w^T(x)$.
\end{lemma}

\begin{proof}

    Fix any $a\in\w^P$. Take any $x\in M$ satisfying Proposition \ref{essential G equivariance of transverse Pesin} with respect to $a$. Suppose $g\in G$ is  such that $g\cdot\w^T(x)\cap \w^T(x)\neq\emptyset$.

    By Proposition \ref{essential G equivariance of transverse Pesin} we have that there exists $g_1\in G$ such that $g_1g\w^T(x)=\w^u_{\lo(a^{-1})}(g_1gx;a^{-1})$ and $g_1\w^T(x)=\w^u_{\lo(a^{-1})}(g_1x;a^{-1})$. By our assumption on $g$ we get that $g_1g\w^T(x)\cap g_1\w^T(x)\neq\emptyset$, i.e., $\w^u_{\lo(a^{-1})}(g_1gx;a^{-1})\cap \w^u_{\lo(a^{-1})}(g_1x;a^{-1})\neq \emptyset$. 
    
    The last non-empty intersection implies, from the definition of $\w^u_{\lo(a^{-1})}(-,a^{-1})$, see (\ref{definition unstable set}), that  $\w^u_{\lo(a^{-1})}(g_1gx;a^{-1})= \w^u_{\lo(a^{-1})}(g_1x;a^{-1})$. In turn, from our choice of $g_1$, the last equality implies that $g\w^T(x)=\w^T(x)$. 
\end{proof}

\begin{lemma}\label{G-flow from generic point}
    For $\mu^P$-a.e. $x$, the family $\{g\cdot\w^T(x)\}_{g\in G}$ forms a $C^2$ foliation of $M$.
\end{lemma}
\begin{proof}

    Because of Lemmas \ref{equality from intersections} and \ref{generic transverse Pesin manifold intersects all Gorbits}, we obtain that $\{g\cdot\w^T(x)\}_{g\in G}$ forms a partition of $M$.

    Assume $x$ satisfies Lemma \ref{transverse Pesin manifold is uniformly locally nice}, then for any $y\in\w^T(x)$ and any $g\in G$ there exist a neighborhoov $U$ of $y$, with respect to the immersed topology of $\w^T(x)$, along with a neighborhood $V$ of $e_G$ in $G$ such that the map
    \begin{equation*}
    \begin{split}
        U\times V&\rightarrow M\\
        (h,z)&\mapsto hgz
    \end{split}
    \end{equation*}
    is a diffeomorphism onto its image; these maps form a family of charts of $M$ making the partition $\{g\cdot\w^T(x)\}_{g\in G}$ a $C^2$ foliation of $M$ because of Lemma \ref{equality from intersections} and as $G\acts M$ is $C^2$. 
\end{proof}

We need to prove that under Assumption \ref{assumption no invariant probability}, for a $\mu^P$-generic $x$, $\w^T(x)$ is diffeomorphic to $\R$.

\begin{lemma}\label{generic transverse pesin is R}
    For $\mu^P$-a.e. $x$, $\w^T(x)$ is diffeomorphic to $\R$. 
\end{lemma}

\begin{proof}

    For $\mu^P$-a.e. $x$ we have the following two properties:
    \begin{enumerate}[label=(\roman*)]
        \item\label{yurtj1} $\w^T(x)$ exists and it intersects $G$-orbits  transversely, by Lemma \ref{transverse Pesin manifold is uniformly locally nice};
        \item\label{yurtj2} for every $y\in\w^T(x)$ and every $a\in\w^P$ we have $\lim_{n\rightarrow\infty}d(a^nx,a^ny)=0$, by Definition \ref{definition transverse Pesin manifolds} and Lemma \ref{negativity of transversal lyapunov exponent on positive weyl chambers};
        \item\label{yurtj3} $\w^T(x)$ is a connected injectively immersed one-dimensional submanifold of $M$
    \end{enumerate}

    Take any $x_0\in M$ satisfying these three properties. We claim that $\w^T(x_0)$ is diffeomorphic to $\R$.

    Suppose not, then, because of \ref{yurtj3}, $\w^T(x_0)$ is diffeomorphic to a circle, then, by \ref{yurtj2} for some large enough $n_0>0$ we have
    \begin{equation*}
        \operatorname{diam}_M(a^n\w^T(x_0))<\operatorname{injrad}(M),
    \end{equation*}
   i.e., $a^n\w^T(x_0)$ is contained in an open subset of $M$ diffeomorphic to $\R^{\dim(M)}$, this implies that since  $a^n\w^T(x_0)$ is diffeomorphic to a circle, there exists a point of it at which it is tangent to an orbit of $G\acts M$, which contradicts \ref{yurtj1}.
\end{proof}

\subsection{$G$-flows}

\begin{lemma}\label{sliding generic transverse pesin}
    For $\mu^P$-a.e. $x$ 
    \begin{equation}
    \begin{split}
        G\times\w^T(x)&\xrightarrow{\sigma} M\\
        (g,y) &\mapsto gy
    \end{split}
    \end{equation}
    is a covering space map.
\end{lemma}

\begin{proof}

    Choose $\delta>0$ such that for all $g\in G^{<\delta}$ and all $x\in M$ we have 
    \begin{equation*}
        d(gx,x)<\operatorname{injrad}(M)/2.
    \end{equation*}
    Pick any $x_0\in M$ for which Lemma \ref{G-flow from generic point} holds. Consider on $\w^T(x_0)$ the Riemannian metric induced from the one on $M$. From our choice of $\delta$ it follows that for any $y\in M$, any $g\in G$ and any $z\in \w^T(x)$ such that $gz=y$, 
    \begin{equation*}
    \begin{split}
        G^{<\delta}\cdot g\times B^{\w^T(x_0)}(z,\operatorname{injrad}(M)/2)&\rightarrow M\\
        (hg,w)\mapsto hgw
    \end{split}
    \end{equation*}
    is a local diffeomorphism onto its image. This proves that $\sigma$ as in our lemma is a covering space map. 
\end{proof}

To finish the proof of Theorem \ref{main} when $G$ is not locally isomorphic to $\SL(3,\R)$, a key tool is the notion of $G$-flows, here we use the following framing of it appearing in the work of Carri\'ere, see \cite{carriere1}.

\begin{definition}[$\mathfrak g$-flow]

    Given a simply connected Lie group $G$ with Lie algebra $\mathfrak g$, a $\mathfrak g$-flow on a manifold $N$ is a one-dimensional foliation $\mathcal H$ of $N$ satisfying the following properties:
    \begin{enumerate}[label=(\roman*)]
        \item there exists a fibration $D:\tilde N\rightarrow G$, where $\tilde N$ is the universal cover of $N$, such that the leaves of the lifted foliation of $\mathcal H$ on $\tilde N$ are the fibers of $D$.
        \item There exists a representation $h:\pi_1(N)\rightarrow G$ for which $D$ is equivariant, i.e., $D(\gamma. x)=h(\gamma)\cdot D(x)$, for all $x\in \tilde N$ and all $\gamma\in\pi_1(N)$.
    \end{enumerate}
\end{definition}

Denote by $\tilde G$ the universal cover of $G$, and denote by $\tilde p:\tilde G\rightarrow G$ the corresponding universal cover, which is a Lie groups morphism. Denote by $(g,x)\mapsto gx$ the lifted action $\tilde G\acts M$ of our original action $G\acts M$.

Denote by $\tilde\sigma$ the composition of the two maps 
\begin{equation*}
    \tilde G\times\w^T(x)\xrightarrow{\tilde p\times\operatorname{Id}} G\times\w^T(x)\xrightarrow{\sigma} M.
\end{equation*}

\begin{remark}\label{universal cover of M}
    Lemmas \ref{generic transverse pesin is R} and \ref{sliding generic transverse pesin} imply that for $\mu^P$-a.e. $x$,  $\tilde G\times\w^T(x)$ is the universal cover of $M$, and that $\tilde\sigma$ is the corresponding covering map, so in particular this holds for $\mu^P$-a.e. $x$. We consider $\pi_1(M)$ as the corresponding deck transformations group, and denote its action on $\tilde G\times\w^T(x)$ by $\gamma.(g,y)$.
\end{remark}

\begin{lemma}\label{g-flow from generic pesin}
    For $\mu^P$-a.e. $x$, the one-dimensional foliation $\{g\cdot\w^T(x)\}_{g\in G}$ is a $\lie(G)$-flow on $M$.
\end{lemma}

\begin{proof}
      Fix $x_0\in M$ for which both Lemma \ref{equality from intersections} and Remark \ref{universal cover of M} hold, the set of all such points has full $\mu^P$ measure. 

     For any $\gamma\in\pi_1(M)$, we have that $\gamma.(e_{\tilde G},x_0)=(g_\gamma,x_\gamma)$ for some $x_\gamma\in\w^T(x_0)$.  
     
     Because of the construction of $\tilde\sigma$ and of the action of $\pi_1(M)$ above, we get that $g_\gamma x_\gamma=x_0$; in particular $\w^T(x_0)\cap g_\gamma\w^T(x_0)\neq\emptyset$, and thus $g_\gamma\w^T(x_0)=\w^T(x_0)$ by Lemma \ref{equality from intersections}.

     Let us see that for any $\gamma\in\pi_1(M)$ and any $(g,y)\in\w^T(x)$ we have
     \begin{equation}\label{deck transforms}
         \gamma.(g,y)=(gg_\gamma,g_\gamma^{-1}y).
     \end{equation}

     Let $\alpha$ be any continuous path in $\tilde G\times\w^T(x_0)$ going from $(g,y)$ to $(e_{\tilde G},x_0)$; let $c_\gamma$ be any continuous path in $\tilde G\times\w^T(x_0)$ going from $(e_{\tilde G},x_0)$ to $(g_\gamma,g_\gamma^{-1}x_0)$, note that $c_\gamma$ projects to $\gamma\in\pi_1(M,x_0)$ under $\tilde\sigma$. Let $\beta$ be the path in $\tilde G\times\w^T(x_0)$ opposite to $\alpha$.
     
     Let $(\beta_1,\beta_2)$ be the coordinates of $\beta$ in $\tilde G\times\w^T(x_0)$. Consider the map $\tilde \beta$ whose coordinates in $\tilde G\times\w^T(x_0)$ are $(\beta_1\cdot g_\gamma,g_\gamma^{-1}\cdot\beta_2)$, so that $\beta$ goes from $(g_\gamma,g_\gamma^{-1}x_0)$ to $(gg_\gamma,g_\gamma^{-1}y)$. Then $\tilde \sigma_\ast(\tilde \beta)=\tilde \sigma_\ast(\beta)$, and in particular, $\tilde \sigma_\ast(\tilde \beta)$ is the path opposite to $\tilde \sigma_\ast(\alpha)$ in $M$. Thus, the concatenated path $\tilde \sigma_\ast(\alpha)\cdot[\tilde\sigma_\ast(c_\gamma)\cdot \tilde \sigma_\ast(\tilde \beta)]$ is a loop at $gy$ in $M$ representing $\gamma$ in $\pi_1(M)$. This implies (\ref{deck transforms}).

     Note that (\ref{deck transforms}) implies that $h:\pi_1(M)\rightarrow \tilde G$ given by $h(\gamma)=g_\gamma^{-1}$, is a group homomorphism; it further implies that if we set $D:\tilde G\times\w^T(x)\rightarrow \tilde G$ to be $(g,y)\mapsto g^{-1}$, then
     \begin{equation}\label{equivariance of D}
         D(\gamma.(g,y))=g_\gamma^{-1}g^{-1}=h(\gamma)g^{-1}=h(\gamma)D(g,y).
     \end{equation}
     Since $D$ is a fiber bundle, and its fibers are the lift under $\tilde\sigma$ of the foliation $\{g\w^T(x_0)\}_{g\in M}$ of $M$, we obtain our lemma. 
\end{proof}

Set $\Gamma:=h(\pi_1(M))$, which is a subgroup of $\tilde G$. 

\begin{corollary}\label{map psi_x}
    For $\mu^P$-a.e. $x$, there exists a unique smooth map $\psi_x:M\rightarrow\tilde G/\Gamma$ making the diagram
    \[
\begin{tikzcd}
\tilde G\times\w^T(x) \arrow[r, "D"] \arrow[d, "\tilde \sigma"]
  & \tilde G \arrow[d, "\pi"] \\
M \arrow[r, dotted, "\psi_x"']
  & \tilde G/\Gamma
\end{tikzcd}
\]
commute, where $\pi:\tilde G\rightarrow \tilde G/\Gamma$ is the canonical quotient map. 
\end{corollary}

\begin{proof}
    This follows immediately from the proof of Lemma \ref{g-flow from generic pesin} since we consider $\pi_1(M)$ as the group of deck transformations of $\tilde\sigma$, and because of (\ref{equivariance of D}).
\end{proof}

\subsection{Proof of Theorem \ref{main} when $c(G)>2$}\label{proof of main when not SL3R Section}

To prove Theorem \ref{main} in the $c(G)>2$ case, we need the following observation.

\begin{lemma}\label{nice preimage}
    For $\mu^P$-a.e. $x$, we have that for the map $\psi_x$, appearing in Corollary \ref{map psi_x}, $\psi_x^{-1}(\{e_{\tilde G}\Gamma \})=\w^T(x)$
\end{lemma}

\begin{proof}

    Pick any $x_0\in M$ for which Lemma \ref{g-flow from generic pesin} holds. To ease notation, denote $\psi_{x_0}$ by $\psi$ in this proof.

    For any $y\in\w^T(x_0)$ we obtain from Corollary \ref{map psi_x}
    \begin{equation*}
        \psi(y)=\psi(\tilde\sigma(e_{\tilde G},y))=\pi\circ D(e_{\tilde G},y)=\pi(e_{\tilde G})=e_{\tilde G}\Gamma,
    \end{equation*}
    and hence, $\w^T(x_0)\subseteq\psi^{-1}(\{e\Gamma\})$.

    Let $y\in\psi^{-1}(\{e_{\tilde G}\Gamma\})$. Since $\tilde\sigma$ is a covering map, there exists $(g,z)\in \tilde G\times\w^T(x)$ such that $\tilde\sigma(g,z)=y$, which implies by Lemma \ref{g-flow from generic pesin} that 
    \begin{equation*}
        \pi\circ D(g,z)=\psi\circ\tilde\sigma(g,z)=g\Gamma
    \end{equation*}
    so as $\psi(y)=e_{\tilde G}\Gamma$, we get that $g\Gamma=e_{\tilde G}\Gamma$, i.e., $g\in\Gamma$. Hence, $g=h(\gamma)$ for some $\gamma\in\pi_1(M)$, i.e., $g=g_\gamma^{-1}$ as in the proof of Lemma \ref{g-flow from generic pesin}.

    However, because of (\ref{deck transforms}), we have
    \begin{equation*}
        \gamma.(e,g_\gamma^{-1}z)=(g_\gamma^{-1},z)=(g,z).
    \end{equation*}

    Hence, as we are using the deck transformation action with respect to the covering map $\tilde\sigma$, this implies that $(e,g_\gamma^{-1}z)$ and $(g,z)$ have the same image under $\tilde\sigma$, thus $y=\tilde\sigma(g,z)=\tilde\sigma(e,g_\gamma^{-1}z)=g_\gamma^{-1}z\in\w^T(x_0).$ 

    This proves that $\psi^{-1}(\{e_{\tilde G}\Gamma\})\subseteq\w^T(x_0)$. Combining this with what we did at the beginning of this proof yields the sought after equality.
\end{proof}

\begin{proof of main when c(G)>2}
    Choose $x_0$ for which Lemmas \ref{nice preimage} and \ref{generic transverse pesin is R} both hold. Hence $\phi^{-1}_{x_0}(\{e_{\tilde G}\Gamma\})$ is closed in $M$, and thus, it is a compact subset of it, which contradicts \ref{generic transverse pesin is R}.
\end{proof of main when c(G)>2}

\section{Unique $P$-ergodicity when $c(G)=2$}\label{unique ergodicity sl3 section}

For the proof of Theorem \ref{main} when $G$ is locally isomorphic to $\SL(3,\R)$, we need to establish unique ergodicity of the restricted action $P\acts M$ with respect to the minimal parabolic $P\leq G$, under Assumption \ref{assumption no invariant probability}. This was done in Section \ref{Unique ergodicity of P section} when $G$ is not locally isomorphic to $\SL(3,\R)$, using Proposition \ref{Strict Maximal Inequality of transverse Lyapunov exponent}. However, the inequalities appearing in Proposition \ref{Strict Maximal Inequality of transverse Lyapunov exponent}, a priori, do not necessarily hold in the $c(G)=2$, hence, we need to study what is the behavior of $P$-invariant and $P$-ergodic probabilities measures in this case.

Because of Remark \ref{special case of cG=2 case}, to prove unique ergodicity of $P\acts M$, it is enough to work under the following assumption:

\begin{assumption}\label{combinatorial assumption on G rG=2 case}
    Under Assumption \ref{assumption no invariant probability},  assume further that the following properties hold for $G$ and some fixed $P$-invariant and $P$-ergodic $\mu$: 
    \begin{enumerate}[label=(\roman*)]
    \item\label{yacasi1} $\operatorname{Lie}(G)\simeq\mathfrak{sl}(3,\R)$;
    \item\label{yacasi2} under the isomorphism above, $\operatorname{Lie}(P)$ corresponds to the Lie subalgebra of $\mathfrak{sl}(3,\R)$ of upper triangular matrices of $\SL(3,\R)$;
    \item\label{yacasi3} the inequality of Proposition \ref{Strict Maximal Inequality of transverse Lyapunov exponent} fails on $\ker(\beta_{1,2})$, i.e., 
    \begin{equation*}
        \lo(s^{-1})\leq\beta_{2,3}(s)
    \end{equation*}
    for all $s\in\ker(\beta_{1,2})\cap\overline{\w^P}$. Fix such an $s_0$.
\end{enumerate}
\end{assumption}

The other possibility, where the strict inequality in Proposition \ref{Strict Maximal Inequality of transverse Lyapunov exponent} fails on $\ker(\beta_{2,3})$, is analogous. 

\subsection{Failure of Strict Maximal Inequality}

Fix $\mu$ any $P$-invariant and $P$-ergodic probability measure on $M$. 

\begin{lemma}\label{non resonance when failure of strict}
    If the conclusion of Proposition \ref{Strict Maximal Inequality of transverse Lyapunov exponent} does not hold for $\mu$, then $c(G)=2$, and $\lo$ is not resonant to a root in $\Sigma_{(G,A)}$.
\end{lemma}

\begin{proof}

    Because of Proposition \ref{Strict Maximal Inequality of transverse Lyapunov exponent} we have that $c(G)=2$.

    We have $\operatorname{Lie}(G)\cong\mathfrak{sl}(3,\R)$, because of Remark \ref{locally isometric to SL(3,R) when Strict Inequality fails}; so we may assume, without loss of generality, that $P$ is the minimal parabolic of $G$ whose Lie algebra in $\mathfrak{sl}(3,\R)$ corresponds to the Lie subalgebra of the  subgroup of $\SL(3,\R)$ of upper triangular matrices. We may also assume, without loss of generality, that the inequality of Proposition \ref{Strict Maximal Inequality of transverse Lyapunov exponent} fails on $\ker(\beta_{1,2})$, i.e., that
    \begin{equation}\label{qoilki}
        \lo(s^{-1})\leq\beta_{2,3}(s)=\beta_{1,3}(s)
    \end{equation}
    for some (and hence for all) $s\in\ker(\beta_{1,2})\cap\overline{\w^P}$. Fix $s_0\in\ker(\beta_{1,2})\cap\overline{\w^P}\setminus\{e\}$.

    Suppose, to get a contradiction, that $\lo$ is resonant to a root. By Lemma \ref{kernel of transverse exponent touches not weyl chamber} we have that $\lo$ is resonant to $\beta_{3,1}$, and denote the corresponding coarse Lyapunov exponent by $\chi:=\{\lo,\beta_{3,1}\}$. However, because of Proposition \ref{Maximal negativity of transverse Lyapunov exponent}, we have that since $s_0\in\ker(\beta_{1,2})\cap\overline{\w^P}$, $\beta_{2,3}(s_0)\leq\lo(s_0^{-1})$, so that from (\ref{qoilki}), we get that $\beta_{2,3}(s_0)=\lo(s_0^{-1})$. However, $\beta_{1,3}(s_0)=\beta_{2,3}(s_0)$ as $s_0\in\ker(\beta_{1,2})$, and thus, since $\lo$ and $\beta_{3,1}$ are resonant, we conclude that $\lo=\beta_{3,1}$.

    We now finish the proof of our lemma in two steps. 

\begin{step no resonance of lo when strict fails}\label{step no resonance of lo when strict fails 1}
$\mu$ is $G^{\beta_{3,2}}$-invariant
\end{step no resonance of lo when strict fails}

Since $\mu$ is $P$-invariant, from the codimension one assumption on $G\acts M$ and from Lemmas \ref{kernel of transverse exponent touches not weyl chamber} and \ref{negativity of transversal lyapunov exponent on positive weyl chambers}, we obtain that $h_\mu(s_0)=\beta_{2,3}(s_0)+\beta_{1,3}(s_0)=2\beta_{1,3}(s_0)$, since $s_0\in\ker(\beta_{1,2})\cap\overline{\w^P}$. However, since $\lo$ is resonant to $\beta_{3,1}$, it follows from Theorem \ref{entropy product structure} that 
\begin{equation*}
    h_\mu(s_0^{-1})=h_\mu(s_0^{-1}|\w^{\chi})+h_\mu(s_0^{-1}|\w^{\beta_{3,2}}),
\end{equation*}
as $\beta_{1,2}(s_0)=0$. Thus, by using that $h_\mu(s_0)=h_\mu(s_0^{-1})$, we obtain
\begin{equation*}
    2\beta_{1,3}(s)=h_\mu(s^{-1}|\w^\chi)+h_\mu(s^{-1}|\w^{\beta_{3,2}})\leq \lo(s^{-1})+\beta_{2,3}(s)=2\beta_{1,3}(s);
\end{equation*}
where the last inequality follows from Theorem \ref{Ruelle-Margulis inequality} and Lemma \ref{no full entropy along coarse Pesin manifold}. Then
\begin{equation*}
    h_\mu(s^{-1}|\w^{\beta_{3,2}})=\beta_{1,3}(s)=\beta_{2,3}(s), \text{ } h_\mu(s^{-1}|\w^\chi)=\beta_{1,3}(s);
\end{equation*}
in particular, $\mu$ is $G^{\beta_{3,2}}$-invariant by applying Theorem \ref{ledrappier young conditional entropy}.
\begin{step no resonance of lo when strict fails}

Contradiction to $\lo$ being resonant to a root.

\end{step no resonance of lo when strict fails}

For any $a\in\w^P$ we have
\begin{equation*}
    h_{\mu}(a^{-1})=h_\mu(a^{-1}|\w^\chi)+h_\mu(a^{-1}|\w^{\beta_{3,2}})=h_\mu(a^{-1}|\w^\chi)+\beta_{2,3}(a)
\end{equation*}
by applying Theorem \ref{entropy product structure} with Step 1 above; since $h_\mu(a)=\sum_{\beta\in\Sigma_P}\beta(a)=2\beta_{1,3}(a)$, we get that 
\begin{equation*}
    h_\mu(a^{-1}|\w^\chi)=2\beta_{1,3}(a)-\beta_{2,3}(a),
\end{equation*}
for all $a\in\w^P$. However, since $\chi$ does not change sign on $\w^P$, we get from Theorem \ref{ledrappier young conditional entropy} that $a\mapsto h_\mu(a^{-1}|\w^\chi)$ is linear on $\w^P$ and resonant to $\beta_{1,3}$, so we can extend this map to a linear function $\rho$ on all of $A$ resonant to $\beta_{1,3}$, and such that $\rho=2\beta_{1,3}-\beta_{2,3}$. However, this contradicts the fact that the set $\{\beta_{1,3},\beta_{2,3}\}$ is linearly independent.
\end{proof}

\begin{lemma}\label{maximal parabolic stabilizer when SMI fails}
    The following properties hold:
    \begin{itemize}
        \item $Q_0=\operatorname{Stab}(\mu)$ equals the maximal parabolic subgroup $\langle P,G^{\beta_{3,2}}\rangle$ of $G$;
        \item $\mu$ disintegrates along $\w^T_\mu(x)$ equivalently to Riemannian volume for $\mu$-a.e. $x$;
        \item $\lo=\beta_{2,1}+\beta_{3,1}$.
    \end{itemize}
\end{lemma}
\begin{proof}
From Proposition \ref{Maximal negativity of transverse Lyapunov exponent} and \ref{yacasi3} we obtain that $\lo(s_0^{-1})=\beta_{2,3}(s_0)$. 

Since $\lo(s_0)<0$, because of Corollary \ref{negativity on faces of weyl chamber}, we get that $h_\mu(s_0)=\beta_{2,3}(s_0)+\beta_{1,3}(s_0)=2\beta_{2,3}(s_0)$; from Theorem \ref{entropy product structure} and since $\lo$ is not resonant to a root, because of Lemma \ref{non resonance when failure of strict}, we get that
\begin{equation*}
    h_\mu(s_0^{-1})=h_\mu(s_0^{-1}|\w^{\beta_{3,2}})+h_\mu(s_0^{-1}|\w^{\lo})+h_\mu(s_0^{-1}|\w^{\beta_{3,1}}).
\end{equation*}

However, $h_\mu(s_0^{-1}|\w^{\beta_{3,1}})=0$, for otherwise, by the High Entropy Method, Theorem \ref{high entropy method}, using the $P$-invariance of $\mu$, we would obtain that $\mu$ would be invariant by $G^{\beta_{2,1}}$ and $G^{\beta_{3,2}}$, and then $\mu$ would be invariant under all of $G$, contradicting Assumption \ref{assumption no invariant probability}. Thus, using that $h_\mu(s_0)=h_\mu(s_0^{-1})$, we obtain 
\begin{equation*}
    h_\mu(s_0^{-1}|\w^{\beta_{3,2}})+h_\mu(s_0^{-1}|\w^{\lo})=2\beta_{2,3}(s_0).
\end{equation*}
However, $h_\mu(s_0^{-1}|\w^{\beta_{3,2}})\leq\beta_{3,2}(s_0^{-1})=\beta_{2,3}(s_0)$, and $h_\mu(s_0^{-1}|\w^{\lo})\leq \lo(s_0^{-1})=\beta_{2,3}(s_0)$ because of the Ruelle-Margulis inequality, Theorem \ref{Ruelle-Margulis inequality}, we obtain that $h_\mu(s_0^{-1}\mid\w^{\beta_{3,2}})=\beta_{2,3}(s_0^{-1})$ and $h_\mu(s_0^{-1}\mid\w^{\lo})=\beta_{2,3}(s_0)=\lo(s_0^{-1})$, and therefore, by Theorem \ref{ledrappier young conditional entropy}, we obtain that $\mu$ is $G^{\beta_{3,2}}$-invariant, and in particular $\mu$ is invariant under a maximal parabolic of $G$; and that $\mu^P$ disintegrates along $\w^T(x)$ for $\mu^P$-a.e. $x$. 

Finally, for $a\in\w^P$ we have 
\begin{equation*}
    h_{\mu}(a)=(\beta_{1,2}+\beta_{2,3}+\beta_{1,3})(a); \text{ }h_{\mu}(a^{-1})=\beta_{3,2}(a^{-1})+\lo(a^{-1}),
\end{equation*}
where the last inequality follows from what we did above and Theorem \ref{ledrappier young conditional entropy}. Hence, since $h_\mu(a)=h_\mu(a^{-1})$, we conclude that 
$\lo(a^{-1})=(\beta_{1,2}+\beta_{1,3})(a)$ for all $a\in\w^P$, which implies that $\lo=\beta_{2,1}+\beta_{3,1}$.
\end{proof}

\begin{corollary}\label{maximal negativity case sl3r}
     Proposition \ref{Maximal negativity of transverse Lyapunov exponent} holds for $\mu$.
\end{corollary}
\begin{proof}

    By Lemma \ref{maximal parabolic stabilizer when SMI fails} it follows that for $a\in\w^P$, 
    \begin{equation*}
        \lo(a^{-1})=(\beta_{2,1}+\beta_{3,1})(a^{-1})=\beta_{1,2}(a)+\beta_{1,3}(a)>\beta_{1,3}(a).
    \end{equation*}
    Hence, since $\beta_{1,3}$ is the highest root of $\Sigma_P$, we obtain the inequality in Proposition \ref{Maximal negativity of transverse Lyapunov exponent}.
\end{proof}

\subsection{Some Observations about semisimple Lie groups}

\begin{lemma}\label{Pertubations by generic elements of $P$ have $LU$ factorization}
Let $G$ be any connected semisimple Lie group with finite center, and $P\leq G$ a minimal parabolic. For any $g\in G$ and for $\haar_P$-a.e. $p$, $pg=up'$ for unique $u\in U(P^{\operatorname{op}})$ and $p'\in P$. Furthermore, the $\haar_P$-a.e. defined function $p\mapsto p'$ is absolutely continuous.    
\end{lemma}

\begin{proof}
    Set $U'$ to be $U(\Pop)$. 

    From the ``big cell" lemma, see \cite[Theorem~21.84]{milne} and \cite[Lemma~5.1.4]{zimmerbook}, we have that the multiplication map $U'\times P\rightarrow G$, $(u,p)\mapsto up$ is injective, and has open image in $G$ in both the Hausdorff and Zariski topologies. Since multiplication by $g$ is birational(with coefficients in $\R$), we get that $U'Pg^{-1}$ is also open in $G$ in the Zariski topology.

    $P$ is a closed subset of $G$ in the Zariski topology, hence $U'Pg^{-1}\cap P$ is open in $P$ with its Zariski topology. Since $\haar_P$ is equivalent to a volume form on $P$, we get that any proper closed Zariski subset of $P$ has zero $\haar_P$ measure. Hence, our claim will follow if we show $U'Pg^{-1}\cap P\neq\emptyset$, for all $g\in G$, which is equivalent to showing that $G=PU'P$

     Applying again the ``big cell" lemma, we get that both $PU'$ and $U'P$ are open dense in $G$; however the product of a dense and an open subset of $G$ equals this whole group, so $G=PU'U'P=PU'P$, as $U'U'=U'$, as $U'$ is a subgroup of $G$.

     The map $U'\times P\rightarrow G$ given by $(u,p)\mapsto up$ is a local diffemorphism, and since it is injective, by what we did above, it follows that it is a diffeomorphism onto its image, so it has a smooth inverse defined on its image. Since the map $p\mapsto pg$ is smooth, we obtain that the map $p\mapsto p'$, defined in the statement of our claim, is smooth, and therefore absolutely continuous.
\end{proof}

\begin{definition}\label{definition non-trivial component}
    Let $\alpha\in\Sigma_{P}$. We say $g\in G$ \emph{has \emph{non-trivial} $G^{-\alpha}$ component with respect to $P$} if $g$ can be factored as $up$, for some $p\in P$ and $u\in U(P^{\operatorname{op}})$, and $u$ is such that if $v$ is the element of $\operatorname{Lie}(U(P^{\operatorname{op}}))$ such that $u=\exp(v)$, then $v$ has non trivial $\mathfrak g_{-\alpha}$ component in the decomposition $\operatorname{Lie}(U(P^{\operatorname{op}}))=\bigoplus_{\beta\in\Sigma_{P}}\mathfrak g_{-\beta}$.
\end{definition}

\begin{lemma}\label{almost surely nice perturbations}
Let $P\leq Q\leq G$ be a proper parabolic. Suppose $g\in G\setminus Q$. Then there exists a simple positive root $\alpha\in\Delta_P$ with $-\alpha\notin\Sigma_Q$, such that if $pg=up'$, as in Lemma \ref{Pertubations by generic elements of $P$ have $LU$ factorization}, then $u$ has non-trivial $G^{-\alpha}$ component for $\haar_P$-a.e. $p$.
\end{lemma}

\begin{proof}
    By Lemma \ref{Pertubations by generic elements of $P$ have $LU$ factorization}, choose $p_0\in P$ such that $p_0g=u_1p_1$ for some $u_1\in U(P^{\operatorname{op}})$ and $p_1\in P$. Since $g\notin Q$, it follows that $u_1\in U(P^{\operatorname{op}})\setminus\{e\}$. Then $u_1$ has non-trivial $G^{-\alpha}$ component for some $\alpha\in\Sigma_P$, and let us assume $\alpha$ is minimal with this property with respect to the ordering on $\Sigma_P$.
    
    If $\alpha$ is a simple root, then for $\haar_P$-a.e. $p$, $pu$ has non-trivial $G^{-\alpha}$ component by using the Baker-Campbell-Hausdorff formula, and using the minimality of $\alpha$. 

    If $\alpha=\gamma_1+\gamma_2$ for some $\gamma_1,\gamma_2\in\Sigma_P$, then using the minimality of $\alpha$, it follows that for $\haar_P$-a.e. $p$, $pu$ has non-trivial $G^{-\gamma_1}$ component, applying again the Baker-Campbell-Hausdorff formula. 

    Continuing the last argument by induction, and the argument for when $\alpha\in\Delta_P$, we obtain a proof of our Lemma.  
\end{proof}


\subsection{Unique ergodicity of $P\acts M$ in the $c(G)=2$ case}

\begin{lemma}\label{transverse pesin intersects all G-orbits when rG=2}
     For $\mu$-a.e. $x$, $\w^T(x)$ intersects all $G$-orbits.   
\end{lemma}

\begin{proof}

    The statement of our lemma is equivalent to showing that for $\mu$-a.e. $x$, the set $Y_x:=G\cdot\w^T(x)$, equals $M$.

    Fix $x_0\in M$, whose genericity with respect to $\mu$ is indicated during the proof.

    Suppose $x_0$ is  such that Lemma \ref{transverse Pesin manifold is uniformly locally nice} holds for it. This implies that $Y_{x_0}$ is open. Let us prove that $Y_{x_0}$ is a closed subset of $M$, by taking $x_0$ to be generic enough with respect to $\mu$.

    Take any $y_0\in\overline{Y_{x_0}}$. Since $Y_{x_0}$ is $G$-invariant, then so is $\overline{Y_{x_0}}$. Let $\nu$ be a weak-$\star$ limit of the sequence 
    \begin{equation*}
        \left\{\frac{1}{\lambda_P(F_n)}\int_{F_n}\delta_{py_0}d\lambda_P(p)\right\}_{n\geq 0},
    \end{equation*}
    as in Lemma \ref{weak stars of folner}, built with respect to $y_0$. Choose $\mu_1$ a $P$-ergodic component of $\nu$ for which Corollary \ref{ergodic components of weak star} holds with respect to $y_0$. $\mu_1$ is $a$-ergodic for all $a\in\w^P$, see Remark \ref{mu is A-ergodic}. Since $\overline{Y_{x_0}}$ is $G$-invariant, it follows that $\nu$ is supported on $\overline{Y_{x_0}}$, and we may assume $\mu_1$ is such that $\supp(\mu_1)\subseteq\overline{Y_{x_0}}$. 

    We claim that $\mu_1=\mu$, once we prove this, we can prove that $y_0\in Y_{x_0}$ using the argument given in the proof of Lemma \ref{generic transverse Pesin manifold intersects all Gorbits}, which proves $Y_{x_0}$ is a closed subset of $M$ and finishes the proof of our lemma.
    
    Fix $w_0\in \overline{Y_{x_0}}$ satisfying Birkhoff's ergodic theorem with respect to some fixed $a_0\in\w^P$ and $\mu_1$, which can be done as $\supp(\mu_1)\subseteq\overline{Y_{x_0}}$. Assume also that Lemma \ref{transverse Pesin manifold is uniformly locally nice} holds for $w_0$ with respect to $\mu_1$. Since $Y_{x_0}$ is open and $w_0\in \overline{Y_{x_0}}$, it follows that $\w^T(w_0)$ intersects $Y_{x_0}$. 

    Therefore, there exists a $G$-orbit of $G\acts M$ such that the Pesin manifolds $\w^T(x_0)$ and $\w^T(w_0)$ both intersect it. Fix $z_0\in\w^T(x_0)$ such that $g_0z_0\in\w^T(w_0)$ for some $g_0\in G$.

    We have two cases.

    \begin{transverse pesin intersects all G-orbits when rG=2 case}
        Proposition \ref{Strict Maximal Inequality of transverse Lyapunov exponent} holds for $\mu_1$.
    \end{transverse pesin intersects all G-orbits when rG=2 case}
    
    Suppose $x_0$ is such that $\mu$ disintegrates equivalently to the Riemannian volume on $\w^T(x_0)$. Since holonomy maps along the partition by $G$-orbits is $C^2$, and in particular absolutely continuous, by perturbing as little as necessary $z_0$ inside of $\w^T(x_0)$ if necessary, we may assume that $z_0$ is such that $pz_0$ satisfies Bikhoff's ergodic theorem for $\mu$ with respect to $a$.

    Assume $w_0$ is  such that $\w^T(w_0)$ and $\w^T(px_0)$ exist for $\haar_P$-a.e. $p$, with $p\cdot\w^T(w_0)=\w^T(pw_0)$.

    Because of Lemma \ref{Pertubations by generic elements of $P$ have $LU$ factorization} and our choice above of $w_0$, we can choose $p_0\in P$ such that the following properties hold:
    \begin{enumerate}[label=(\roman*)]
        \item\label{thisfar1} $p_0g_0=u_1p_1$ for some $u_1\in U(P^{\operatorname{op}})
        $ and $p_1\in P$ such that $p_1z_0$ satisfies Birkhoff's ergodic theorem for $\mu$ with respect to $a$;
        \item\label{thisfar3}$p_0\cdot\w^T(w_0)=\w^T(p_0w_0)$;
        \item $p_0w_0$ satisfies Birkhoff ergodic theorem for $\mu_1$ with respect to $a$.
    \end{enumerate}

    Then $p_0g_0z_0\in p_0\cdot\w^T(x_0)=\w^T(p_0x_0)$. However, $p_0g_0z_0=u_1p_1z_0$, and $u_1p_1z_0$ satisfies Birkhoff's ergodic theorem for $\mu$ with respect to $a$ because $u_1\in U(P^{\operatorname{op}})$ and by item \ref{thisfar1}. But $u_1p_1z_0=p_0g_0z_0\w^T(p_0x_0)$; so by \ref{thisfar3}, it follows that also satisfies Birkhoff's ergodic theorem for $\mu_1$ for $\mu_1$ with respect to $a$.

    Therefore $\mu_1=\mu$.

   \begin{transverse pesin intersects all G-orbits when rG=2 case}
        Proposition \ref{Strict Maximal Inequality of transverse Lyapunov exponent} does not hold for $\mu$ nor $\mu_1$. 
   \end{transverse pesin intersects all G-orbits when rG=2 case}

Fix any $a_0\in\w^P$. We may assume we chose $x_0$ and $w_0$ so that they satisfy the following properties:
\begin{enumerate}[label=(\roman*)]
    \item the disintegrations $\mu^{\w^T_\mu}_{x_0}$ and $[\mu_1]_{w_0}^{\w^T_{\mu_1}}$ are equivalent to the Riemannian volume on $\w^T_\mu(x_0)$ and $\w^T_{\mu_1}(w_0)$, respectively, which holds because of Lemma \ref{maximal parabolic stabilizer when SMI fails}, 
    \item\label{unique erogdicity case ii} for $\mu^{\w^T_\mu}_{x_0}$-a.e. $y_0$, $[\mu_1]_{w_0}^{\w^T_{\mu_1}}$-a.e. $y_1$ and $\haar_P$-a.e. $p$ Birkhoff's ergodic theorem holds for $\mu$ and $\mu_1$, respectively, with respect to $a_0$, and with respect to $py_0$ and $py_1$, respectively. 
\end{enumerate}

Recall that $z_0\in\w^T_\mu(x_0)$ and $g_0z_0\in\w_{\mu_1}^T(w_0)$.

On a small open neighborhood $V$ of $z_0$ in $\w^T_\mu(x_0)$(immersed topology) the holonomy map along the foliation by $G$-orbits from $V$ to an open neighborhood of $g_0z_0$ in $\w^T_{\mu_1}(w_0)$ is defined and $C^2$, as $G\acts M$ is $C^2$, so this holonomy map is in particular absolutely continuous; thus, there exists $v_0\in V$ satisfying Property \ref{unique erogdicity case ii} with respect to $\mu$, which is sent via this holonomy map to a point $v'$ in $\w_{\mu_1}(w_0)$ satisfying Property \ref{unique erogdicity case ii} with respect to $\mu_1$. 

However, $v_1=g_1v_0$ for some $g_1\in G$; by Lemma \ref{Pertubations by generic elements of $P$ have $LU$ factorization} we can choose $p_0\in P$ such that $p_0g_1=u_1p_1$ for some $u_1\in U(P^{\operatorname{op}})$ and some $p_1\in P$ such that $p_1v_0$ and $p_0v_1$ satisfy Property \ref{unique erogdicity case ii} with respect to $\mu$ and $\mu_1$, respectively, i.e., 
\begin{equation}\label{birkhoff for mu and mu1}
    \lim_{n\rightarrow\infty}\frac{1}{n}\sum_{j=0}^{n-1}f(a_0^jp_1v_0)=\int fd\mu\text{ and }\lim_{n\rightarrow\infty}\frac{1}{n}\sum_{j=0}^{n-1}f(a_0^jp_0v_1)=\int fd\mu_1
\end{equation}
for all $f\in C^0(M)$. 

However, since $u_1\in U$, we have that $\lim_{n\rightarrow\infty}d(a_0^nu_1x,a_0^nx)=0$ uniformly for $x\in M$, as $M$ is closed, and thus, $\lim_{n\rightarrow\infty}|f(a_0^nu_1x)-f(a_0^nx)|=0$ uniformly for $x\in M$ for any $f\in C^0(M)$, using again that $M$ is compact. But, $p_0v_1=p_0g_1v_0=u_1p_1v_0$, hence
\begin{equation*}
    \lim_{n\rightarrow\infty}\frac{1}{n}\sum_{j=0}^{n-1}f(a_0^jp_0v_1)=\lim_{n\rightarrow\infty}\frac{1}{n}\sum_{j=0}^{n-1}f(a_0^ju_1p_1v_0)=\lim_{n\rightarrow\infty}\frac{1}{n}\sum_{j=0}^{n-1}f(a_0^jp_1v_0),
\end{equation*}
for any $f\in C^0(M)$,
therefore, from (\ref{birkhoff for mu and mu1}) we obtain that $\mu=\mu_1$.

Using Corollary \ref{ergodic components of weak star} we prove that $y_0\in Y_{x_0}$ the same way as at the end of the proof of Lemma \ref{generic transverse Pesin manifold intersects all Gorbits}, which proves the assertion of our lemma.
\end{proof}

We can now prove unique ergodicity of the restricted action $P\acts M$ when $c(G)=2$. 

\begin{lemma}\label{unique P-ergodicity when rG=2}
     The restricted action $P\acts M$ is uniquely ergodic. 
\end{lemma}
\begin{proof}
    Let $\mu$ be a $P$-invariant and ergodic probability measure on $M$. By Lemma \ref{transverse pesin intersects all G-orbits when rG=2} we have that for $\mu$-a.e. $x$ the transverse Pesin manifold $\w^T_\mu(x)$ intersects all $G$-orbits of our given action $G\acts M$. Hence, because of the arguments given in the proof above of Lemma \ref{transverse pesin intersects all G-orbits when rG=2}, we obtain that all $P$-invariant and ergodic probability measures on $M$ are equal to $\mu$, which yields the assertion of our lemma.
\end{proof}

\begin{remark}\label{notation unique P-invariant cG=2}
    We denote by $\mu^P$ the unique $P$-invariant probability measure on $M$ with respect to the restricted action $P\acts M$ coming from our original action $G\acts M$, when $c(G)=2$. Note that the uniqueness of this $P$-invariant probability depends on our framework given by Assumption \ref{assumption no invariant probability}.  We also denote by $\w^T(x)$ the corresponding transverse Pesin manifolds with respect to $\mu^P$, for $\mu^P$-a.e. $x$. 
\end{remark}

\section{A geometric version of the High Entropy Method}

The main result of this section, which can be thought of as a geometric version of the High Entropy Method developed by Eisiendler and Katok, see \cite{ek} and \cite{el}, is the main tool to proving our main result in the case $r(c)=2$.

In this section, we consider a slightly more general setup than the one we have been working with so far. 

Let $G$ be a semisimple Lie group with finite center, $A$ a maximal torus of it, and $P$ a minimal parabolic of $G$ containing $A$. For the rest of this section, fix a locally free, $C^2$ action $G\acts M$, with $M$ a closed manifold for the rest of this section. In this section we do not make an assumption on the dimension of $M$. 
\subsubsection{Generic basins}

For any $a\in A$, consider the set
\begin{equation*}
\begin{split}
    \operatorname{GB}_M(a):=\bigcup\{\supp(\nu):\nu&\text{ is an $a$-ergodic and $\langle a\rangle\ltimes N_a$-invariant}\\
    &\text{probability measure on $M$}\},
\end{split}
\end{equation*}
and call it the \emph{generic basin of $a$ in $M$}.

\begin{remark}\label{supports of a-invariant measures generic basin}
If $\mu$ is an $a$-invariant probability measure on $M$, then $\supp(\mu)\subseteq\overline{\operatorname{GB}_M(a)}$, using the $a$-ergodic decomposition of $\mu$ and Birkhoff's ergodic theorem.    
\end{remark}

\subsubsection{Geometric high entropy method}

\begin{definition}
    Given an $A$-ergodic and $A$-invariant probability measure $\mu$ on $M$, consider the set $\Lambda_\mu^T$ of all Lyapunov functionals in $A^\ast$ with respect to $\mu$ such that their corresponding Osedelecs distribution are not tangent to the orbits of $G\acts M$.
\end{definition}

\begin{remark}
    The set $\Lambda_\mu^T$ is empty if and only if the action $G\acts M$ is homogenous, that is, if $M=G/\Gamma$ for some torsion-free cocompact lattice $\Gamma$ of $G$, and the action is the coset left action.
\end{remark}

Let $\alpha$ and $\beta$ be nonzero roots of $(G,A)$, such that $\alpha+\beta$ is a root as well.

\begin{proposition}[Geometric High Entropy Method]\label{geometric high entropy method}
    Let $\mu$ be a probability measure on $M$ such that for some $s\in\ker(\alpha+\beta)$ with $\beta(s)<0$, we have the following properties:
    \begin{enumerate}
        \item $\mu$ is $s$-ergodic and $\langle s\rangle\ltimes N_s$-invariant;
        \item all elements in $\Lambda^T_\mu$ evaluated at $s$ are negative, or $\Lambda^T_\mu(s)=\emptyset$;
        \item $\mu$ is $G^\alpha$-invariant, and $\alpha(s)\neq 0$.
    \end{enumerate}
       Then there exists $\Lambda\subseteq M$ of full $\mu$ measure for which the following holds:  if there exists $x\in \Lambda$ such that $ux\in\overline{\operatorname{GB}_M(s)}$ for some $u\in G^{\beta}\setminus\{e\}$, then $\mu$ is invariant under all elements of $[u,G^\alpha]$. 
\end{proposition}

We will prove this fact in Section \ref{section Proof of Geometric High Entropy Method}. 

\subsection{Almost everywhere continuity of Birkhoff averages}

For non-trivial $a\in A$, recall the definition of the subgroups $L_a$, $C_a$, $N_a$ and $Q_a$ of $G$, see (\ref{Na and La}).

The following fact yields a ``continuity property" for the class of measures appearing in Proposition \ref{geometric high entropy method}.

\begin{proposition}\label{almost continuity birkhoff averages}
     Let $\mu$ be an $a$-ergodic, $\langle a\rangle\ltimes N_a$-invariant probability measure on $M$, such that either $\Lambda^T_\mu(a)=\emptyset$, or, all elements in $\Lambda^T_\mu(a)$ are negative. Then for $\mu$-a.e. $x$, all $f\in C^0(M)$ and all $\epsilon>0$, there exists $\delta>0$ such that for every $y\in M$ with $d(x,y)<\delta$ the set of all $u\in N_a$ such that
     \begin{equation*}
         \left |\lim_{n\rightarrow\infty}\frac1n\sum_{k=0}^{n-1} f(a^kuy)-\int fd\mu\right|<\epsilon
     \end{equation*}
     has positive $\haar_{N_a}$ measure.
\end{proposition}

Let $Z_1$ be an $a$-invariant subset of $M$ for the $a$-ergodic measure $\mu$, satisfying Osedelets' Theorem, Theorem \ref{Osedelets theorem}, which is a nonuniformly partially hyperbolic set with respect to $a$, see Remark \ref{nonuniform partially hyperbolic from multiplicative}.

\begin{lemma}\label{choice unstable and centers}

    Let $\mu$ and $a\in A$ be as in Proposition \ref{almost continuity birkhoff averages}. We can take
\begin{equation}\label{obvious distros}
    E^u(x)=T_x(N_a\cdot x)\text{ and } E^0(x)=T_x(C_a\cdot x)
\end{equation}
for all $x\in Z_1$, as the unstable and center distributions, respectively, of the nonuniformly partially hyperbolic subset $Z_1$ of $M$ with respect to $a$.       
\end{lemma}

\begin{proof}

    The hypothesis on $\Lambda^T_\mu(a)$ in Proposition \ref{almost continuity birkhoff averages} implies that if $T_xM=\bigoplus_j E^{\lambda_j}(x)$ is the Osedelets' decomposition of $TM$ with respect to $Z_1$, for any $x\in Z_1$, then $\bigoplus_{\lambda_j<0} E^{\lambda_j}(x)$ is transversal to the $G$-orbit $G\cdot x$ at $x$, and it  contains $T_x(L_a\cdot x)$; thus, the distributions indicated in (\ref{obvious distros}) can be taken as the unstable and central distributions for $Z_1$ with respect to $a$.
\end{proof}

Consider the regular subsets of $Z_1$, $\{Z_1^\ell\}_{\ell>0}$, see . For $x$ in such regular sets $Z_1^\ell$, see Definition \ref{definition regular sets}, there exists a local stable manifold $\locstable(x;a)$, depending on this set, see Definition \ref{local stable manifolds}. 

\begin{remark}\label{Qa orbits are transversal to local stables}
    For any $x\in Z_1$, $Q_a\cdot x$ is transversal to $\locstable(x,a)$ at $x$ in $M$ with complementary dimension. Indeed, since $Q_a$ is the subgroup of $G$ generated by $C_a$ and $N_a$, this follows from Lemma \ref{choice unstable and centers}.
\end{remark}

\begin{lemma}\label{birkhoff averages small central perturbations}
For any $f\in C^0(M)$ and any $\epsilon>0$  there exists $\delta>0$ such that for any $x\in M$ and any $c\in C_a^{\leq\delta}$,
\begin{equation*}
    \limsup_{n\rightarrow\infty}\frac1n \sum_{j=0}^{n-1}|f(a^jcx)-f(a^jx)|<\epsilon.
\end{equation*}
\end{lemma}

\begin{proof}

    There exists $\delta_1>0$ such that $|f(x)-f(y)|<\epsilon$ whenever $d_M(x,y)\leq\delta$ for all $x,y\in M$, as $M$ is compact. Choose $\delta>0$ small enough so that $d(cx,x)<\delta_1$ for all $x\in M$ and all $c\in C_a^{\leq\delta}$, which can be done as $M$ is compact; so that 
    \begin{equation*}
    \limsup_{n\rightarrow\infty}\frac1n \sum_{j=0}^{n-1}|f(a^jcx)-f(a^jx)|=\limsup_{n\rightarrow\infty}\frac1n \sum_{j=0}^{n-1}|f(ca^jx)-f(a^jx)|<\epsilon.
    \end{equation*}
\end{proof}

If $x\in Z_1$ and $\ell>0$, recall there is a family of local stable manifolds, with respect to $Z_1^\ell$
\begin{equation}\label{hkhklg}
    \mathcal L^\ell_1(x)=\left \{\locstable(x;a):w\in Z^\ell_1\cap B^M(x_0,r_\ell)\right \},
\end{equation}
see (\ref{family of local stable manifolds}) in Section \ref{Absolute continuity along stable manifolds section}. 

Recall that for local transversals $T^1$ and $T^2$ to the family $\mathcal L^l_1(x)$, we have a holonomy map $\mathcal{Q}^\ell_1(x) \cap T^1 \to \mathcal{Q}^\ell_1(x) \cap T^2$ along the family $\mathcal L^\ell_1(x)$, where $\mathcal{Q}^\ell_1(x)$ is defined as in (\ref{union of family of local stables}) in Section \ref{Absolute continuity along stable manifolds section} with respect to $Z^\ell_1$, and the holonomy map is defined as in (\ref{holonomy between transversals along stable}) in the same section.


\begin{lemma}\label{nice property of holonomy map}    
Let $x\in Z_1$ be such that $ux\in Z_1$ for $\haar_{N_a}$-a.e. $u$. For any $r>0$ and $\ell>0$, there exists $\delta>0$ such that for any $y\in M$ with $d(x,y)<\delta$ we have: 
    \begin{enumerate}[label=(\roman*)]
        \item\label{fff1} $Q_a^{\leq r}\cdot y$ is a local transversal to $\mathcal L^\ell_1(x)$;
        \item\label{fff2} the corresponding holonomy map along $\mathcal L^\ell_1(x)$ from $\mathcal{Q}^\ell_1(x)\cap Q_a^{\leq r}\cdot x$ to $\mathcal{Q}^\ell_1(x)\cap Q_a^{\leq r}\cdot y$ restricts to an absolutely continuous map $$Z^\ell_1\cap Q_a^{\leq r}\cdot x\to \mathcal{Q}^\ell_1(x)\cap Q_a^{\leq r}\cdot y$$
        with respect to the volumes of the local submanifolds $Q_a^{\leq r}\cdot x$ and  $Q_a^{\leq r}\cdot y$.
    \end{enumerate}
\end{lemma}
\begin{proof}


    $Z_1^\ell$ is compact and  $\locstable(-;a)$ is  continuous, in the $C^1$-norm, on this set, see Definition \ref{local stable manifolds}; and thus, uniformly continuous. Hence, $Q_a^{\leq r}\cdot x$ is uniformly transversal to  the family (\ref{hkhklg}) by Remark \ref{Qa orbits are transversal to local stables}. This observation, along with using again the uniform continuity of $\locstable(-;a)$ on $Z_1^\ell$, implies that if $y$ is close enough to $x$, then $Q_a^{\leq r}\cdot y$ is uniformly transversal to the family (\ref{hkhklg}). This proves item \ref{fff1}. 

    
    Since $Q_a^{\leq r}\cdot x$ is uniformly transversal to the family (\ref{hkhklg}), then for any $y\in M$ close enough to $x$, the local manifold $Q_a^{\leq 2r}\cdot y$(containing the closure of $Q_a^{\leq r}\cdot y$) intersects all the elements of the family (\ref{hkhklg}) that intersect $Q_a^{\leq r}\cdot x$; because of how the holonomy map $$\mathcal{Q}^\ell_1(x) \cap  Q_a^{\leq r}\cdot x \to \mathcal{Q}^\ell_1(x) \cap Q_a^{\leq r}\cdot y$$ is defined in (\ref{holonomy between transversals along stable}) in Section \ref{Absolute continuity along stable manifolds section}, this implies that such holonomy map restricts to a map on $Z^\ell_1\cap Q_a^{\leq r}\cdot x$; since such holonomy maps are absolutely continuous by Theorem \ref{absolute continuity along stable}, this proves \ref{fff2}.
\end{proof}

\begin{proof almost continuity}
Suppose we chose the full $\mu$ measure, $a$-invariant, Osedelets set $Z_1$ so that the following two properties hold for every $x\in Z_1$ and $\haar_{N_a}$-a.e. $u$:
    \begin{enumerate}[label=(\roman*)]
        \item\label{eee1} $ux\in Z_1$;
        \item\label{eee2} Birkhoff's ergodic theorem holds for $ux$ with respect to $\mu$ and $a$; 
    \end{enumerate}
such $Z_1$ exists by applying Proposition \ref{haarconditionals} as $\mu$ is $N_a$-invariant, $a$-invariant and ergodic; and because of Remark \ref{nonuniform partially hyperbolic from multiplicative}. Fix any $x_0\in Z_1$. 

Choose $r>0$ such that $|f(gx)-f(x)|<\epsilon/2$ for all $x\in M$ and all $g\in G^{\leq r}$.

Shrink $r$ if necessary so that $Q_a^{\leq r}$ is a small enough neighborhood of $e_G$ in $Q_a$ for which Lemma \ref{multiplication maps} holds. Then from Fubini's theorem, Lemma \ref{birkhoff averages small central perturbations}, our choice of $r$, and Item \ref{eee2} above, we obtain that 
\begin{equation}\label{good birkhoff averages}
    \lim_{n\rightarrow\infty}\frac1n \sum_{j=0}^{n-1}f(a^jqx_0)\text{ exists and }\left|\lim_{n\rightarrow\infty}\frac1n \sum_{j=0}^{n-1}f(a^jqx_0)-\int fd\mu\right|<\epsilon/2.
\end{equation}
for $\haar_{Q_a}$-a.e. $q\in Q_a^{\leq r}$.

By Lemma \ref{nice property of holonomy map} we get that for any $y\in M$ with $d(x_0,y)<\delta$, there is a holonomy map $\phi:Z^\ell_1\cap Q_a^{\leq r}\cdot x_0\to\mathcal{Q}^\ell_1(x_0)\cap Q_a^{\leq r}\cdot y$ along the family (\ref{hkhklg}), satisfying item \ref{fff2} in Lemma \ref{nice property of holonomy map}

Let $Y$ denote the image under  $\phi$ of the set of all points $qx_0$ such that $q\in Q_a^{\leq r}$ satisfies (\ref{good birkhoff averages}). By construction of $\phi$ we have that for every $z\in Y$ there exists some $q\in Q_a^{\leq r}$  satisfying (\ref{good birkhoff averages}), along with some  $z\in\locstable(qx_0,a)$ such that
\begin{equation*}
\lim_{n\rightarrow\infty}\frac1n\sum_{j=0}^{n-1}f(a^jz)=\lim_{n\rightarrow\infty}\frac1n \sum_{j=0}^{n-1}f(a^jqx_0),
\end{equation*}
using that $f$ is uniformly continuous as $M$ is compact. Combining this last equality with (\ref{good birkhoff averages}) we obtain
\begin{equation}\label{fff3}
    \left|\lim_{n\rightarrow\infty}\frac1n \sum_{j=0}^{n-1}f(a^jz)-\int fd\mu\right|<\epsilon/2.
\end{equation}

However, since $\phi$ is absolute continuous, because of Item \ref{fff2} of Lemma \ref{nice property of holonomy map}, we have that $Y$ has positive volume in $Q_a^{\leq r}\cdot y$. Since we chose $r$ above so that Lemma \ref{multiplication maps} holds for $Q_a^{\leq r}$, we obtain from Fubini's Theorem that there exists $c_0\in C_a^{\leq r}$, and $\Lambda\subseteq N_a$ with positive $\haar_{N_a}$ measure such that $c_0\cdot\Lambda\cdot y\subseteq Y$. Therefore, because of Lemma \ref{birkhoff averages small central perturbations}, our choice of $r$, and (\ref{fff3}), we obtain
\begin{equation*}
    \left|\lim_{n\rightarrow\infty}\frac1n \sum_{j=0}^{n-1}f(a^juy)-\int fd\mu\right|<\epsilon
\end{equation*}
for all $u\in\Lambda$. This finishes the proof of our Lemma, as $Z_1$ has full $\mu$ measure.
\end{proof almost continuity}


We have the following consequence of Proposition \ref{almost continuity birkhoff averages}. Recall the definition of the group $L_a$ given in Section \ref{Some Lie group preliminaries}. 

\begin{lemma}\label{corollary almost continuity}
    Under the hypothesis of Proposition \ref{almost continuity birkhoff averages}, the following holds:
    for $\mu$-a.e. $x$, any $v\in L_a$, any $f\in C^0(M)$, and any $\epsilon>0$, there exists $\delta>0$ such that for any $y\in M$ with $d_M(vx,y)<\delta$, the set of all $u\in N_a$ such that
    \begin{equation*}
         \left |\lim_{n\rightarrow\infty}\frac1n\sum_{k=0}^{n-1} f(a^kuy)-\int fd\mu\right|<\epsilon
     \end{equation*}
     has positive $\haar_{N_a}$ measure.
\end{lemma}
\begin{proof}
    Choose $\delta_1>0$ for which Proposition \ref{almost continuity birkhoff averages} holds for $f$ and $\epsilon$.

    Because of Poincar\'e's recurrence theorem we have that for $\mu$-a.e. $x$ there exists a large enough $N>0$ such that
    \begin{equation*}
        d(a^Nvx,a^Nx), d(a^Nx,x)<\delta_1/3.
    \end{equation*}
    Choose $\delta>0$ such that $d(a^Ny,a^Nx)<\delta_1/3$ for any $y\in M$ with $d(x,y)<\delta$; for such $y$ we have $d(a^Ny,x)<\delta_1$, and thus, by our choices, this implies
    \begin{equation*}
        \left |\lim_{n\rightarrow\infty}\frac1n \sum_{j=0}^{n-1}f(a^jua^Ny)-\int fd\mu\right|<\epsilon.
    \end{equation*}
    for all $u\in X$, for some $X\subseteq N_a$ of positive $\haar_{N_a}$ measure. However, $a^{-N}Xa^N$ has positive $\haar_{N_a}$ measure as well, this implies from the previous equality that  
    \begin{equation*}
        \left |\lim_{n\rightarrow\infty}\frac1n \sum_{j=0}^{n-1}f(a^ja^Nuy)-\int fd\mu\right|<\epsilon,
    \end{equation*}
    for all $u$ in  a set of positive $\haar_{N_a}$ measure. It is easily checked that this implies the inequality of our Corollary.
\end{proof}

\begin{lemma}\label{corollary almost continuity along center}
    Let $\mu$ and $a\in A$ be as in Proposition \ref{almost continuity birkhoff averages}.  Let $x\in M$ be such that Lemma \ref{corollary almost continuity} holds for $x$ with respect to $\mu$ and $a$. If $w\in C_a$, then the conlusion of Lemma \ref{corollary almost continuity} is true for $wx$ with respect to $w_\ast\mu$ and $a$. 
\end{lemma}

\begin{proof}
    Using the root decomposition of $G$ with respect to $A$ we get that $C_a$ is a subgroup of both $N_G(L_a)$ and $N_G(N_a)$. In particular, elements of $C_a$ preserve both $L_a$ and $N_a$ orbits, this implies the assertion.
\end{proof}


\subsection{Proof of Geometric High Entropy Method}\label{section Proof of Geometric High Entropy Method}

\begin{lemma}\label{invariance pushforward along center}
    Let $a\in A$, and let $\alpha$ be a root of $(G,A)$ such that $\alpha(a)\neq 0$. Suppose $\nu$ is an $\langle a\rangle\ltimes G^\alpha$-invariant probability measure on $M$. Then for any $w\in C_a$, $w_\ast\nu$ is $\langle a\rangle\ltimes G^\alpha$-invariant as well. 
\end{lemma}
\begin{proof}
    By using the root decomposition of $G$ with respect to $A$ and that $\alpha(a)\neq 0$, we get that $C_a\leq N_G(G^\alpha)$, i.e., for any $w\in C_a$, $N_a\cdot w=w\cdot N_a$ as sets. Therefore, the statement now follows from Proposition \ref{haarconditionals}. 
\end{proof}

\begin{proof of geometric high entropy method}
Assume $x_0\in M$ is such that $u_0x_0\in \overline{\operatorname{GB}_M(s)}$, for some $u_0\in G^\beta\setminus\{e\}$. We indicate throughout the proof what $\mu$-generic properties $x_0$ has to satisfy.

We prove that $\mu$ is $[u_0,v]$-invariant for $\haar_{G^{\alpha}}$-a.e. $v$. This is enough to get the conclusion of Theorem \ref{geometric high entropy method}.

Suppose $x_0$ satisfies the following property:
    \begin{enumerate}[label=(\roman*)]
        \item\label{cvbmg} Lemma \ref{corollary almost continuity} holds for $vx_0$ with respect to $s$ and $\mu$, for $\haar_{G^\alpha}$-a.e. $v$, 
    \end{enumerate}
    the set of $x$ satisfying this property has full $\mu$ measure as $\mu$ is $G^\alpha$-invariant and by Proposition \ref{haarconditionals}.

    Set $w_0=[u_0,v]$.

    From property \ref{cvbmg}, since $\beta(s)<0$, and $u_0\in G^\beta$,  we get that for any $f\in C^0(M)$ and any $\epsilon>0$, there exists $\delta>0$ such that for any $y\in M$ with  $d_M(y,u_0x_0)<\delta$ the set of all $v\in N_a$ such that
    \begin{equation}\label{ygglhh}
         \left |\lim_{n\rightarrow\infty}\frac1n\sum_{k=0}^{n-1} f(s^kvy)-\int fd\mu\right|<\epsilon/2
     \end{equation}
     has positive $\haar_{N_s}$ measure. 
     
     Similarly, for any $y\in M$ with  $d_M(y,u_0w_0v_0^{-1}x_0)<\delta$, the set of all $v\in N_a$ such that
    \begin{equation}\label{ygglhh2}
         \left |\lim_{n\rightarrow\infty}\frac1n\sum_{k=0}^{n-1} f(s^kvy)-\int fd[(w_0)_\ast\mu]\right|<\epsilon/2
     \end{equation}
     has positive $\haar_{N_s}$ measure.

     As $u_0x_0\in\overline{\operatorname{GB}_M(s)}$, there exists $y_0\in M$ close enough to $u_0x_0$ so that $d(y_0,u_0x_0)<\delta$ and  $d(v^{-1}_0y_0,u_0w_0v^{-1}_0x_0)=d(v^{-1}_0y_0,v^{-1}_0u_0x_0)<\delta$(as $w_0=[u_0,v_0]$), which implies (\ref{ygglhh}) and (\ref{ygglhh2}) hold for $y_0$ and $v_0^{-1}y_0$, respectively; and for some $\langle a\rangle\ltimes N_a$-invariant and $a$-ergodic probability measure $\nu_0$, $y_0$ is generic enough for this measure so that 
     \begin{equation}\label{birkhoff for vy0}
     \lim_{n\rightarrow\infty}\frac1n\sum_{k=0}^{n-1}f(s^kvy_0)=\int fd\nu_0
     \end{equation}
    for $\haar_{N_s}$-a.e. $v$.

    Combining (\ref{ygglhh}), (\ref{ygglhh2}), for $y_0$ and $v_0^{-1}y_0$ respectively, along with (\ref{birkhoff for vy0}) yields
    \begin{equation*}
        \left |\int fd[(w_0)_\ast\mu]-\int f d\mu\right |<\epsilon.
    \end{equation*}

    Since $f$ and $\epsilon$ were arbitrary, this implies 
    \begin{equation*}
        (w_0)_\ast\mu=\mu.
    \end{equation*}

    Since $w_0=[u_0,v]$, and $v$ can vary in a subset of $G^\alpha$ of full $\haar_{G^\alpha}$ measure, Theorem \ref{geometric high entropy method} is proved.
\end{proof of geometric high entropy method}

\section{Proof of Theorem \ref{main} when  $c(G)=2$}\label{section main result non resonance}

We proved Theorem \ref{main} when $c(G)\geq 3$ in Section \ref{proof of main when not SL3R Section}.

Using the work of Section \ref{Proof of main result when G is not locally isomorphic to SL(3,R) Section}, we can prove Theorem \ref{main} when $c(G)=2$, when the unique $P$-invariant probability measure  $\mu^P$ on $M$, see Lemma \ref{unique P-ergodicity when rG=2} and Remark \ref{notation unique P-invariant cG=2}, satisfies the inequalities of Proposition \ref{Strict Maximal Inequality of transverse Lyapunov exponent}.

\begin{proposition}\label{no strict inequalities when G is SL3R}
   Suppose $G\acts M$ satisfies the hypothesis of Theorem \ref{main}, and that $c(G)=2$. If $G\acts M$ further satisfies Assumption \ref{assumption no invariant probability}, then the inequalities of Proposition \ref{Strict Maximal Inequality of transverse Lyapunov exponent} do not hold for $\mu^P$.
\end{proposition}

\begin{proof}

    Suppose, to get a contradiction, that the inequalities of Proposition \ref{Strict Maximal Inequality of transverse Lyapunov exponent} hold for $\mu^P$. 
    
    Then we are in the framework of Assumption \ref{assumption strict maximal}, and hence, in this case, we can repeat the argument we gave to prove Theorem \ref{main} in the $c(G)>2$ given in Section \ref{proof of main when not SL3R Section} for this case, and thus, we obtain a contradiction to Assumption \ref{assumption no invariant probability}
\end{proof}

Because of Proposition \ref{no strict inequalities when G is SL3R} and since we have proved Theorem \ref{main} for $c(G)=2$, to finish the prove of Theorem \ref{main}, and to ease notation as well, we work in this section, without loss of generality, under the following framework:

\begin{remark}\label{assumption final argument section}
    Because of Proposition \ref{no strict inequalities when G is SL3R}, to ease notation, from this point of the section on until the end of it, we work under Assumption \ref{combinatorial assumption on G rG=2 case}. 
\end{remark}

\subsection{Coherence property of Pesin manifolds}\label{Coherence property of Pesin manifolds section}

The first key tool towards proving Theorem \ref{main} in the $r(c)=2$ case, is the following ``coherence" property of $\mu^P$-generic Pesin manifolds.

\begin{proposition}\label{coherence proposition}
    For $\mu^P$-a.e. $x$, if $g\in G$ is such that $g\cdot \w^T(x)\cap \w^T(x)\neq\emptyset$, then $g\in Q_0=\operatorname{Stab}(\mu^P)$. 
\end{proposition}

We use Proposition \ref{geometric high entropy method} to prove this Proposition. We first need some technical observations.

\begin{lemma}\label{refined error recurrence}
  Suppose $Q\geq P$ is a parabolic subgroup of $G$. For $\mu^P$-a.e. $x$ we have that if $$g\cdot\w^T(x)\cap  \w^T(x)\neq\emptyset$$ for some $g\in G\setminus Q$, then there exist $\alpha\in\Delta_P$, with $-\alpha\notin\Sigma_Q$, such that for all  $p\in P$ in a set of positive $\haar_P$ measure, there is $u\in G^{-\alpha}\setminus\{e\}$ with $upx\in \supp(\mu^P)$.
\end{lemma}

\begin{proof}
    Suppose $g_0\in G\setminus Q$ and $x_0\in M$ are such that 
    \begin{equation}\label{translated self intersection}
        g_0\cdot\w^T(x_0)\cap  \w^T(x_0)\neq\emptyset.
    \end{equation}

    Fix $s_\alpha\in\operatorname{Int}(\ker(\alpha)\cap\overline{\w^P})$ for each $\alpha\in\Delta_P$.

    Suppose $x_0$ satisfies the following properties:

    \begin{enumerate}[label=(\roman*)]
        \item\label{rw1} $\mu^P$ disintegrates equivalently to Riemannian volume along $\w^T(x_0)$, which can be done because of Lemma \ref{maximal parabolic stabilizer when SMI fails};
        \item\label{rw1.5} $\w^T(x_0)\subseteq\supp(\mu^P)$, by applying Lemma \ref{maximal parabolic stabilizer when SMI fails} again;
        \item\label{rw2} for $[\mu^P]_{x_0}^{\w^T(x_0)}$-a.e. $y$ and $\haar_P$-a.e. $p$, $py$ satisfies Poincar\'e recurrence with respect to $s_\alpha$, for all $\alpha\in\Delta_P$.
    \end{enumerate}

    By (\ref{translated self intersection}), there exists $y_0\in\w^T(x_0)$ such that $g_0y_0\in\w^T(x_0)$. Because of \ref{rw1} and \ref{rw2} there is $y_1\in\w^T(x_0)$ close enough to $y_0$ satisfying \ref{rw2}, and for which there exists $g_1\in G\setminus Q$ such that $g_1y_1\in\w^T(x_0)$.

    By Lemma \ref{Pertubations by generic elements of $P$ have $LU$ factorization}, we have that for $\haar_P$-a.e. $p$, $pg_1=u'p'$ for some $u'\in U(P^{\operatorname{op}})$, $p'\in P$ is such that $p'y_1$ satisfies \ref{rw2}; and for some $\alpha\in\Delta_P$, $u'$ has non-trivial $G^{-\alpha}$ component, by applying Lemma \ref{almost surely nice perturbations}. Note that by Lemma \ref{Pertubations by generic elements of $P$ have $LU$ factorization} the set of all possible such $p'$ have positive $\haar_P$ measure.

    Since $p'y_1$ satisfies \ref{rw2}, there exists an increasing sequence $\{n_k\}_{k\in\N}\subseteq\N$ such that
    \begin{equation}\label{rw3}
        \lim_{k\rightarrow\infty}d_M(s_\alpha^{n_k}p'y_1,p'y_1)=0.
    \end{equation}

    However, since $u'\in U(P^{\operatorname{op}})$ has non-trivial $G^{-\alpha}$ component and $s_\alpha\in\operatorname{Int}(\ker(\alpha)\cap\overline{\w^P})$, there exists $v'\in G^{-\alpha}\setminus\{e\}$ such that 
    \begin{equation}\label{rw4}
        \lim_{n\rightarrow\infty}d_G(s_\alpha^{n}u's_\alpha^{-n},v')=0.
    \end{equation}
    We get that (\ref{rw3}) and (\ref{rw4}) together imply that 
   \begin{equation}\label{rw5}
       \lim_{k\rightarrow\infty}d_M(s_\alpha^{n_k}pg_1y_1,v'p'y_1)=\lim_{k\rightarrow\infty}d_M(s_\alpha^{n_k}u'p'y_1,v'p'y_1)=0.
   \end{equation}

   However, $g_1y_1\in\w^T(x_0)$, so that $g_1y_1\in\supp(\mu^P)$ by \ref{rw1.5} above, and thus $s_\alpha^npg_1y_1\in\supp(\mu^P)$, for any $n>0$ as $\mu^P$ is $P$-invariant. Hence, from (\ref{rw5}), it follows that $v'p'y_1\in\overline{\supp(\mu^P)}=\supp(\mu^P)$, and these points $p'$ constitute a subset of $P$ of positive $\haar_P$ measure. 
\end{proof}

Given $a\in A\setminus\{e\}$, consider the $a$-ergodic decomposition of $\mu^P$:
\begin{equation}\label{a ergodic decomposition of muP}
    \mu^P=\int\nu_x^a d\mu^P(x).
\end{equation}

Recall the defintion of the subgroup $N_a$ of $G$, see Section \ref{Stable and Unstable subgroups section}.

\begin{lemma}\label{enough from invariance of a component}
    Let $a\in A$ be such that $\lop(a)<0$ and $\mu^P$ is $N_a$-invariant. Let $\beta$ be a root in $(G,A)$ such that $\beta(a)=0$. For any $w\in G^{\beta}$ we have that  for $\mu^P$-a.e. $x$, if $\nu_x^a$ is invariant under the action of $w$, then $\mu^P$ is also invariant under $w$.
\end{lemma}

\begin{proof}

   Fix any $b\in \w^P$, in particular, $\mu^P$ is $b$-ergodic. 

   Fix $x_0$  such that Proposition \ref{almost continuity birkhoff averages} holds for $\nu_{x_0}^a$ with  respect to $a$ and $x_0$, the set of all such $x_0$'s has full $\mu^P$ measure. Assume $\nu_{x_0}^a$ is $w$-invariant, for some $w\in G^{\beta}$. Suppose also that Birkhoff's ergodic theorem holds for $\mu^P$ with respect to $x_0$ and $b$.

   Suppose $x_1\in M$ is such that $ux_1$ satisfies Birkhoff ergodic theorem for $a$ and $\nu_{x_1}^a$ for $\haar_{N_a}$-a.e. $u$. The set of all such $x_1$'s has full $\mu^P$ measure, as $\mu^P$ is $N_a$-invariant by the hypothesis on $a$, and by Lemma \ref{birkhoff for almost all px}, hence, in particular, we may assume $x_1\in\supp(\mu^P)$.
   
   We show, using the $w$-invariance of $\nu_{x_0}^a$, that for any such $x_1$, $\nu_{x_1}^a$ is $w$-invariant. This is enough to  prove our lemma. 

    Pick arbitrary $f\in C^0(M)$ and $\epsilon>0$.

    By Proposition \ref{almost continuity birkhoff averages} and Lemma \ref{corollary almost continuity along center} and since $b, w\in C_a$,  we get that there exists $\delta>0$ such that Proposition \ref{almost continuity birkhoff averages} holds for both $b^Nx_0$ and $wb^Nx_0$ with respect to $f$, $\epsilon$ and $a$, for the measures $b^N_\ast\nu_{x_0}^a$ and $(wb^N)_\ast\nu_{x_0}^a$, respectively. 

    Let $N>0$ be such that $b^Nx_0$ is close enough to $x_1$, which can be done as $x_1\in\supp(\mu^P)$ and our choice of $x_0$ above, so that both $d_M(b_0^Nx_0,x_1)<\epsilon$ and $d_M(wb^Nx_0,wx_1)<\epsilon$ hold. By Proposition \ref{almost continuity birkhoff averages}, this implies that there are $X_1,X_2\subseteq N_a$ both of positive $\haar_{N_a}$ measure such that 
    \begin{equation}\label{inequalities invariant component1}
          \left |\limsup_{n\to\infty}\frac1n\sum_{j=0}^{n-1}f(a^ju_1x_1)-\int f d[b^N_\ast\nu_{x_0}^a]\right |<\epsilon, \text{  and}
    \end{equation}
    \begin{equation}\label{inequalities invariant component2}
          \left |\limsup_{n\to\infty}\frac1n\sum_{j=0}^{n-1}f(a^ju_2wx_1)-\int f d[(wb^N)_\ast\nu_{x_0}^a]\right |<\epsilon,
    \end{equation}
    for all $u_1\in X_1$ and $u_2\in X_2$. 

    Since $\haar_{N_a}(X_1)>0$, our choice of $x_1$ and (\ref{inequalities invariant component1}) imply that
    \begin{equation*}
        \left |\int fd\nu_{x_1}^a-\int f d[b^N_\ast\nu_{x_0}^a]\right |<\epsilon.
    \end{equation*}

    From the definitions of Section \ref{section parabolic subgroups} it is not hard to see that $C_a\leq N_G(N_a)$. Thus, $w^{-1}X_2w$ has positive $\haar_{N_a}$ measure, so there exists $u_2\in X_2$ such that for $u_3:=w^{-1}u_2w$, Birkhoff's ergodic theorem holds for $u_3x_1$ with respect to $\nu_{x_1}^a$ and $a$. Since $a^ju_2w=a^jwu_3=wa^ju_3$, from (\ref{inequalities invariant component2}) we obtain
    \begin{equation*}
    \begin{split}
        \epsilon>&\left|\lim_{n\to\infty}\frac1n\sum_{j=0}^{n-1}f(wa^ju_3x_1)-\int f d[(wb^N)_\ast\nu_{x_0}^a]\right|\\
        =&\left |\int f d[w_\ast\nu_{x_1}^a]-\int f d[(wb^N)_\ast\nu_{x_0}^a]\right |.
    \end{split}
    \end{equation*}

    Set $w_1:=b^{-N}wb^N$, so that $w_1\in G^{\beta}$. Then $(wb^N)_\ast\nu_{x_0}^a=b^N_\ast(w_1)_\ast\nu_{x_0}^a=b^N_\ast\nu{x_0}^a$. Therefore, the last two inequalities imply that
    \begin{equation*}
        \left |\int fd\nu_{x_1}^a-\int f d[w_\ast\nu_{x_1}^a]\right |<2\epsilon. 
    \end{equation*}

    Therefore, $w_\ast\nu_{x_1}^a=\nu_{x_1}^a$, as $f\in C^0(M)$ and $\epsilon>0$ were arbitrary.

    Finally, since the set of points $x_1$ in $M$ as above has full $\mu^P$ measure, this proves $\mu^P$ is $w$-invariant.
\end{proof}

To begin with the proof of Proposition \ref{coherence proposition}, we need to fix some notation. 

Recall we are working in this section under Assumption \ref{combinatorial assumption on G rG=2 case}, see Remark \ref{assumption final argument section}. Hence, because of Lemma \ref{maximal parabolic stabilizer when SMI fails}, we have that $\lop=\beta_{2,1}+\beta_{3,1}$. 

Let 
\begin{equation*}
    \mu^P=\int\nu_x^{s_1}d\mu^P(x)
\end{equation*}
be the $s_1$-ergodic decomposition of $\mu^P$. 

Fix $s_1\in \ker(\beta_{1,3})$ such that $\beta_{2,1}(s_1)<0$, so that $\lop(s_1)<0$; we also have that $N_{s_1}=\langle G^{\beta_{1,2}},G^{\beta_{3,2}}\rangle$, so, because of Lemma \ref{maximal parabolic stabilizer when SMI fails}, we obtain that $\mu^P$ is $N_{s_1}$-invariant. 

\begin{remark}\label{conditions for geometric high entropy method}
    Take $\alpha=\beta_{3,2}$ and $\beta=\beta_{2,1}$, then the hypothesis of Proposition \ref{geometric high entropy method} hold for $s_1$ and for $\nu^{s_1}_x$ for $\mu^P$-a.e. $x$,  for this choice of $\alpha$ and $\beta$.
\end{remark}

\begin{proof of coherence proposition}
    Suppose, to get a contradiction, that there exists some $\Lambda\subseteq M$ of positive $\mu^P$ measure such that for all $x\in\Lambda$ there exists $g\in G\setminus Q_0$ such that 
    \begin{equation*}
        g\cdot\w^T(x)\cap\w^T(x)\neq\emptyset.
    \end{equation*}

    Fix $a\in\w^P$ and  $x_0\in\Lambda$ satisfying the following properties:
    \begin{enumerate}[label=(\roman*)]
        \item\label{xw1} $\mu^P$ disintegrates equivalently to Riemannian volume along $\w^T(x_0)$, which can be done because of Lemma \ref{maximal parabolic stabilizer when SMI fails};   
        \item\label{xw2} Lemma \ref{refined error recurrence} holds for $[\mu^P]_{x_0}^{\w^T(x_0)}$-a.e. $y$ with respect to $\w^T(x_0)$;
        \item\label{xw2.5} Lemma \ref{enough from invariance of a component} holds for $py$ for $[\mu^P]_{x_0}^{\w^T(x_0)}$-a.e. $y$ and $\haar_P$-a.e. $p$, which can be done as $\mu^P$ is $P$-invariant;
        \item\label{xw3} for $[\mu^P]_{x_0}^{\w^T(x_0)}$-a.e. $y$ and $\haar_P$-a.e. $p$, Theorem \ref{geometric high entropy method} holds for $\nu_{py}^{s_1}$ and $s_1$,  with $\alpha=\beta_{3,2}, \beta=\beta_{2,1}$.  
    \end{enumerate}

    Take $g_0\in G\setminus Q_0$ such that $g_0\cdot\w^T(x_0)\cap\w^T(x_0)\neq\emptyset$, so there is $y_0\in\w^T(x_0)$ such that $g_0y_0\in\w^T(x_0)$. Because of \ref{xw1} there exist $y\in\w^T(x_0)$ arbitrarily close to $y_0$ satisfying \ref{xw2}, \ref{xw2.5} and \ref{xw3}, and thus, we may choose such a $y_1$ very close to $y$ so as to have some $g_1\in G\setminus Q_0$ close to $g_0$ such that $g_1y_1\in\w^T(x_0)$.

    Since $g_1y_1\in\w^T(x_0)$, and $Q_0=\langle P,G^{\beta_{3,2}}\rangle$(because of Lemma \ref{maximal parabolic stabilizer when SMI fails}), we obtain from Lemma \ref{refined error recurrence} and \ref{xw2} that there exists a set of $p$'s in $P$ of positive measure for which there exists some $u\in G^{\beta_{2,1}}$ such that $upx\in\supp(\mu^P)$. Fix two such $u_0$ and $p_0$ such that $p_0y_1$ satisfies both \ref{xw2.5} and \ref{xw3}.
    
    However, by Remark \ref{supports of a-invariant measures generic basin}, we have that $\supp(\mu^P)\subseteq \overline{\operatorname{GB}_M(s_1)}$. Hence, because of \ref{xw3}, we obtain that $\nu^{s_1}_{p_0y_1}$ is invariant under the set $[u_0,G^{\beta_{3,2}}]$, but $[u_0,G^{\beta_{3,2}}]=G^{2,1}\setminus\{e\}$ as $G$ is locally isomorphic to $\SL(3,\R)$. Thus, $\nu^{s_1}_{p_0y_1}$ is $G^{2,1}$-invariant.

    From \ref{xw2.5} we obtain that $\mu^P$ is $G^{2,1}$-invariant, however, $\mu^P$ is $G^{\beta_{3,2}}$-invariant and $P$-invariant, which implies $\mu^P$ is $G$-invariant, which contradicts Assumption \ref{assumption strict maximal}, as we are working under Assumption \ref{assumption final argument section}.

    This contradiction proves Proposition \ref{coherence proposition}.
\end{proof of coherence proposition}

\subsection{Construction of equivariant maps $M\rightarrow G/H$}

We need the following differential topology observation.

\begin{proposition}\label{bundle from very nice transversal}
    Let $G$ be a Lie group, $M$ a closed connected manifold, and $G\acts M$ be a $C^1$ action. Let $H\leq G$ be a closed subgroup, and assume there exists a $C^1$, injectively immersed submanifold $T$ of $M$, satisfying the following properties:
    \begin{enumerate}[label=(\roman*)]
        \item\label{intersection all Gorbits} $T$ intersects all $G$-orbits, everywhere transversely, and these intersections are of complementary dimension in $M$;
        \item\label{self intersections of T only along H} for any $g\in G$, $gT\cap T\neq\emptyset$ only if $g\in H$;
    \end{enumerate}
    Then there exists a $C^1$ and $G$-equivariant map $M\rightarrow G/H$.
    \end{proposition}
\begin{proof}
    For any $x\in M$ there is $g\in G$ with $gx\in T$ because of \ref{intersection all Gorbits}; set $\phi(x):=g^{-1}H$. Then by \ref{self intersections of T only along H}, this map is well-defined, and yields a $G$-equivariant function $\phi:M\rightarrow G/H$. We only need to prove $\phi$ is smooth.

    Take any $x\in M$, and pick $g_x\in G$ with $\phi(x)=g_xH$. 

    The quotient map $G\xrightarrow{\pi}G/H$ is a fiber bundle, thus, there exists $X$ embedded submanifold of $G$ containing $g_x$, transversal everywhere to left cosets of $H$ in $G$, such that the multiplication map $X\times H\xrightarrow{m} G$ is a trivialization of $\pi$, and the map $X\xrightarrow{\sigma}G/H$, $g\mapsto gH$ is a diffeomorphism onto an open subset of $G/H$.

    By shrinking $X$ about $g_x$ if necessary, there exists $V\subset H$ open neighborhood of $e$ so that $U:=m(X\times V)$ is open in $G$.

    As $T$ intersects all $G$ orbits transversely, with complementary dimension, we get that for a small enough open subset $W$ of $T$, diffeomorphic to $\R^{\dim(T)}$, the map $\psi:U\times W\rightarrow M$, given by $(g,w)\mapsto g^{-1}w$, is a diffeomorphism onto an open neighborhood $B$ of $x$ in $M$. Then the diagram
     \[   \xymatrix{B \ar[r]^{\phi\upharpoonright} &G/H\\
U\times W \ar[u]_\psi &X\times V\times W\ar[u]_{\sigma\circ\operatorname{proj}_1}\ar[l]_{m\times\operatorname{Id}_W}}\]
commutes, where $\operatorname{proj}_1$ is the projection map of the product $X\times V\times W$ onto the first factor. Since all the maps in this diagram below $\phi\upharpoonright$ are all $C^1$, it follows that this map itself is $C^1$.

Thus, since $x\in M$ was arbitrary, we obtain that $\phi:M\rightarrow G/H$ is $C^1$.
\end{proof}

\subsection{Proof of Theorem \ref{main}}\label{proof of main theorem final section}

We can now finish the proof of Theorem \ref{main}. Since we have proved Theorem \ref{main} in the $r(G)>2$ case in Section \ref{proof of main when not SL3R Section}, we must obtain a contradiction under Assumption \ref{combinatorial assumption on G rG=2 case}, because of Proposition \ref{no strict inequalities when G is SL3R}, and thus, the argument below is given in the framework of Assumption \ref{combinatorial assumption on G rG=2 case}, see Remark \ref{assumption final argument section}.
\s

\begin{proof  main}
    By applying and combining Propositions \ref{coherence proposition} and \ref{bundle from very nice transversal}, along with  yields a smooth and $G$-equivariant map $\varphi:M\rightarrow G/Q_0$.

    From $G$-equivariance of $\varphi$, we get that it is a surjective submersion. Hence $e_GQ_0$ is a regular value of $\varphi$, and thus $N:=\varphi^{-1}(\{eQ_0\})$ is a closed embedded submanifold of $M$, of dimension $\dim(Q_0)+1$. Furthermore, from the $G$-equivariance of $\varphi$, it follows that $Q_0$ acts on $N$(this action is locally free since because our original action $G\acts M$ has this property). 
    
    In particular $P\acts N$, and hence there is a $P$-invariant (Radon) probability measure $\mu$ on $N$, and such a measure gives a (Radon) probability measure $\mu'$ on all of $M$ as $N$ is a closed subset of $M$ by defining for any Borel $X\acts M$, $\mu'(X)=\mu(X\cap N)$, and $\mu'$ is $P$-invariant as well, and thus $\mu'=\mu^P$ by unique ergodicity of the restricted action $P\acts M$, see Lemma \ref{unique P-ergodicity when rG=2}. Therefore, $\mu^P$ is supported on $N$.

    Hence, since $Q_0=\operatorname{Stab}(\mu^P)$, it follows that $Q_0\acts N$ has a $Q_0$-invariant probability measure. This implies that the associated codimension one orbit foliation $\mathcal F_{Q_0}$ of $Q_0\acts N$ has a transverse invariant measure by work of Plante, see \cite[Lemma~9.1]{plante}, see also the  proof of  \cite[Theorem~5.1]{sz}. By \cite[Theorem~6.3]{plante}, $\mathcal F_{Q_0}$ has a leaf with polynomial growth in $N$, hence, since $G/Q$ is compact, by using the $G$-equivariance of the map $\varphi:M\rightarrow G/Q$, we obtain that the foliation $\mathcal F_G$ by orbits of $G\acts M$ has a leaf with polynomial growth as well.

    The proof of  \cite[Theorem~5.1]{sz} shows, using the work of Stuck, see \cite[Theorem~C]{stuck}, that $G\acts M$ having an orbit with polynomial growth in $M$ implies that $G\acts M$ has an orbit closed in $M$, and in particular, this shows $G\acts M$ has a $G$-invariant probability measure, as such leaf is diffeormorphic to a manifold of the form $G/\Gamma$ for some $\Gamma\leq G$ uniform lattice. This contradicts Assumption \ref{assumption no invariant probability}, which in turn contradicts Assumption \ref{combinatorial assumption on G rG=2 case}. 

    By what we discussed at the beginning of this subsection, this shows that Assumption \ref{combinatorial assumption on G rG=2 case} does not hold, which proves Theorem \ref{main} in the $c(G)=2$ case, and in turn, because of what we did in Section \ref{proof of main when not SL3R Section}, it finishes the proof of Theorem \ref{main}.
\end{proof  main}

\section{Some consequences for codimension one actions of Lie groups}

\subsection{Codimension one actions of parabolic subgroups of higher rank simple Lie groups}

Using induced actions, see \cite[Section~2.3]{MR2807830}, we can prove that codimension one actions of parabolic subgroups of higher rank simple Lie groups do not exist, because of Theorem \ref{main}. Contrast this with the many smooth codimension one actions on closed manifolds of the parabolic subgroup of $\SL(2,\R)$, see \cite{ghys85} and \cite{asaokaannals}.

\begin{theorem}
    Let $Q\leq G$ be a proper parabolic subgroup, where $G$ is a higher rank simple Lie group, with finite center. There are no locally free smooth actions of $Q$ on closed manifolds with codimension one orbits.
\end{theorem}

\begin{proof}

    Suppose, to get a contradiction, such action $Q\acts N$ exists. Then we can induce it to an action $G\acts M$ that would be locally free, smooth, and with codimension one orbits. However, such induced actions have a smooth $G$-equivariant $M\rightarrow G/Q$. By Theorem \ref{main}, $G\acts M$ has a $G$-invariant probability measure, which pushes-forward to a $G$-invariant probability measure on $G/Q$, but this contradicts that $G\acts G/Q$ has no such measures, see, \cite[Example~9.1]{ZimmerMorris2008}. 
\end{proof}

\printbibliography

@article {MR1380646,
    AUTHOR = {Katok, A. and Lewis, J.},
     TITLE = {Global rigidity results for lattice actions on tori and new
              examples of volume-preserving actions},
   JOURNAL = {Israel J. Math.},
  FJOURNAL = {Israel Journal of Mathematics},
    VOLUME = {93},
      YEAR = {1996},
     PAGES = {253--280},
      ISSN = {0021-2172,1565-8511},
   MRCLASS = {22E40 (58F15)},
  MRNUMBER = {1380646},
MRREVIEWER = {Scot\ Adams},
       DOI = {10.1007/BF02761106},
       URL = {https://doi-org.proxy.lib.umich.edu/10.1007/BF02761106},
}

@article {MR4021068,
    AUTHOR = {Stover, Matthew},
     TITLE = {Lattices in {${\rm PU}(n,1)$} that are not profinitely rigid},
   JOURNAL = {Proc. Amer. Math. Soc.},
  FJOURNAL = {Proceedings of the American Mathematical Society},
    VOLUME = {147},
      YEAR = {2019},
    NUMBER = {12},
     PAGES = {5055--5062},
      ISSN = {0002-9939,1088-6826},
   MRCLASS = {20E18 (11F06 11G18 14G35 20F67 20H10)},
  MRNUMBER = {4021068},
MRREVIEWER = {B.\ Sury},
       DOI = {10.1090/proc/14763}
}

@article {MR2154667,
    AUTHOR = {Benveniste, E. Jerome and Fisher, David},
     TITLE = {Nonexistence of invariant rigid structures and invariant
              almost rigid structures},
   JOURNAL = {Comm. Anal. Geom.},
  FJOURNAL = {Communications in Analysis and Geometry},
    VOLUME = {13},
      YEAR = {2005},
    NUMBER = {1},
     PAGES = {89--111},
      ISSN = {1019-8385,1944-9992},
   MRCLASS = {53C24 (57S20)},
  MRNUMBER = {2154667},
MRREVIEWER = {Igor\ Belegradek},
       DOI = {10.4310/cag.2005.v13.n1.a2},
       URL = {https://doi-org.proxy.lib.umich.edu/10.4310/cag.2005.v13.n1.a2},
}

@misc{dh,
      title={Non left-orderability of lattices in higher rank semi-simple Lie groups}, 
      author={Bertrand Deroin and Sebastian Hurtado},
      year={2020},
      eprint={2008.10687},
      archivePrefix={arXiv},
      primaryClass={math.GR}
}

@book {bp,
    AUTHOR = {Barreira, Luis and Pesin, Yakov B.},
     TITLE = {Lyapunov exponents and smooth ergodic theory},
    SERIES = {University Lecture Series},
    VOLUME = {23},
 PUBLISHER = {American Mathematical Society, Providence, RI},
      YEAR = {2002},
     PAGES = {xii+151},
   MRCLASS = {37D25 (34D08 37C40 37D10)},
  MRNUMBER = {1862379},
MRREVIEWER = {Viviane Baladi},
       DOI = {10.1090/ulect/023},
}

@book{ZimmerMorris2008,
  author    = {Zimmer, Robert J. and Morris, Dave Witte},
  title     = {Ergodic Theory, Groups, and Geometry},
  series    = {CBMS Regional Conference Series in Mathematics},
  volume    = {109},
  publisher = {American Mathematical Society},
  address   = {Providence, RI},
  year      = {2008},
  doi       = {10.1090/cbms/109}
}

@article {bfh,
    AUTHOR = {Brown, Aaron and Fisher, David and Hurtado, Sebastian},
     TITLE = {Zimmer's conjecture: subexponential growth, measure rigidity,
              and strong property ({T})},
   JOURNAL = {Ann. of Math. (2)},
  FJOURNAL = {Annals of Mathematics. Second Series},
    VOLUME = {196},
      YEAR = {2022},
    NUMBER = {3},
     PAGES = {891--940},
       DOI = {10.4007/annals.2022.196.3.1},
}

@article {brhw,
    AUTHOR = {Brown, Aaron and Rodriguez Hertz, Federico and Wang, Zhiren},
     TITLE = {Invariant measures and measurable projective factors for
              actions of higher-rank lattices on manifolds},
   JOURNAL = {Ann. of Math. (2)},
  FJOURNAL = {Annals of Mathematics. Second Series},
    VOLUME = {196},
      YEAR = {2022},
    NUMBER = {3},
     PAGES = {941--981},
       DOI = {10.4007/annals.2022.196.3.2},      
}

@article {sz,
    AUTHOR = {Stuck, Garrett and Zimmer, Robert J.},
     TITLE = {Stabilizers for ergodic actions of higher rank semisimple
              groups},
   JOURNAL = {Ann. of Math. (2)},
  FJOURNAL = {Annals of Mathematics. Second Series},
    VOLUME = {139},
      YEAR = {1994},
    NUMBER = {3},
     PAGES = {723--747},
       DOI = {10.2307/2118577},
}

@article {nzann,
    AUTHOR = {Nevo, Amos and Zimmer, Robert J.},
     TITLE = {A structure theorem for actions of semisimple {L}ie groups},
   JOURNAL = {Ann. of Math. (2)},
  FJOURNAL = {Annals of Mathematics. Second Series},
    VOLUME = {156},
      YEAR = {2002},
    NUMBER = {2},
     PAGES = {565--594},
       DOI = {10.2307/3597198},
}

@article {brhwzd,
    AUTHOR = {Brown, Aaron and Rodriguez Hertz, Federico and Wang, Zhiren},
     TITLE = {Smooth ergodic theory of {$\mathbb Z^d$}-actions},
   JOURNAL = {J. Mod. Dyn.},
  FJOURNAL = {Journal of Modern Dynamics},
    VOLUME = {19},
      YEAR = {2023},
     PAGES = {455--540},
       DOI = {10.3934/jmd.2023014},
}

@article {MR422501,
    AUTHOR = {Millson, John J.},
     TITLE = {On the first {B}etti number of a constant negatively curved
              manifold},
   JOURNAL = {Ann. of Math. (2)},
  FJOURNAL = {Annals of Mathematics. Second Series},
    VOLUME = {104},
      YEAR = {1976},
    NUMBER = {2},
     PAGES = {235--247},
      ISSN = {0003-486X},
   MRCLASS = {22E40},
  MRNUMBER = {422501},
MRREVIEWER = {\`E.\ Vinberg},
       DOI = {10.2307/1971046},
}

@incollection {el,
    AUTHOR = {Einsiedler, M. and Lindenstrauss, E.},
     TITLE = {Diagonal actions on locally homogeneous spaces},
 BOOKTITLE = {Homogeneous flows, moduli spaces and arithmetic},
    SERIES = {Clay Math. Proc.},
    VOLUME = {10},
     PAGES = {155--241},
 PUBLISHER = {Amer. Math. Soc., Providence, RI},
      YEAR = {2010},
   MRCLASS = {22F30 (22D40 22F10 28D10 28D15 37A15 58J51)},
  MRNUMBER = {2648695},
MRREVIEWER = {Thomas Ward},
}

@article {plante,
    AUTHOR = {Plante, J. F.},
     TITLE = {Foliations with measure preserving holonomy},
   JOURNAL = {Ann. of Math. (2)},
  FJOURNAL = {Annals of Mathematics. Second Series},
    VOLUME = {102},
      YEAR = {1975},
    NUMBER = {2},
     PAGES = {327--361},
      ISSN = {0003-486X},
   MRCLASS = {57D30},
  MRNUMBER = {391125},
MRREVIEWER = {Robert Roussarie},
}

@book {bpintro,
    AUTHOR = {Barreira, Lu\'{\i}s and Pesin, Yakov},
     TITLE = {Introduction to smooth ergodic theory},
    SERIES = {Graduate Studies in Mathematics},
    VOLUME = {231},
      NOTE = {Second edition [of  3076414]},
 PUBLISHER = {American Mathematical Society, Providence, RI},
      YEAR = {2023},
     PAGES = {xv+336},
      ISBN = {9781470470654},
   MRCLASS = {37-01 (37Axx)},
  MRNUMBER = {4574839},
}

@incollection {ek,
    AUTHOR = {Einsiedler, Manfred and Katok, Anatole},
     TITLE = {Rigidity of measures---the high entropy case and non-commuting
              foliations},
      NOTE = {Probability in mathematics},
   JOURNAL = {Israel J. Math.},
  FJOURNAL = {Israel Journal of Mathematics},
    VOLUME = {148},
      YEAR = {2005},
     PAGES = {169--238},
      ISSN = {0021-2172},
   MRCLASS = {37C85 (37A15 37C40 37D30)},
  MRNUMBER = {2191228},
MRREVIEWER = {S. G. Dani},
       DOI = {10.1007/BF02775436},
}

@book {bm,
    AUTHOR = {Bekka, M. Bachir and Mayer, Matthias},
     TITLE = {Ergodic theory and topological dynamics of group actions on
              homogeneous spaces},
    SERIES = {London Mathematical Society Lecture Note Series},
    VOLUME = {269},
 PUBLISHER = {Cambridge University Press, Cambridge},
      YEAR = {2000},
     PAGES = {x+200},
      ISBN = {0-521-66030-0},
   MRCLASS = {37A15 (11E99 22E40 22F30 37A30 37C85)},
  MRNUMBER = {1781937},
MRREVIEWER = {S. G. Dani},
       DOI = {10.1017/CBO9780511758898},
}

@article {osedelec,
    AUTHOR = {Oseledec, V. I.},
     TITLE = {A multiplicative ergodic theorem. {C}haracteristic {L}japunov,
              exponents of dynamical systems},
   JOURNAL = {Trudy Moskov. Mat. Ob\v{s}\v{c}.},
  FJOURNAL = {Trudy Moskovskogo Matemati\v{c}eskogo Ob\v{s}\v{c}estva},
    VOLUME = {19},
      YEAR = {1968},
     PAGES = {179--210},
      ISSN = {0134-8663},
   MRCLASS = {28.70 (34.00)},
  MRNUMBER = {240280},
MRREVIEWER = {J\'{o}zsef Sz\H{u}cs},
}

@article {ledrappierstrelcyn,
    AUTHOR = {Ledrappier, Fran\c{c}ois and Strelcyn, Jean-Marie},
     TITLE = {A proof of the estimation from below in {P}esin's entropy
              formula},
   JOURNAL = {Ergodic Theory Dynam. Systems},
  FJOURNAL = {Ergodic Theory and Dynamical Systems},
    VOLUME = {2},
      YEAR = {1982},
    NUMBER = {2},
     PAGES = {203--219 (1983)},
      ISSN = {0143-3857},
   MRCLASS = {58F11 (28D20 34C35)},
  MRNUMBER = {693976},
MRREVIEWER = {M. I. Brin},
       DOI = {10.1017/S0143385700001528}
}

@article {ledrappieryoung2,
    AUTHOR = {Ledrappier, F. and Young, L.-S.},
     TITLE = {The metric entropy of diffeomorphisms. {II}. {R}elations
              between entropy, exponents and dimension},
   JOURNAL = {Ann. of Math. (2)},
  FJOURNAL = {Annals of Mathematics. Second Series},
    VOLUME = {122},
      YEAR = {1985},
    NUMBER = {3},
     PAGES = {540--574},
      ISSN = {0003-486X},
   MRCLASS = {58F11 (58F15)},
  MRNUMBER = {819557},
MRREVIEWER = {D. Newton},
       DOI = {10.2307/1971329}
}

@article {ledrappieryoung1,
    AUTHOR = {Ledrappier, F. and Young, L.-S.},
     TITLE = {The metric entropy of diffeomorphisms. {I}. {C}haracterization
              of measures satisfying {P}esin's entropy formula},
   JOURNAL = {Ann. of Math. (2)},
  FJOURNAL = {Annals of Mathematics. Second Series},
    VOLUME = {122},
      YEAR = {1985},
    NUMBER = {3},
     PAGES = {509--539},
      ISSN = {0003-486X},
   MRCLASS = {58F11 (58F15)},
  MRNUMBER = {819556},
MRREVIEWER = {D. Newton},
       DOI = {10.2307/1971328}
}

@book {wittemorrisbook,
    AUTHOR = {Morris, Dave Witte},
     TITLE = {Introduction to arithmetic groups},
 PUBLISHER = {Deductive Press},
      YEAR = {2015},
     PAGES = {xii+475},
      ISBN = {978-0-9865716-0-2; 978-0-9865716-1-9},
   MRCLASS = {22E40 (11F70 20G20 20G30 22D10 22D40)},
  MRNUMBER = {3307755},
}

@article {stuck,
    AUTHOR = {Stuck, Garrett},
     TITLE = {Growth of homogeneous spaces, density of discrete subgroups
              and {K}azhdan's property ({T})},
   JOURNAL = {Invent. Math.},
  FJOURNAL = {Inventiones Mathematicae},
    VOLUME = {109},
      YEAR = {1992},
    NUMBER = {3},
     PAGES = {505--517},
      ISSN = {0020-9910},
   MRCLASS = {22E40 (22E50 28D15 57S20)},
  MRNUMBER = {1176201},
MRREVIEWER = {Alain Valette},
       DOI = {10.1007/BF01232036},
}

@article {stuckminimal,
    AUTHOR = {Stuck, Garrett},
     TITLE = {Minimal actions of semisimple groups},
   JOURNAL = {Ergodic Theory Dynam. Systems},
  FJOURNAL = {Ergodic Theory and Dynamical Systems},
    VOLUME = {16},
      YEAR = {1996},
    NUMBER = {4},
     PAGES = {821--831},
      ISSN = {0143-3857},
   MRCLASS = {57S20 (22E46)},
  MRNUMBER = {1406436},
MRREVIEWER = {A. L. Onishchik},
       DOI = {10.1017/S0143385700009123},
}

@article {spatziervinhage,
    AUTHOR = {Spatzier, Ralf and Vinhage, Kurt},
     TITLE = {Cartan actions of higher rank abelian groups and their
              classification},
   JOURNAL = {J. Amer. Math. Soc.},
  FJOURNAL = {Journal of the American Mathematical Society},
    VOLUME = {37},
      YEAR = {2024},
    NUMBER = {3},
     PAGES = {731--859},
      ISSN = {0894-0347},
   MRCLASS = {37C85 (37C15 37C40 37D20)},
  MRNUMBER = {4736528},
MRREVIEWER = {Sanghoon Kwon},
       DOI = {10.1090/jams/1033},
}

@misc{brown2019entropylyapunovexponentsrigidity,
      title={Entropy, Lyapunov exponents, and rigidity of group actions}, 
      author={Aaron W. Brown and Sébastien Alvarez and Dominique Malicet and Davi Obata and Mario Roldán and Bruno Santiago and Michele Triestino},
      year={2019},
      eprint={1809.09192},
      archivePrefix={arXiv},
      primaryClass={math.DS},
      url={https://arxiv.org/abs/1809.09192}, 
}

@book {knapprep,
    AUTHOR = {Knapp, Anthony W.},
     TITLE = {Representation theory of semisimple groups},
    SERIES = {Princeton Mathematical Series},
    VOLUME = {36},
      NOTE = {An overview based on examples},
 PUBLISHER = {Princeton University Press, Princeton, NJ},
      YEAR = {1986},
     PAGES = {xviii+774},
      ISBN = {0-691-08401-7},
   MRCLASS = {22E46 (22-01 22E30)},
  MRNUMBER = {855239},
}

@misc{damjanovic2025zimmerprogrampartiallyhyperbolic,
      title={The Zimmer Program for partially hyperbolic actions}, 
      author={Danijela Damjanovic and Ralf Spatzier and Kurt Vinhage and Disheng Xu},
      year={2025},
      eprint={2211.08195},
      archivePrefix={arXiv},
      primaryClass={math.DS},
      url={https://arxiv.org/abs/2211.08195}, 
}

@book {zimmerbook,
    AUTHOR = {Zimmer, Robert J.},
     TITLE = {Ergodic theory and semisimple groups},
    SERIES = {Monographs in Mathematics},
    VOLUME = {81},
 PUBLISHER = {Birkh\"{a}user Verlag, Basel},
      YEAR = {1984},
     PAGES = {x+209},
       DOI = {10.1007/978-1-4684-9488-4},
}

@article {bfhinventiones,
    AUTHOR = {Brown, Aaron and Fisher, David and Hurtado, Sebastian},
     TITLE = {Zimmer's conjecture for actions of {${\rm SL}(m, \mathbb Z)$}},
   JOURNAL = {Invent. Math.},
  FJOURNAL = {Inventiones Mathematicae},
    VOLUME = {221},
      YEAR = {2020},
    NUMBER = {3},
     PAGES = {1001--1060},
      ISSN = {0020-9910},
   MRCLASS = {37C85 (20H25 57S20)},
  MRNUMBER = {4132960},
MRREVIEWER = {Rafael Oswaldo Ruggiero},
       DOI = {10.1007/s00222-020-00962-x},
}

@inproceedings {zimmericm,
    AUTHOR = {Zimmer, Robert J.},
     TITLE = {Actions of semisimple groups and discrete subgroups},
 BOOKTITLE = {Proceedings of the {I}nternational {C}ongress of
              {M}athematicians, {V}ol. 1, 2 ({B}erkeley, {C}alif., 1986)},
     PAGES = {1247--1258},
 PUBLISHER = {Amer. Math. Soc., Providence, RI},
      YEAR = {1987},
}

@article {chateletrosenberg,
    AUTHOR = {Chatelet, G. and Rosenberg, H.},
     TITLE = {Manifolds which admit {${\bf R}\sp{n}$} actions},
   JOURNAL = {Inst. Hautes \'Etudes Sci. Publ. Math.},
  FJOURNAL = {Institut des Hautes \'Etudes Scientifiques. Publications
              Math\'ematiques},
    NUMBER = {43},
      YEAR = {1974},
     PAGES = {245--260},
}

@article {rrw,
    AUTHOR = {Rosenberg, H. and Roussarie, R. and Weil, D.},
     TITLE = {A classification of closed orientable {$3$}-manifolds of rank
              two},
   JOURNAL = {Ann. of Math. (2)},
  FJOURNAL = {Annals of Mathematics. Second Series},
    VOLUME = {91},
      YEAR = {1970},
     PAGES = {449--464},
       DOI = {10.2307/1970633},
}

@article {ghys85,
    AUTHOR = {Ghys, \'{E}tienne},
     TITLE = {Actions localement libres du groupe affine},
   JOURNAL = {Invent. Math.},
  FJOURNAL = {Inventiones Mathematicae},
    VOLUME = {82},
      YEAR = {1985},
    NUMBER = {3},
     PAGES = {479--526},
       DOI = {10.1007/BF01388867},
}

@article {asaokaannals,
    AUTHOR = {Asaoka, Masayuki},
     TITLE = {Nonhomogeneous locally free actions of the affine group},
   JOURNAL = {Ann. of Math. (2)},
  FJOURNAL = {Annals of Mathematics. Second Series},
    VOLUME = {175},
      YEAR = {2012},
    NUMBER = {1},
     PAGES = {1--21},
       DOI = {10.4007/annals.2012.175.1.1},
}

@book {milne,
    AUTHOR = {Milne, J. S.},
     TITLE = {Algebraic groups},
    SERIES = {Cambridge Studies in Advanced Mathematics},
    VOLUME = {170},
      NOTE = {The theory of group schemes of finite type over a field},
 PUBLISHER = {Cambridge University Press, Cambridge},
      YEAR = {2017},
     PAGES = {xvi+644},
      ISBN = {978-1-107-16748-3},
   MRCLASS = {14L15 (14-01 17B45 20-01 20G15)},
  MRNUMBER = {3729270},
MRREVIEWER = {Boris\ \`E.\ Kunyavski\u i},
       DOI = {10.1017/9781316711736},
}

@incollection {MR2807830,
    AUTHOR = {Fisher, David},
     TITLE = {Groups acting on manifolds: around the {Z}immer program},
 BOOKTITLE = {Geometry, rigidity, and group actions},
    SERIES = {Chicago Lectures in Math.},
     PAGES = {72--157},
 PUBLISHER = {Univ. Chicago Press, Chicago, IL},
      YEAR = {2011},
      ISBN = {978-0-226-23788-6; 0-226-23788-5},
   MRCLASS = {22F05 (22E40 37C85 57S20)},
  MRNUMBER = {2807830},
MRREVIEWER = {Dave\ Witte\ Morris},
       DOI = {10.7208/chicago/9780226237909.001.0001},
}

@book {margulisbook,
    AUTHOR = {Margulis, G. A.},
     TITLE = {Discrete subgroups of semisimple {L}ie groups},
    SERIES = {Ergebnisse der Mathematik und ihrer Grenzgebiete (3) [Results
              in Mathematics and Related Areas (3)]},
    VOLUME = {17},
 PUBLISHER = {Springer-Verlag, Berlin},
      YEAR = {1991},
     PAGES = {x+388},
      ISBN = {3-540-12179-X},
}

@incollection {carriere1,
    AUTHOR = {Carri\`ere, Yves},
     TITLE = {Variations sur les flots riemanniens},
 BOOKTITLE = {S\'eminaire de {T}h\'eorie {S}pectrale et {G}\'eom\'etrie,
              {N}o.\ 5, {A}nn\'ee 1986--1987},
     PAGES = {43--55},
 PUBLISHER = {Univ. Grenoble I, Saint-Martin-d'H\`eres},
      YEAR = {1987},
   MRCLASS = {58F18 (53C12 57R30)},
  MRNUMBER = {1046241},
}

@misc{AIM2022,
  title        = {Global rigidity of actions by higher-rank groups},
  howpublished = {Workshop at the American Institute of Mathematics,
                  San Jose, California},
  year         = {2022},
  note         = {May 16--20, 2022. Problem session},
  shorthand = {AIM22}
}

@misc{morris2009latticesslnractcircle,
      title={Can lattices in SL(n,R) act on the circle?}, 
      author={Dave Witte Morris},
      year={2009},
      eprint={0811.0051},
      archivePrefix={arXiv},
      primaryClass={math.RT},
      url={https://arxiv.org/abs/0811.0051}, 
}

@article {fishermelnick,
    AUTHOR = {Fisher, David and Melnick, Karin},
     TITLE = {Smooth and analytic actions of {${\rm SL}(n,{\bf R})$} and
              {${\rm SL}(n,{\bf Z})$} on closed {$n$}-dimensional manifolds},
   JOURNAL = {Kyoto J. Math.},
  FJOURNAL = {Kyoto Journal of Mathematics},
    VOLUME = {64},
      YEAR = {2024},
    NUMBER = {4},
     PAGES = {873--904},
      ISSN = {2156-2261,2154-3321},
   MRCLASS = {22F50 (20F65 37C85 53C24 57S20)},
  MRNUMBER = {4794146},
MRREVIEWER = {Marja\ K.\ Kankaanrinta},
       DOI = {10.1215/21562261-2024-0009},
       URL = {https://doi.org/10.1215/21562261-2024-0009},
}

@misc{son2025realanalyticmathrmslnmathbbractionsclosed,
      title={Real analytic $\mathrm{SL}(n,\mathbb{R})$-actions on closed manifolds}, 
      author={Miri Son},
      year={2025},
      eprint={2508.04499},
      archivePrefix={arXiv},
      primaryClass={math.DS},
      url={https://arxiv.org/abs/2508.04499}, 
}

@misc{serrato2025globalrigiditycodimensionactions,
      title={Global Rigidity of Codimension One Actions}, 
      author={Arosemena Serrato, Camilo},
      year={2025},
      eprint={2512.00842},
      archivePrefix={arXiv},
      primaryClass={math.DS},
      url={https://arxiv.org/abs/2512.00842}, 
}

@article{BFHgafa,
  author    = {Brown, Aaron and Fisher, David and Hurtado, Sebastian},
  title     = {{Zimmer's} Conjecture for Non-uniform Lattices: Escape of Mass and Growth of Cocycles},
  journal   = {Geometric and Functional Analysis},
  year      = {2026},
  month     = nov,
  issn      = {1420-8970},
  doi       = {10.1007/s00039-026-00759-9},
  url       = {https://doi.org/10.1007/s00039-026-00759-9},
  publisher = {Springer}
}

@article {MR3104553,
    AUTHOR = {Agol, Ian},
     TITLE = {The virtual {H}aken conjecture},
      NOTE = {With an appendix by Agol, Daniel Groves, and Jason Manning},
   JOURNAL = {Doc. Math.},
  FJOURNAL = {Documenta Mathematica},
    VOLUME = {18},
      YEAR = {2013},
     PAGES = {1045--1087},
      ISSN = {1431-0635,1431-0643},
   MRCLASS = {20F67 (57Mxx)},
  MRNUMBER = {3104553},
MRREVIEWER = {Thomas\ Koberda},
       URL = {https://elibm.org/article/10000267},
}
\end{document}